\documentclass[amscd, leqno]{amsart}
\usepackage{amscd}
\usepackage{amsmath}
\usepackage{amssymb}
\usepackage{url}

\usepackage{color}
\definecolor{blue}{rgb}{0.0, 0.0, 1.0}
\definecolor{green}{rgb}{0.0, 0.5, 0.0}

\def\red#1{\textcolor{red}{#1}}

\theoremstyle{plain}
\newtheorem{theorem}{Theorem}[section]
\newtheorem{proposition}[theorem]{Proposition}
\newtheorem{lemma}[theorem]{Lemma}

\newtheorem{corollary}[theorem]{Corollary}
\theoremstyle{definition}

\theoremstyle{remark}
\newtheorem{remark}[theorem]{Remark}
\newtheorem{example}[theorem]{Example}

\newtheorem{conjecture}[theorem]{Conjecture}
\newtheorem{problem}[theorem]{Problem}
\newenvironment{pf}{\begin{proof}}{\end{proof}}
\makeatletter

\@addtoreset{equation}{section}
\makeatother

\begin{document}

	\title[Degeneration of Riemann surfaces and small eigenvalues]
	{Degeneration of Riemann surfaces and small eigenvalues of the Laplacian}

	\author{Xianzhe Dai}
	\address{
		Department of Mathematics, 
		University of Californai, Santa Barbara
		CA93106,
		USA}
	\email{dai@math.ucsb.jedu}
	\author{Ken-Ichi Yoshikawa}
	\address{
		Department of Mathematics,
		Faculty of Science,
		Kyoto University,
		Kyoto 606-8502,
		Japan}
	\email{yosikawa@math.kyoto-u.ac.jp}
	
	\begin{abstract}
		For a one-parameter degeneration of compact Riemann surfaces endowed with the K\"ahler metric induced from a K\"ahler metric on the total space of the family, we determine the exact magnitude of the small eigenvalues of the Laplacian as a function on the parameter space, 
		under the assumption that the singular fiber is reduced. 
		The novelty in our approach is that we compute the asymptotic behavior of certain difference of logarithm of analytic torsions in the degeneration in two ways. On the one hand, via heat kernel estimates, it is shown that the leading asymptotic is determined by the product of the small eigenvalues. On the other hand, using Quillen metrics, the leading asymptotic is connected with the period integrals, which we explicitly evaluate.
	\end{abstract}
	
	\maketitle
	\tableofcontents

	\section*{Introduction}
	\label{sect:0}
	Let $M$ be a compact Riemann surface of genus $g>1$ endowed with a Riemannian metric. 
	Let $C$ be \red{a} disjoint union of simple closed geodesics of $M$ such that
	$M\setminus C$ consists of $n+1$ components. Let ${\mathcal C}_{n}$ be the set of all those $C$. 
	For $C\in{\mathcal C}_{n}$, write $L(C)$ for the length of $C$. Set $\ell_{n} = \inf\{ L(C);\, C\in{\mathcal C}_{n} \}$.
	Let $0=\lambda_{0}<\lambda_{1}\leq\cdots\leq\lambda_{n}\leq$ be the eigenvalues of the Laplacian acting on the functions on $M$.
	Then the classical Schoen-Wolpert-Yau theorem \cite{SWY80} says the following:

	\begin{theorem}[Schoen-Wolpert-Yau]
		\label{thm:SWY}
		Let $k>0$ be a constant. Assume that the Gauss curvature $K$ satisfies $-1 \leq K \leq -k$.
		Then there exist positive constants $\alpha_{1},\alpha_{2}>0$ depending only on $g$ such that 
		$$
		\alpha_{1}k^{3/2} \ell_{n} \leq \lambda_{n} \leq \alpha_{2}\ell_{n}
		$$
		for $1\leq n \leq 2g-3$ and $\alpha_{1}k \leq \lambda_{2g-2} \leq \alpha_{2}$.
	\end{theorem}

	Furthermore, when $k=1$, namely, $M$ is a hyperbolic Riemann surface, Burger \cite{Burger90} proved that the small eigenvalues of $M$ are asymptotically calculated by those of the combinatorial Laplacian of a certain weighted graph associated to $M$ and the set of short geodesics of $M$.
	\par
	By Masur \cite{Masur76}, for a degeneration of compact Riemann surfaces to a stable curve, the hyperbolic metric on the fiber is comparable to 
	the hyperbolic metric on the annulus near the singular points. Namely, on a neighborhood of the vanishing cycles, the hyperbolic metric
	is bounded below and above by some constant multiple of the metric $dzd\bar{z}/(|z|^{2}(\log |z|)^{2})$ on the annulus.
	In particular, the magnitude of the length of any short geodesic is given by $1/\log (|s|^{-1})$, where the fiber is given locally by the equation $xy=s$
	near the nodes.
	From the Schoen-Wolpert-Yau theorem and Masur's theorem, for degenerations of compact Riemann surfaces to a stable curve, 
	it follows easily that there exist constants $C_{0}, C_{1}>0$ such that 
	\begin{equation}
		\label{eqn:small:ev:stable:c}
		\frac{C_{0}}{\log (|s|^{-1})} \leq \lambda_{1}(s) \leq\cdots\leq \lambda_{N-1}(s) \leq \frac{C_{1}}{\log (|s|^{-1})},
	\end{equation}
	where $N$ is the number of irreducible components of the singular fiber. (See Section~\ref{sect:1} for the details.)
	
	\medskip
	On the other hand, when the singularity of the singular fiber is more complicated and the K\"ahler metric of the fiber is no longer hyperbolic, 
	very little is known about the exact magnitude of the small eigenvalues of the Laplacian. 
	The goal of this article is to reveal the asymptotic behavior of the small eigenvalues of the Laplacian when the metric on the fibers 
	are induced from the K\"ahler metric on some ambient space. To state our results, let us introduce some notation and assumptions, which we keep throughout 
	this article.
	
	\medskip
	{\bf Set up }
	Let $f\colon X\to S$ be a proper surjective holomorphic map from a complex surface $X$ to a Riemann surface $S$ isomorphic to the unit disc of ${\mathbf C}$.
	We assume that $f$ has connected fibers and that $X_{0}$ is a unique singular fiber of $f$.
	Hence $\{0\} \subset S$ is the discriminant locus of $f$.
	We set $X_{s}=f^{-1}(s)$ for $s\in S$. We set $S^{o}=S\setminus\{0\}$, $X^{o}=X\setminus X_{0}$ and $f^{o}=f|_{X^{o}}$. Then $f^{o}\colon X^{o}\to S^{o}$
	is a family of compact Riemann surfaces. 
	Assume that $X$ carries a positive line bundle. In particular, $X$ is K\"ahler. Let $g^{X}$ be a K\"ahler metric on $X$. 
	We set $g_{s} = g^{X}|_{X_{s}}$. Then $(X_{s}, g_{s})$ $(s\not=0)$ is a compact Riemann surface endowed with a K\"ahler metric. We further make the following:
	
	\medskip
	\noindent
	{\bf Assumption }
	$X_{0}$ is a {\em reduced} and {\em reducible} divisor of $X$. In particular, $f$ has only isolated critical points on $X_{0}$.
	
	\medskip
	\par
	Let $0<\lambda_{1}(s)\leq\lambda_{2}(s)\leq\cdots$ be the eigenvalues of the Hodge-Kodaira Laplacian $\square_{s}=\bar{\partial}^{*}\bar{\partial}$ 
	counted with multiplicities, where $\square_{s}$ acts on the smooth functions on $X_{s}$ with respect to the induced metric $g_{s}$.
	For $s=0$, we regard $\square_{0}$ as the Friedrichs extension of the Laplacian acting on the smooth functions on 
	$X_{0,{\rm reg}} = X_{0}\setminus{\rm Sing}\,X_{0}$ with compact support. By Br\"uning-Lesch \cite{BruningLesch96}, 
	the spectrum of $\square_{0}$ consists of discrete eigenvalues and the heat operator of $\square_{0}$ is of trace class. 
	Moreover, $\ker \square_{0} \cong H^{0}(X_{0}\setminus{\rm Sing}\,X_{0}, {\mathbf C})$.
	For each $k\in{\mathbf N}$, the $k$-th eigenvalue $\lambda_{k}(s)$ is a continuous function on $S$ by Ji-Wentworth \cite{JiWentworth92} 
	when $X_{0}$ is a stable curve and by the second author \cite{Yoshikawa97} when $X_{0}$ is general. We set
	$$
	N = \dim H^{0}(X_{0}\setminus{\rm Sing}\,X_{0}, {\mathbf C}) = \#\{ \hbox{ irreducible components of } X_{0} \}. 
	$$
	By our assumption, $N>1$. From the continuity of $\lambda_{k}(s)$ as a function on $S$, it follows that
	$$
	\lim_{s\to0} \lambda_{k}(s) = 0
	\qquad
	(1 \leq k\leq N-1)
	$$
	and that $\lambda_{k}(s)$ is uniformly bounded from below by a positive constant for $k\geq N$. In this article, we investigate the asymptotic behavior of 
	the small eigenvalues $\lambda_{k}(s)$ for $1 \leq k\leq N-1$ as $s\to0$.
	\par
	In \cite{Gromov92}, Gromov gave an estimate for $\lambda_{1}(s)$ of the form
	\begin{equation}
		\label{eqn:Gromov}
		\lambda_{1}(s) \geq C\,|s|^{\alpha},
	\end{equation}
	where $C>0$ and $\alpha>0$ are constants. 
	It seems likely that a similar estimate can also be obtained by Cheeger's inequality \cite{Cheeger71}.
	By comparing \eqref{eqn:Gromov} with \eqref{eqn:small:ev:stable:c}, a natural question arises if the estimate \eqref{eqn:Gromov} is optimal or not.
	\par
	Since $X_{0}$ is not assumed to be a stable curve, there is no control of the critical points of $f$ except they consist of isolated points.
	In particular, any plane curve singularity can appear as a singularity of $X_{0}$ 
	as long as it is defined by a reduced equation. The following is the main result of this article.

	\begin{theorem}
		\label{Main:Thm:1}
		There exist constants $C_{0}, C_{1}>0$ such that for all $s \in S^{o}$, 
		$$
		\frac{C_{0}}{ \log (|s|^{-1}) } \leq \lambda_{1}(s) \leq \cdots \leq \lambda_{N-1}(s) \leq \frac{C_{1}}{ \log (|s|^{-1}) }.
		$$
	\end{theorem}

	This is in striking contrast to the rate of convergence of the small eigenvalues of Schr\"odinger operators
	when the central fiber is non-singular, which, restricted to any real analytic curve of $S$, is given by $|s|^{\nu}$ for some $\nu \in {\mathbf N}$ (cf. \cite{Farber95}).
	By Theorem~\ref{Main:Thm:1}, the estimate \eqref{eqn:small:ev:stable:c} obtained from the Schoen-Wolpert-Yau theorem and Masur's theorem
	holds true for general degenerations of compact Riemann surfaces if the singular fiber is reduced. 
	In fact, it is not difficult to prove the estimate $\lambda_{k}(s) \leq C/\log (|s|^{-1})$ for $1 \leq k \leq N-1$. 
	(See Section~\ref{sect:An upper bound of the small eigenvalues}.) Under this estimate, Theorem~\ref{Main:Thm:1} is deduced from the following
	(see Section~\ref{sect:7}):

	\begin{theorem}
		\label{Main:Thm:2}
		There exists a constant $c\in{\mathbf R}_{>0}$ such that as $s \to 0$,
		$$
		\prod_{k=1}^{N-1}  \lambda_{k}(s) = \frac{c + o(1)}{(\log (|s|^{-1}) )^{N-1}}.
		$$
		In particular, if $X_{0}$ consists of two irreducible components, then as $s \to 0$,
		$$
		\lambda_{1}(s) = \frac{c + o(1)}{ \log (|s|^{-1})}.
		$$
	\end{theorem}

	We remark that for degenerations of hyperbolic Riemann surfaces to stable curves, the corresponding result was obtained by 
	Grotowski-Huntley-Jorgenson \cite{GHJ01}\footnote{We are grateful to Professor Jorgenson for bringing this to our attention}. Also see Conjecture 9.8 for related discussion.
	
	\par
	Since the length $l(s)$ of any vanishing cycle of $(X_{s}, g_{s})$ is bounded above by $C|s|^{\nu}$ and from below by $C'|s|^{\nu'}$ for some constants
	$\nu,\nu',C,C'>0$, contrary to the Schoen-Wolpert-Yau theorem \cite{SWY80}, we conclude the following:

	\begin{corollary}
		As $s \to 0$, the small eigenvalue $\lambda_{k}(s)$ is comparable to $1/\log l(s)^{-1}$.
	\end{corollary}
	
	In \cite[Remark 5.10]{JiWentworth92}, Ji-Wentworth conjecture the second statement of Theorem~\ref{Main:Thm:2} with an explicit value of $c$, 
	when ${\rm Sing}\,X_{0}$ consists of a unique node. By Theorem~\ref{Main:Thm:2}, we have an affirmative answer to a generalization of their conjecture 
	without a comparison of the constant $c$ in Theorem~\ref{Main:Thm:2} with the constant in \cite{JiWentworth92}.

	Let us explain the strategy to prove Theorem~\ref{Main:Thm:2}.
	We choose a holomorphic line bundle $L$ on $X$ so that $L^{-1}$ is ample, and 
	a Hermitian metric $h$ on $L$ with semi-negative curvature such that $(L,h)$ is flat on a neighborhood of ${\rm Sing}\,X_{0}$.
	Let $\tau(X_{s},{\mathcal O}_{X_{s}})$ be the analytic torsion of the trivial Hermitian line bundle on $X_{s}$ and let $\tau(X_{s}, L_{s})$ be the analytic torsion of
	$(L, h)|_{X_{s}}$ (both defined using the metric $g_{s}$ on $X_{s}$).	
	We then compute the asymptotics in the degeneration of the difference $\log \tau(X_{s},{\mathcal O}_{X_{s}}) - \log \tau(X_{s},L_{s})$ in two different ways. On the one hand, using heat kernel estimates, we show that the leading asymptotic is given by the logarithm of the product of the small eigenvalues. On the other hand, we compute the asymptotics using the Quillen metrics and period integrals and show that the leading asymptotic is given by the logarithm of the right hand side of Theorem~\ref{Main:Thm:2}. 
	
	We emphasize that the curvature may well diverge to negative infinity in the degeneration (see Appendix). 
	Instead we rely crucially on the results of Li-Tian \cite{LiTian95}, Carlen-Kusuoka-Stroock \cite{CarlenKusuokaStroock87}, and Grigor'yan \cite{Grigor'yan97} for the uniform heat kernel upper bound. We make use of the partial analytic torsions introduced in \cite{DaiYoshikawa22} which localizes the analytic torsion in space and time.
	The fact that we are working with the difference of the analytic torsions also plays a critical role in dealing with the small time contribution near the singularity.
	More precisely, in Section~\ref{sect:3}, computing the behavior of the partial analytic torsions, we prove that as $s \to 0$,
	\begin{equation}
		\label{eqn:ratio:tor:1}
		\log\frac{\tau(X_{s},{\mathcal O}_{X_{s}})}{ \tau(X_{s}, L_{s})} = -\log \prod_{k=1}^{N-1} \lambda_{k}(s) + C + o(1).
	\end{equation}
	
	For the second way of calculating the asymptotics, we make critical use of the result of Bismut-Bost \cite{BismutBost90}, which gives the asymptotics of the Quillen metrics under degeneration. As the Quillen metric is the combination of the analytic torsion and the $L^2$-metric on the determinant of the cohomology, this leads to the period integrals, which can be computed using semi-stable reduction. It should be pointed out that the leading asymptotic arising in \cite{BismutBost90} gets cancelled out for $\log \tau(X_{s},{\mathcal O}_{X_{s}}) - \log \tau(X_{s},L_{s})$. We obtain our leading asymptotic, which is different, from the period integrals. Also, different metrics are needed in different steps, but that can be dealt with by using the anomaly formula of Bismut-Gillet-Soule \cite{BGS88}. 
	
	To explain in more detail, let $\widetilde{f} \colon Y \to T$ be a semi-stable reduction of $f \colon X \to S$ associated to a finite map $\mu \colon T \to S$. 
	Let $F \colon Y \to X$ be the corresponding map of total spaces sending $Y_{t} = \widetilde{f}^{-1}(t)$ to $X_{\mu(t)}$ for $t \in T \setminus\{0\}$. 
	Let $K_{Y/T}$ be the relative canonical bundle of $\widetilde{f}$. We assume that the line bundle $L^{-1}$ is defined by an effective divisor on $X$ whose support is disjoint from ${\rm Sing}\,X_{0}$ and intersects each irreducible component of $X_{0}$ transversally at a single point. (Such an assumption can always be satisfied in our situation since our family is over the unit disc with a reduced central fiber; Cf. the beginning of \S 4 for further clarification.)
	Then the direct image sheaves $\widetilde{f}_{*}K_{Y/T}$ and $\widetilde{f}_{*}K_{Y/T}(F^{*}L^{-1})$
	are locally free of rank $g$ and $g-1+N$, respectively, such that $\widetilde{f}_{*}K_{Y/T} \subset \widetilde{f}_{*}K_{Y/T}(F^{*}L^{-1})$.
	Let $\{\omega_{1},\ldots,\omega_{g+N-1}\}$ be a free basis of $\widetilde{f}_{*}K_{Y/T}(F^{*}L^{-1})$ around $0\in T$ such that
	$\{\omega_{1},\ldots,\omega_{g}\}$ is a free basis of $\widetilde{f}_{*}K_{Y/T}$. In Section~\ref{sect:4}, we prove that as $t \to 0$,
	\begin{equation}
		\label{eqn:ratio:tor:2}
		\log \frac{\tau(X_{\mu(t)},{\mathcal O}_{X_{\mu(t)}}) }{ \tau(X_{\mu(t)}, L_{\mu(t)})}
		=
		\log 
		\left[ \frac{ \det\left( \int_{Y_{t}} F^{*}h^{-1}(\omega_{i}(t) \wedge \overline{\omega_{j}(t)}) \right) }
		{\det\left(\int_{Y_{t}} \omega_{i}(t) \wedge \overline{\omega_{j}(t)} \right)}
		\right]
		+ C' + o(1),
	\end{equation}
	where $h^{-1}$ is the metric on $L^{-1}$ induced from $h$ and $C'$ is a constant. Also $1\leq i, j \leq g+N-1$ for the numerator and $1\leq i, j \leq g$ for the denominator; same for the equation below. Equations \eqref{eqn:ratio:tor:1}, \eqref{eqn:ratio:tor:2} yield the following key identity:
	\begin{equation}
		\label{eqn:McKean:Singer:mult}
		\prod_{k=1}^{N-1} \lambda_{k}( \mu(t) )^{-1} \equiv \frac{\tau(X_{\mu(t)},{\mathcal O}_{X_{\mu(t)}}) }{ \tau(X_{\mu(t)}, L_{\mu(t)})} \equiv
		\frac{ \det\left( \int_{Y_{t}} h_{F^{*}H}(\omega_{i}(t) \wedge \overline{\omega_{j}(t)}) \right) }
		{\det\left(\int_{Y_{t}} \omega_{i}(t) \wedge \overline{\omega_{j}(t)} \right)}.
	\end{equation}
	Here $f(t) \equiv g(t)$ if $f(t)/g(t)$ extends to a nowhere vanishing continuous function on $T$.
	By Hodge theory \cite{EFM21} and the theory of fiber integrals \cite{Barlet82}, \cite{Takayama21}, we prove that
	the right hand side of \eqref{eqn:McKean:Singer:mult} is of the form $(\log |t|^{-1})^{(N -1)}$, up to a nowhere vanishing continuous function on $T$ (Section~\ref{sect:5}), 
	which implies Theorem~\ref{Main:Thm:2}. 
	
	\par
	It is worth mentioning that, replacing the time parameter with the deformation parameter of the family $\widetilde{f} \colon Y \to T$,
	the role played by the ratio of the analytic torsions $\tau(X_{s},{\mathcal O}_{X_{s}}) / \tau(X_{s}, L_{s})$ in \eqref{eqn:McKean:Singer:mult} is similar to the one played by
	the difference of the heat traces in the McKean-Singer formula in the Atiyah-Singer index theorem in the sense that the ratio of analytic torsions provides a direct link 
	between the spectral quantity $\prod_{k=1}^{N-1} \lambda_{k}( \mu(t) )$ and the cohomological quantity, i.e., the ratio of the determinants of the period integrals.

	\par
	Recently, J. Cao \cite{Cao26} generalized Theorem~\ref{Main:Thm:1} by establishing the same result
	without assuming the reducedness and the relative dimension one of the fibers, using a different method of proof. We also note that in an earlier work \cite{Cao26.1}, Cao finds very interesting application of our result for the Archimedean pairing of parabolic dynamics.

	This article is organized as follows. In Section 1, we give a direct proof for Theorem~\ref{Main:Thm:1} for the case of semistable degeneration, using the Schoen-Wolpert-Yau theorem and Masur's theorem. Section 2 concerns the uniform heat kernel estimates. Then in Section 3, we compute the asymptotics of $\log \tau(X_{s},{\mathcal O}_{X_{s}}) - \log \tau(X_{s},L_{s})$ using the heat kernel estimates.  
	In Section 4, we recall semi-stable reductions and prove \eqref{eqn:ratio:tor:2}. Section 5 is devoted to the computation of the period integrals appearing in \eqref{eqn:ratio:tor:2} and we finally prove Theorem~\ref{Main:Thm:2}.
	Then, in Section 6, an upper bound is established for the small eigenvalues using the mini-max principle. This enables us to give the proof of Theorem~\ref{Main:Thm:1} in Section 7. In Section 8, we discuss some illustrating examples concerning small eigenvalues of Laplacian for degenerating families of Riemann surfaces. And finally, in Section 9, we end with some problems and conjectures. 
	In the appendix, we explain why the curvature diverges to negative infinity in our situation.

	\medskip
	
	{\bf Acknowledgements } 
	The authors thank the referee for their careful reading of the manuscript and helpful comments.
	The first author was partially supported by the Simons Foundation.
	The second author was partially supported by JSPS KAKENHI Grant Numbers 21H00984, 21H04429.

	\section{The small eigenvalues: semistable degeneration case}
	\label{sect:1}
	In this section, combining the Schoen-Wolpert-Yau theorem and Masur's theorem, we prove Theorem~\ref{Main:Thm:1} when $X_{0}$ is a stable curve of genus $g>1$. We maintain the notation from the introduction.

	\begin{lemma}
		\label{lemma:minimax}
		Let $M$ be a compact Riemann surface and let $g$, $g'$ be K\"ahler metrics on $M$. Let $\lambda_{1}$ (resp. $\lambda'_{1}$) be the first nonzero eigenvalue
		of the Laplacian of $(M,g)$ (resp. $(M,g')$). Then 
		$$
		\lambda'_{1} / \lambda_{1} \geq \min_{M} g / g'.
		$$
	\end{lemma}
	
	\begin{pf}
		Let $A^{0,1}=A^{0,1}(M)$ be the space of smooth $(0,1)$-forms and let ${\mathcal H}$ be the space of Abelian differentials on $M$. 
		Since the $\bar{\partial}$-operator induces an isomorphism from the eigenspace $E(\lambda; \square^{0,0})$ to $E(\lambda; \square^{0,1})$ for $\lambda>0$, 
		it follows from the mini-max principle that
		$$
		\begin{aligned}
			\lambda'_{1} 
			&= 
			\inf_{\phi \in A^{0,1}\cap\overline{\mathcal H}^{\perp}} \left| \frac{\int_{M} (\bar{\partial}\phi \otimes \partial\bar{\phi})/g'}{\int_{M}\phi\wedge\overline{\phi}} \right|
			=
			\inf_{\phi \in A^{0,1}\cap\overline{\mathcal H}^{\perp}} 
			\left| \frac{\int_{M} \{ (\bar{ \partial}\phi \otimes \partial\bar{\phi})/g \} (g/g')}{\int_{M}\phi\wedge\overline{\phi}} \right|
			\\
			&\geq
			\min_{M} g / g'
			\inf_{\phi \in A^{0,1}\cap\overline{\mathcal H}^{\perp}} 
			\left| \frac{\int_{M} (\bar{ \partial}\phi \otimes \partial\bar{\phi})/g }{\int_{M}\phi\wedge\overline{\phi}} \right|
			=
			(\min_{M} g / g') \cdot \lambda_{1}.
		\end{aligned}
		$$
		This completes the proof.
	\end{pf}

	\begin{theorem}
		Suppose that $X_{0}$ is a stable curve of genus $g>1$. Then there exist constants $C_{0}, C_{1}>0$ such that for all $s \in S^{o}$,
		$$
		\frac{C_{0}}{\log (|s|^{-1})} \leq \lambda_{1}(s) \leq \cdots \leq \lambda_{N-1}(s) \leq \frac{C_{1}}{\log (|s|^{-1})}.
		$$
	\end{theorem}
	
	\begin{pf}
		Let $p \in {\rm Sing}\,X_{0}$. Let $(U_{p},(z,w))$ be a coordinate neighborhood of $X$ centered at $p$ such that $f(z,w) = zw$ and $U_{p}=\{|z|<1, |w|<1\}$.
		Hence $X_{s}\cap U_{p} = \{ (z,w)\in\varDelta^{2};\, zw = s \}$ can be identified with the annulus $\{ z\in {\mathbf C};\, |s| < |z| <1 \}$.
		Let $g^{\rm hyp}_{s}$ be the hyperbolic metric on $X_{s}$. By Masur \cite{Masur76}, there exist constants $A_{1},A_{2}>0$ independent of $s \in S^{o}$
		such that for all $p \in {\rm Sing}\,X_{0}$,
		\begin{equation}
			\label{eqn:Masur:1}
			\frac{A_{1} dzd\bar{z}}{|z|^{2}(\log |z|)^{2}} \leq g^{\rm hyp}_{s}|_{X_{s}\cap U_{p}} \leq \frac{A_{2} dzd\bar{z}}{|z|^{2}(\log |z|)^{2}}
		\end{equation}
		and such that
		\begin{equation}
			\label{eqn:Masur:2}
			A_{1} g_{s} |_{X_{s}\setminus\bigcup_{p\in {\rm Sing}\,X_{0}} U_p}  \leq g^{\rm hyp}_{s} |_{X_{s}\setminus\bigcup_{p\in {\rm Sing}\,X_{0}} U_p}
			\leq
			A_{2} g_{s} |_{X_{s}\setminus\bigcup_{p\in {\rm Sing}\,X_{0}} U_p}.
		\end{equation}
		Let $\lambda^{\rm hyp}_{1}(s)$ be the first nonzero eigenvalue of the Laplacian of $(X_{s}, g^{\rm hyp}_{s})$. 
		Since there exists by \eqref{eqn:Masur:1}, \eqref{eqn:Masur:2} a constant $K>0$ with $g^{\rm hyp}_{s} \geq K g_{s}$ for all $s \in S^{o}$,
		it follows from Lemma~\ref{lemma:minimax} that
		\begin{equation}
			\label{eqn:comp:1st:ev}
			\lambda_{1}(s) \geq \min_{X_{s}} (g^{\rm hyp}_{s}/g_{s}) \cdot \lambda_{1}^{\rm hyp}(s) \geq K \lambda_{1}^{\rm hyp}(s)
			\qquad
			(s \in S^{o}).
		\end{equation}
		\par
		Write $\ell(s)$ for the length of the shortest simple geodesic of $X_{s}$. Then $\ell(s)$ is the $\ell_{1}$ for $(X_{s}, g^{\rm hyp}_{s})$. 
		By \eqref{eqn:Masur:1}, \eqref{eqn:Masur:2}, there exist constants $B_{1},B_{2}>0$  independent of $s \in S^{o}$ such that for all $s \in S^{o}$,
		\begin{equation}
			\label{eqn:short:geod}
			\frac{B_{1}}{\log (|s|^{-1})} \leq  \ell(s) \leq \frac{B_{2}}{\log (|s|^{-1})}.
		\end{equation}
		By \eqref{eqn:short:geod} and Theorem \ref{thm:SWY}, there exists a constant $C_{1}>0$ independent of $s \in S^{o}$ such that
		\begin{equation}
			\label{eqn:1st:ev:hyp}
			\lambda_{1}^{\rm hyp}(s) \geq \frac{C_{1}}{\log (|s|^{-1})}.
		\end{equation}
		By \eqref{eqn:comp:1st:ev}, \eqref{eqn:1st:ev:hyp}, there exists a constant $C_{2}>0$ independent of $s \in S^{o}$ such that
		\begin{equation}
			\label{eqn:1st:ev:ind}
			\lambda_{1}(s) \geq \frac{C_{2}}{\log (|s|^{-1})}.
		\end{equation}
		In Proposition~\ref{prop:lb:small:ev} below, we prove the existence of a constant $C_{3}>0$ with
		\begin{equation}
			\label{eqn:lb:ev}
			\lambda_{N-1}(s) \leq \frac{C_{3}}{\log (|s|^{-1})}
			\quad
			(s \in S^{o}).
		\end{equation}
		The result follows from \eqref{eqn:1st:ev:ind} and \eqref{eqn:lb:ev}.
	\end{pf}

	\section{Some estimates for the heat kernels}
	\label{sect:heat:kernel:estimates}
	In this section, we obtain some technical results concerning heat kernel estimates, which will play crucial roles to study the asymptotic behavior of 
	partial analytic torsions in the later section.
	\par
	Let $(M, g)$ be a compact K\"ahler manifold of complex dimension $n$. We assume that $M$ is projective. Namely, $M$ admits a holomorphic embedding
	into a projective space. Let $(L, h)$ be a holomorphic Hermitian line bundle on $M$. 
	Let $K^{L}(t,x,y)$ be the heat kernel of the Hodge-Kodaira Laplacian $\square^{L} =\bar{\partial}^{*}\bar{\partial}$ acting on the sections of $L$. 
	For $(x,y)\in M\times M$ and $t>0$, we have $K^{L}(t,x,y)\in{\rm Hom}(L_{y},L_{x})$.
	In what follows, the norm and inner product at each point are denoted by $|\cdot|$ and $\langle\cdot,\cdot\rangle$ respectively,
	while the $L^{p}$-norm and the $L^{2}$-inner product are denoted by $\| \cdot \|_{p}$ and $(\cdot, \cdot)$, respectively.

	\subsection{Gaussian type upper bounds}

	\begin{lemma}
		\label{lemma:2nd:derv:HK}
		Set $B:=\sup_{x\geq0}x^{2}e^{-x/2}$. Then for all $t>0$ and $x,y\in M$,
		$$
		| \square^{L}_{x} K^{L}(t,x,y) | \leq B^{\frac{1}{2}} t^{-1} \left\{ K^{L}( t/2,x,x) \right\}^{\frac{1}{2}}\left\{ K^{L}(t,y,y) \right\}^{\frac{1}{2}}.
		$$
	\end{lemma}
	
	\begin{pf}
		Let $\lambda_{1} \leq \lambda_{2} \leq \cdots$ be the eigenvalues of $\square^{L}$ counted with multiplicities. Let $\{ \phi_{i}(x) \}_{i\in{\mathbf N}}$ be a unitary
		basis of the Hilbert space of the $L^{2}$-sections of $L$ consisting of the eigenfunctions of $\square^{L}$ such that $\square^{L}\phi_{i} = \lambda_{i}\phi_{i}$.
		Since we have $K^{L}(t,x,y) = \sum_{i} e^{-t\lambda_{i}} \phi_{i}(x) \otimes \langle \cdot, \phi_{i}(y) \rangle_{y}$, we get 
		$K(t,x,x) = \sum_{i}e^{-t\lambda_{i}} |\phi_{i}(x)|^{2}$. By the Cauchy-Schwarz inequality and the definition of $B$, we get
		$$
		\begin{aligned}
			| \square^{L}_{x}K^{L}(t,x,y) | 
			&\leq 
			\sum_{i} \lambda_{i} e^{-\frac{t\lambda_{i}}{2}} |\phi_{i}(x)| \cdot e^{-\frac{t\lambda_{i}}{2}} |\phi_{i}(y)|
			\\
			&\leq
			\left\{ \sum_{i} \lambda_{i}^{2} e^{-t\lambda_{i}} |\phi_{i}(x)|^{2} \right\}^{\frac{1}{2}} \left\{ \sum_{j} e^{-t\lambda_{j}} |\phi_{j}(y)|^{2} \right\}^{\frac{1}{2}}
			\\
			&\leq
			B^{\frac{1}{2}}t^{-1} \left\{ \sum_{i} e^{-\frac{t\lambda_{i}}{2}} |\phi_{i}(x)|^{2} \right\}^{\frac{1}{2}} K^{L}(t,y,y)^{\frac{1}{2}}
			\\
			&=
			B^{\frac{1}{2}}t^{-1} K^{L}(t/2,x,x)^{\frac{1}{2}} K^{L}(t,y,y)^{\frac{1}{2}}.
		\end{aligned}
		$$
		This completes the proof.
	\end{pf}

	\par
	Let $\Omega$ and $\Omega'$ be domains of $M$ such that $\overline{\Omega} \subset \Omega'$.
	Let $\chi \in C^{\infty}_{0}(M)$ be a smooth function such that 
	$\chi\geq0$, $\chi=1$ on $\Omega$ and $\chi=0$ on $M\setminus \Omega'$.
	Let $A>0$ be a constant such that $|d\chi|_{g} \leq A$. (We can take $A=2/{\rm dist}(\partial\Omega', \partial\Omega)$.)
	Let $\nabla=\nabla^{L} = \partial^{L} + \bar{\partial}^{L}$ be the Chern connection of $(L,h^{L})$.
	Let $dv$ be the volume form of $(M,g)$, whose value at $z \in M$ is denoted by $dv_{z}$.

	\begin{lemma}
		\label{lemma:integral:nabla:HK}
		If $(L,h)|_{\Omega'}$ is a trivial Hermitian line bundle on $\Omega'$, then for all $y\in M$ and $0 < t \leq 1$, 
		the following inequality holds:
		$$
		\begin{aligned}
			\,&
			\int_{\Omega} |\nabla_{z}K^{L}(t,z,y)|^{2} dv_{z}
			\\
			&\leq
			4B^{\frac{1}{2}} t^{-1} \left\{ K^{L}(t,y,y)\right\}^{\frac{1}{2}}
			\left\{\int_{\Omega'} K^{L}(\frac{t}{2},z,z) dv_{z}\right\}^{\frac{1}{2}}
			\left\{ \int_{\Omega'} |K^{L}(t,z,y)|^{2} dv_{z} \right\}^{\frac{1}{2}}
			\\
			&\qquad+ 
			4A^{2} \int_{\Omega'} |K^{L}(t,z,y)|^{2} dv_{z}.
		\end{aligned}
		$$
	\end{lemma}
	
	\begin{pf}
		We proceed in the same way as in \cite[proof of Th.\,6]{ChenLiYau81}. Let $\sigma_{y}\in L_{y}$ be such that $|\sigma_{y}|=1$.
		Since $(L,h)|_{\Omega'}$ is a flat trivial line bundle on $\Omega'$ by assumption, there is a flat section of $L$ with length $1$
		defined on $\Omega'$. Trivializing $L$ with this flat section, we get $\nabla^{*}\nabla = 2\square$ on $\Omega'$, where $\nabla$
		is the Chern connection of $(L,h)$. Then 
		$$
		\begin{aligned}
			\,&
			\int_{\Omega'} \chi(z)^{2}|\nabla_{z}K^{L}(t,z,y)|^{2} dv_{z}
			=
			\int_{\Omega'} \langle \nabla_{z}^{*}\{\chi(z)^{2} \nabla_{z}K^{L}(t,z,y)\} \sigma_{y}, K^{L}(t,z,y)\sigma_{y}\rangle dv_{z}
			\\
			&\leq
			\int_{\Omega'} 
			\{ \chi(z)^{2} |\nabla_{z}^{*}\nabla_{z}K^{L}(t,z,y)| + 2 \chi(z)|d\chi(z)| |\nabla_{z}K^{L}(t,z,y)| \} |K^{L}(t,z,y)| dv_{z}
			\\
			&\leq
			\int_{\Omega'} 
			\chi(z)^{2} |2\square_{z} K^{L}(t,z,y)| |K^{L}(t,z,y)| dv_{z} 
			\\
			&\quad+ 
			\frac{1}{2}\int_{\Omega'} \chi(z)^{2} |\nabla_{z}K^{L}(t,z,y)|^{2} dv_{z}
			+
			2\int_{\Omega'} |d\chi(z)|^{2} |K^{L}(t,z,y)|^{2} dv_{z}
			\\
			&\leq
			\left\{ \int_{\Omega'} |2\square_{z} K^{L}(t,z,y)|^{2} dv_{z} \right\}^{\frac{1}{2}}
			\left\{ \int_{\Omega'} |K^{L}(t,z,y)|^{2} dv_{z} \right\}^{\frac{1}{2}}
			\\
			&\quad+ 
			\frac{1}{2}\int_{\Omega'} \chi(z)^{2} |\nabla_{z}K^{L}(t,z,y)|^{2} dv_{z}
			+
			2A^{2} \int_{\Omega'} |K^{L}(t,z,y)|^{2} dv_{z}.
		\end{aligned}
		$$
		Hence we get
		\begin{equation}
			\label{eqn:gradient:HK}
			\begin{aligned}
				\,&
				\int_{\Omega} |\nabla_{z}K^{L}(t,z,y)|^{2} dv_{z}
				\leq
				\int_{\Omega'} \chi(z)^{2}|\nabla_{z}K^{L}(t,z,y)|^{2} dv_{z}
				\\
				&\leq
				2\left\{ \int_{\Omega'} |2\square_{z} K^{L}(t,z,y)|^{2} dv_{z} \right\}^{\frac{1}{2}}
				\left\{ \int_{\Omega'} K^{L}(t,z,y)^{2} dv_{z} \right\}^{\frac{1}{2}}
				+ 
				4A^{2} \int_{\Omega'} K^{L}(t,z,y)^{2} dv_{z}.
			\end{aligned}
		\end{equation}
		Substituting the inequality in Lemma~\ref{lemma:2nd:derv:HK} into \eqref{eqn:gradient:HK}, we get 
		\begin{equation}
			\label{eqngradient:HK:2:}
			\begin{aligned}
				\,&
				\int_{\Omega} |\nabla_{z}K^{L}(t,z,y)|^{2} dv_{z}
				\leq
				\\
				&
				\frac{4B^{\frac{1}{2}}}{t} \sqrt{K^{L}(t,y,y)}
				\left\{ \int_{\Omega'} K^{L}(t/2,z,z) dv_{z} \right\}^{\frac{1}{2}}
				\left\{ \int_{\Omega'} |K^{L}(t,z,y)|^{2} dv_{z} \right\}^{\frac{1}{2}}
				\\
				&+ 
				4A^{2} \int_{\Omega'} |K^{L}(t,z,y)|^{2} dv_{z}.
			\end{aligned}
		\end{equation}
		This completes the proof.
	\end{pf}

	Recall that $n$ is the complex dimension of $M$. Hence $2n$ is the real dimension of $M$. Since $M$ is projective by assumption,
	we have a projective embedding $\iota \colon M \hookrightarrow {\mathbb P}^{N}$. Let $g^{M}_{\rm FS}$ be the restriction of the Fubini-Study metric on ${\mathbb P}^{N}$ to $M$ via $\iota$. Then there exists a constant $\gamma>0$ such that
	\begin{equation}
		\label{eqn:quasi:isom}
		\gamma^{-1} \, g^{M}_{\rm FS} \leq g^{M} \leq \gamma \, g^{M}_{\rm FS}.
	\end{equation}
	Let $d(\cdot, \cdot) = d_{M}(\cdot, \cdot)$ be the distance function on $M$ with respect to the metric $g^{M}$.
	In the next lemma, we prove a uniform Gaussian type upper bound for the heat kernel of $(M,g^{M})$, generalizing an analogous estimate by Li-Tian \cite{LiTian95} for the Fubini-Study metric.

	\begin{lemma}
		\label{lemma:Gauss:sc:h:k}
		Let $k(t,x,y)$ be the heat kernel of $(M, g^{M})$ acting on the functions on $M$.
		Then there exists a constant $C_{1}=C_{1}(n, \gamma)$ depending only on $n$ and $\gamma$ such that
		$$
		k(t,x,y) \leq C_{1}(n,\gamma)t^{-n} e^{-\frac{d(x,y)^{2}}{8t}}
		$$
		for all $t \in (0,1]$ and $x,y \in M$.
	\end{lemma}
	
	\begin{pf}
		Let $k_{\rm FS}^{M}(t,x,y)$ be the heat kernel of $(M, g^{M}_{\rm FS})$ and let $r(x,y) = d_{{\mathbb P}^{N}}(\iota(x), \iota(y))$ be the distance of the two points $\iota(x), \iota(y) \in {\mathbb P}^{N}$ with respect to the Fubini-Study metric on ${\mathbb P}^{N}$. 
		By Li-Tian \cite[Main Result]{LiTian95}, we have the following Gaussian type upper bound
		$$
		k_{\rm FS}^{M}(t,x,y) \leq C(n)t^{-n} e^{-\frac{r(x,y)^{2}}{8t}}
		$$
		for all $t\in (0,1]$ and $x,y \in M$, where $C(n)>0$ is a constant depending only on $n$. 
		By Carlen-Kusuoka-Stroock \cite[Th.\,2.1]{CarlenKusuokaStroock87}, this inequality implies the Nash inequality for $(M, g^{M}_{\rm FS})$:
		\begin{equation}
			\label{eqn:Nash:1}
			\| f \|_{2,{\rm FS}}^{2+\frac{2}{n}} \leq \alpha(n) ( \| df \|_{2,{\rm FS}}^{2} + \| f \|_{2,{\rm FS}}^{2} ) \cdot \| f \|_{1,{\rm FS}}^{\frac{2}{n}},
			\qquad
			f \in C^{\infty}(M),
		\end{equation}
		where $\alpha(n) > 0$ is a constant depending only on $n$. Here $\| \cdot \|_{p,{\rm FS}}$ denotes the $L^{p}$-norm with respect to $g^{M}_{\rm FS}$.
		By \eqref{eqn:quasi:isom}, \eqref{eqn:Nash:1}, there exists a constant $\alpha(n,\gamma)>0$ depending only on $n$ and $\gamma$ such that
		\begin{equation}
			\label{eqn:Nash:2}
			\| f \|_{2}^{2+\frac{2}{n}} \leq \alpha(n,\gamma) ( \| df \|_{2}^{2} + \| f \|_{2}^{2} ) \cdot \| f \|_{1}^{\frac{2}{n}},
			\qquad
			f \in C^{\infty}(M),
		\end{equation}
		where all the norms are those with respect to $g^{M}$. Then again by \cite[Th.\,2.1]{CarlenKusuokaStroock87}, we have the following upper bound
		for all $t \in (0,1]$ and $x \in M$
		\begin{equation}
			\label{eqn:weak:Gauss:bd}
			k(t,x,x) \leq C_{0}(n,\gamma) t^{-n},
		\end{equation}
		where $C_{0}(n,\gamma)>0$ is a constant depending only on $n$ and $\gamma$. Since $k(t,x,x)$ is decreasing in $t$, 
		we deduce from \eqref{eqn:weak:Gauss:bd} that $k(t,x,x) \leq C_{0}(n,\gamma)(t^{-n} + 1)$ for all $t>0$ and $x \in M$.
		By Grigor'yan \cite[Th.\,1.1]{Grigor'yan97}, this implies the desired Gaussian type upper bound.
	\end{pf}

	Let $c_{1}(L,h^{L})$ be the Chern form of $(L,h^{L})$ and let $\Lambda$ be the adjoint of the multiplication of the K\"ahler form of $M$.
	By the Bochner-Kodaira-Nakano formula \cite[Chap.\,VII, Cor.\,1.3]{Demailly12}, we have the following identity of differential operators on $A^{0}_{M}(L)$:
	$$
	\begin{aligned}
		(\nabla^{L})^{*}\nabla^{L} 
		&= 
		(\partial^{L}+\bar{\partial}^{L})^{*}(\partial^{L}+\bar{\partial}^{L})
		=
		(\partial^{L})^{*}\partial^{L} + (\bar{\partial}^{L})^{*}\bar{\partial}^{L}
		\\
		&=
		2\square^{L} + 2\pi \Lambda_{g^{M}} c_{1}(L,h^{L}),
	\end{aligned}
	$$
	where $\Lambda_{g^{M}}$ is the Lefschetz operator with respect to $g^{M}$. We define 
	$$
	Q^{L}:=2\pi \Lambda_{g^{M}} c_{1}(L,h^{L}) \in C^{\infty}(M).
	$$ 
	Hence 
	$$
	(\nabla^{L})^{*}\nabla^{L} = 2\square^{L} + Q^{L}.
	$$
	We set 
	$$
	\kappa = \kappa^{L} := \sup_{x\in M} |Q^{L}(x)|.
	$$
	When $L={\mathcal O}_{M}$ and $h=|\cdot|$ is the trivial metric on ${\mathcal O}_{M}$, we have $K^{L}(t,x,y) = k(t,x,y)$.

	\begin{lemma}
		\label{lemma:Gaussian:nabla:HK}
		Suppose that $(L,h)|_{\Omega'}$ is a trivial Hermitian line bundle on $\Omega'$. Then
		the following inequalities hold:
		\begin{itemize}
			\item[(1)]
			For all $x,x'\in M$ and $t\in(0,1]$,
			$$
			|K^{L}(t,x,x')| \leq C_{1}(n,\gamma)\,e^{\kappa} t^{-n}e^{-\frac{d(x,x')^{2}}{8t}}.
			$$
			\item[(2)]
			For all $y\in M\setminus\Omega'$ and $t\in(0,1]$,
			$$
			\int_{\Omega} |\nabla_{z}K^{L}(t,z,y)|^{2} dv_{z}
			\leq
			C_{2}(n,\gamma,A)\, e^{2\kappa}{\rm Vol}(M)t^{-(2n+1)}\,e^{-\frac{d(y, \Omega')^{2}}{8t}},
			$$
			where $C_{2}(n,\gamma,A) := 4C_{1}(n,\gamma)^{2}(2^{\frac{n}{2}}B^{\frac{1}{2}}+A^{2})$ and
			$A$ is the same constant as in Lemma~\ref{lemma:integral:nabla:HK}. 
		\end{itemize}
	\end{lemma}
	
	\begin{pf}
		By \cite[p.32 l.4-l.5]{HessSchraderUhlenbrock80}, the following inequality holds for all $t >0$ and $x,x' \in M$:
		$$
		|K^{L}(t,x,x')| \leq e^{\kappa t} k(t,x,x'). 
		$$
		This, together with Lemma~\ref{lemma:Gauss:sc:h:k}, yields (1). Write $C_{1}$ for $C_{1}(n, \gamma)$.
		If $z\in\Omega'$, then $d(z,y)\geq d(y,\Omega')$. Hence
		$|K^{L}(t,z,y)| \leq C_{1}(n, \gamma) e^{\kappa}t^{-n} e^{-\frac{d(y,\Omega')^{2}}{8t}}$ by (1).
		In particular, $|K^{L}(t,x,x)| \leq C_{1}e^{\kappa}t^{-n}$ for all $t\in(0,1]$ and $x\in M$.
		Substituting these inequalities into the inequality in Lemma~\ref{lemma:integral:nabla:HK}, we get
		$$
		\begin{aligned}
			\,&
			\int_{\Omega} |\nabla_{z}K^{L}(t,z,y)|^{2} dv_{z}
			\leq
			\\
			&
			\frac{4B^{\frac{1}{2}}}{t} C_{1}e^{\kappa} t^{-n}2^{n/2}
			{\rm Vol}(\Omega')
			C_{1}e^{\kappa} t^{-n}e^{-\frac{d(y,\Omega')^{2}}{8t}}
			+ 
			4A^{2} {\rm Vol}(\Omega') C_{1}^{2}e^{2\kappa} t^{-2n}e^{-\frac{d(y,\Omega')^{2}}{4t}}
			\\
			&\leq
			4{\rm Vol}(M)C_{1}^{2}e^{2\kappa}(B^{\frac{1}{2}}2^{\frac{n}{2}} + A^{2}te^{-\frac{d(y,\Omega')^{2}}{8t}})
			t^{-(2n+1)}e^{-\frac{d(y,\Omega')^{2}}{8t}}
			\\
			&\leq
			4C_{1}^{2}(2^{\frac{n}{2}}B^{\frac{1}{2}}+A^{2})
			e^{2\kappa}{\rm Vol}(M)t^{-(2n+1)}e^{-\frac{d(y,\Omega')^{2}}{8t}}.
		\end{aligned}
		$$
		We get the second inequality by setting $C_{2}(n,\gamma,A)=4C_{1}(n,\gamma)^{2}(2^{\frac{n}{2}}B^{\frac{1}{2}}+A^{2})$.
	\end{pf}

	\subsection{Estimates for the difference of two heat kernels}
	\par
	Let $\rho$ be a smooth function on $M$ and set
	$$
	\Omega_{c} := \{ x\in M;\, \rho(x)<c \}.
	$$
	We assume the following:
	\begin{itemize}
		\item
		For $1\leq c \leq 3$, $\Omega_{c}$ is a relatively compact domain of $M$.
		\item
		$d\rho\not=0$ on $\overline{\Omega}_{3}\setminus\Omega_{1}$.
		\item
		$S:=\partial\Omega_{1}=\rho^{-1}(1)$ is a compact manifold.
	\end{itemize}
	Then $\Omega_{r}=\Omega_{1}\cup \rho^{-1}([1,r))$ and $\Omega_{r}\setminus\Omega_{1}\cong S \times[1,r)$ for $1 < r \leq 3$.
	We set $S_{r}:=\rho^{-1}(r) \cong S$. Let $d\sigma_{r}$ be the volume form on $S_{r}$ induced by $g^{M}$.
	There are constants $K_{1},K_{2}>0$ such that under the diffeomorphism
	$\overline{\Omega}_{3}\setminus\Omega_{1} \cong S \times [1,3]$,
	\begin{equation}
		\label{eqn:comparison:volume:form}
		K_{1}\,d\rho\wedge d\sigma_{\rho} |_{S_{\rho}\times\{\rho\}} 
		\leq 
		dv|_{S_{\rho}\times\{ \rho \}}
		\leq K_{2} \,d\rho\wedge d\sigma_{\rho} |_{S_{\rho}\times\{ \rho \}}
		\qquad
		(\forall\,\rho \in [1,3]).
	\end{equation}
	\par
	We assume that $(L,h^{L})$ is a trivial Hermitian line bundle on $\overline{\Omega}_{3}$. 
	Recall that the constants $C_{1}, C_{2} >0$ were defined in Lemma~\ref{lemma:Gaussian:nabla:HK}, which depends only on $A$, $n$, $S$.
	For $x,y \in \Omega_{1}$, we define
	$$
	\delta(x,y) := \min\{ d(x, \partial\Omega_{1}), d(y, \partial\Omega_{1}) \}>0.
	$$

	\begin{theorem} 
		\label{thm:diff:HK}
		Set $B'_{m}:=\sup_{x\geq0}x^{m}e^{-x/16}$.
		Then for all $x, y \in \Omega_{1}$ and $0<t\leq 1$, the following inequality holds:
		$$ 
		\left| k(t, x, y)- K^{L}(t, x, y) \right| 
		\leq
		D(M)\,\delta(x,y)^{-2(2n+1)}e^{-\frac{\delta(x,y)^{2}}{16t}},
		$$
		where $D(M) := K_{1}^{-1} e^{2\kappa} {\rm vol}(M) \{ B'_{2n+1} C_{2}(n,\gamma,A) + B'_{4n} C_{1}(n,\gamma)^{2} {\rm diam}(M)^{2} \}$.
	\end{theorem}
	
	\begin{pf}
		We write $C_{1}$, $C_{2}$, $\delta$ for $C_{1}(n,\gamma)$, $C_{2}(n,\gamma,A)$, $\delta(x,y)$, respectively. 
		Fix the trivialization $(L,h^{L}) \cong ({\mathcal O}_{M}, |\cdot|)$ on $\overline{\Omega}_{3}$ as above. 
		Since $\square^{L}=\square$ on $\overline{\Omega}_{3}$, $K^{L}$ satisfies the heat equation $(\partial_{t}+\square_{x})K^{L}(t,x,y)=0$ 
		for $x \in \overline{\Omega}_{3}$, $y \in M$ and $t>0$.
		We apply the Duhamel principle \cite[(3.9)]{C} to the function $k(t, x, y)- K^{L}(t, x, y)$ on 
		${\mathbf R}_{>0}\times\Omega_{\rho} \times \Omega_{\rho}$ $(\rho \leq 3)$. Then we obtain
		\begin{eqnarray*}
			\left| k(t, x, y)- K^{L}(t, x, y) \right|
			& \leq & 
			\int_0^t ds \int_{\partial \Omega_{\rho}}  \left| \nabla_{z}  k(t-s, x, z) \right| \cdot \left| K^{L}(s, z, y) \right| \,d\sigma_{\rho}(z)  
			\\
			& & +  \int_0^t ds \int_{\partial \Omega_{\rho}}  k(t-s, x, z)  \left| \nabla_{z}  K^{L}(s, z, y) \right|\,d\sigma_{\rho}(z)
		\end{eqnarray*}
		for all $x,y\in\Omega_{1}$ and $t>0$. For $1 < \rho \leq 3$ and $x,y \in \Omega_{1}$, we get by the Cauchy-Schwarz inequality
		\begin{equation}
			\label{eqn:diff:HK:bdy:1}
			\begin{aligned}
				\,&
				\left| k(t, x, y)- K^{L}(t, x, y) \right|
				\\
				& \leq
				\int_0^t ds \int_{\partial \Omega_{\rho} }  
				\{ |\nabla_{z}k(t-s, x, z)|  |K^{L}(s, z, y)| + k(t-s, x, z)  |\nabla_{z}K^{L}(s, z, y)| \} d\sigma_{\rho}(z)
				\\
				& \leq
				\frac{1}{2} \int_0^t ds \int_{\partial \Omega_{\rho}} ( |\nabla_{z}k(t-s, x, z)|^{2} + |K^{L}(s, z, y)|^{2} ) d\sigma_{\rho}(z)
				\\
				&\quad+ 
				\frac{1}{2} \int_0^t ds \int_{\partial \Omega_{\rho}}  ( k(t-s, x, z)^{2} + |\nabla_{z}K^{L}(s, z, y)|^{2} ) d\sigma_{\rho}(z).
			\end{aligned}
		\end{equation}
		Integrating \eqref{eqn:diff:HK:bdy:1} with respect to the variable $\rho$ over the interval $[2,3]$ and using \eqref{eqn:comparison:volume:form}, 
		we get the following estimate for all $x, y \in \Omega_{1}$ and $t>0$
		\begin{equation}
			\label{eqn:diff:HK:bdy:2}
			\begin{aligned}
				\,&
				\left| k(t, x, y)- K^{L}(t, x, y) \right|
				=
				\int_{2}^{3} \left| k(t, x, y)- K^{L}(t, x, y) \right| d\rho
				\\
				& \leq
				\frac{1}{2} \int_0^t ds \int_{[2,3]\times\partial \Omega_{\rho}}  ( |\nabla_{z}k(t-s, x, z)|^{2} + |K^{L}(s, z, y)|^{2} ) d\rho\wedge d\sigma_{\rho}(z)
				\\
				&\quad+ 
				\frac{1}{2} \int_0^t ds \int_{[2,3]\times\partial \Omega_{\rho}} ( k(t-s, x, z)^{2} + |\nabla_{z}K^{L}(s, z, y)|^{2} ) d\rho\wedge d\sigma_{\rho}(z)
				\\
				& \leq
				\frac{1}{2K_{1}} \int_{0}^{t} ds \int_{\Omega_{3}\setminus\Omega_{2}}  \left| \nabla_{z}k(t-s, x, z) \right|^{2} dv_z 
				+
				\frac{1}{2K_{1}} \int_{0}^{t} ds\int_{\Omega_{3}\setminus\Omega_{2}} \left| K^{L}(s, z, y) \right|^{2} dv_z 
				\\
				&\quad+ 
				\frac{1}{2K_{1}}\int_{0}^{t} ds \int_{\Omega_{3}\setminus\Omega_{2}} k(t-s, x, z)^{2} dv_z 
				+
				\frac{1}{2K_{1}}\int_{0}^{t} ds \int_{\Omega_{3}\setminus\Omega_{2}} \left| \nabla_{z}K^{L}(s, z, y) \right|^{2} dv_z.
			\end{aligned}
		\end{equation}
		\par
		In Lemma~\ref{lemma:Gaussian:nabla:HK}, consider the case $\Omega=M\setminus\Omega_{2}$ and $\Omega'=M\setminus\Omega_{1}$.
		Then $\Omega_{3}\setminus\Omega_{2} \subset M\setminus\Omega_{2}=\Omega \subset \Omega'=M\setminus\Omega_{1}$ and
		$d(w,\Omega')=d(w,M\setminus\Omega_{1})=d(w,\partial\Omega_{1})$ for all $w\in\Omega_{1}$.
		Hence, by definition of $\delta$, we have $d(x,\Omega') \geq \delta$ and $d(y,\Omega') \geq \delta$.
		Similarly, for $z \in \Omega_{3}\setminus\Omega_{2} \subset \Omega'$ and $x,y \in \Omega_{1}$, we have $d(x,z) \geq \delta$ and $d(y,z) \geq \delta$.
		By Lemma~\ref{lemma:Gaussian:nabla:HK} (2) with $\Omega=M\setminus\Omega_{2}$, $\Omega'=M\setminus\Omega_{1}$ and 
		$y\in \Omega_{1} = M \setminus \Omega'$, we get
		\begin{equation}
			\label{eqn:est:nabla:K:1}
			\begin{aligned}
				\int_{\Omega_{3}\setminus\Omega_{2}} \left| \nabla_{z} K^{L}(s,z,y) \right|^{2} dv_{z} 
				&\leq
				\int_{M\setminus\Omega_{2}} \left| \nabla_{z} K^{L}(s,z,y) \right|^{2} dv_{z} 
				\\
				&=
				\int_{\Omega} \left| \nabla_{z} K^{L}(s,z,y) \right|^{2} dv_{z} 
				\\
				&\leq
				C_{2} e^{2\kappa} {\rm Vol}(M) s^{-(2n+1)} e^{-\frac{d(y,\Omega')^{2}}{8s}}
				\\
				&\leq
				e^{2\kappa} C_{2} B'_{2n+1} {\rm Vol}(M) \delta^{-2(2n+1)} e^{-\frac{\delta^{2}}{16s}}.
			\end{aligned}
		\end{equation}
		Similarly, for $x \in \Omega_{1} = M \setminus \Omega$, we get
		\begin{equation}
			\label{eqn:est:nabla:K:2}
			\int_{\Omega_{3}\setminus\Omega_{2}} \left| \nabla_{z} k(t-s,x,z) \right|^{2} dv_{z} 
			\leq
			C_{2} B'_{2n+1} {\rm Vol}(M) \delta^{-2(2n+1)} e^{-\frac{\delta^{2}}{16(t-s)}}. 
		\end{equation}
		Since $d(z,y) \geq \delta$ for $z \in \Omega_{3} \setminus \Omega_{2}$ and $y \in \Omega_{1}$, we get by
		Lemma~\ref{lemma:Gaussian:nabla:HK} (1) with $y\in \Omega_{1} = M \setminus \Omega'$
		$$
		\left| K^{L}(s,z,y) \right|^{2} 
		\leq 
		C_{1}^{2} e^{2\kappa} s^{-2n} e^{-\frac{d(z,y)^{2}}{4s}} 
		\leq 
		C_{1}^{2} e^{2\kappa} s^{-2n} e^{-\frac{\delta^{2}}{4s}}
		\leq
		C_{1}^{2} e^{2\kappa} B'_{4n} \delta^{-4n} e^{-\frac{\delta^{2}}{8s}}.
		$$
		Hence
		\begin{equation}
			\label{eqn:est:K:1}
			\int_{\Omega_{3}\setminus\Omega_{2}} \left| K^{L}(s,z,y) \right|^{2} dv_{z} 
			\leq  
			e^{2\kappa} C_{1}^{2} B'_{4n} {\rm Vol}(M) \delta^{-4n} e^{-\frac{\delta^{2}}{8s}}.
		\end{equation}
		Similarly, 
		\begin{equation}
			\label{eqn:est:K:2}
			\int_{\Omega_{3}\setminus\Omega_{2}} \left| k(t-s,x,z) \right|^{2} dv_{z} 
			\leq  
			e^{2\kappa} C_{1}^{2} B'_{4n} {\rm Vol}(M) \delta^{-4n} e^{-\frac{\delta^{2}}{8(t-s)}}.
		\end{equation}
		By substituting \eqref{eqn:est:nabla:K:1}, \eqref{eqn:est:nabla:K:2}, \eqref{eqn:est:K:1}, \eqref{eqn:est:K:2} into \eqref{eqn:diff:HK:bdy:2}
		and using the inequalities $e^{-\frac{\delta^{2}}{16s}} \leq e^{-\frac{\delta^{2}}{16t}}$ and $e^{-\frac{\delta^{2}}{16(t-s)}} \leq e^{-\frac{\delta^{2}}{16t}}$ 
		for $0< s <t \leq 1$, the following inequality holds for all $x,y\in \Omega_{1} \subset M \setminus \Omega'$ and $t\in(0,1]$:
		\begin{equation}
			\label{eqn:diff:HK:bdy:3}
			\begin{aligned}
				\left| k(t, x, y)- K^{L}(t, x, y) \right|
				&\leq
				K_{1}^{-1} e^{2\kappa} C_{2}B'_{2n+1} {\rm Vol}(M)\, \delta^{-2(2n+1)} e^{-\frac{\delta^{2}}{16t}}
				\\
				&\quad+
				K_{1}^{-1} e^{2\kappa} C_{1}^{2} B'_{4n} {\rm Vol}(M) \delta^{-4n} e^{-\frac{\delta^{2}}{8t}}.
			\end{aligned}
		\end{equation}
		The result follows from \eqref{eqn:diff:HK:bdy:3}.
	\end{pf}

	\subsection{Uniformity of the asymptotic expansion of the heat kernels}
	
	For $x\in M$, let $i_{x}$ be the injectivity radius at $x$ and set $j_{x} := i_{x}/3$. 
	For $0<r<i_{x}$, set 
	$$
	{\mathbb B}(y,r):=\{ x\in M;\, d(x,y)<r\}.
	$$
	There exist $u_{i}(\cdot,y) \in A^{0}({\mathbb B}(y,j_{y}), L \otimes L_{y}^{\lor})$ $(i\geq0)$ such that 
	$$
	p(t,x,y) = (4\pi t)^{-n}\exp\left(-\frac{d(x,y)^{2}}{4t}\right) \sum_{i=0}^{\infty} t^{i} u_{i}(x,y)
	$$
	is a formal solution of the heat equation $(\partial_{t} + \square^{L}_{x}) p(t,x,y) =0$ with $u_{0}(y,y) =1$.
	(See \cite[Th.\,2.26]{BerlineGetzlerVergne92} for an explicit formula for $u_{i}(x,y)$.)
	Let $k>n+4$. We set
	$$
	p_{k}(t,x,y) := (4\pi t)^{-n}\exp\left(-\frac{d(x,y)^{2}}{4t}\right) \left\{ u_{0}(x,y) + tu_{1}(x,y) + \cdot\cdot\cdot + t^{k} u_{k}(x,y) \right\},
	$$
	$$
	F_{k}(t,x,y) := K^{L}(t,x,y)-p_{k}(t,x,y).
	$$
	For any $y\in M$, $F_{k}(t,\cdot,y)$ is defined on ${\mathbb B}(y,j_{y})$. 
	By setting $F_{k}(t,\cdot,y)=0$ for $t\leq 0$, $F_{k}(\cdot,\cdot,y)$ extends to a $C^{2}$-function on ${\mathbf R} \times {\mathbb B}(y,j_{y})$  
	by \cite[Th.~2.23 (2)]{BerlineGetzlerVergne92}. We define 
	$$
	B_{x} := j(x)^{1/2} \circ \square^{L}_{x} \circ j^{1/2}(x).
	$$
	Here, if $x=\exp_{y}({\mathbf x})$ with ${\mathbf x}=({\mathbf x}_{1},\ldots,{\mathbf x}_{2n})$ being the geodesic normal coordinates centered at $y$ and 
	$g(x) = \sum_{i,j}g_{ij}({\mathbf x}) d{\mathbf x}_{i}d{\mathbf x}_{j}$, then $j(x) = \det( g_{ij}( {\mathbf x} ) )^{1/2}$.
	By \cite[Prop.\,2.24 and Th.\,2.26]{BerlineGetzlerVergne92}, for any $(t,x) \in {\mathbf R}_{>0} \times {\mathbb B}(y,j_{y})$, we have
	$$
	\left( \frac{\partial}{\partial t}+\square^{L}_{x}\right) F_{k}(t,x,y)
	=
	(4\pi)^{-n}t^{k-n}\exp\left(-\frac{d(x,y)^{2}}{4t}\right)B_{x}u_{k}(x,y).
	$$
	Set
	
	\begin{equation}
		\label{eqn:constants}
		C_{k}(y):=\sup_{x\in{\mathbb B}(y,j_{y})} \left| u_{k}(x,y) \right|,
		\qquad
		D_{k}(y):=(4\pi)^{-n}\sup_{x\in{\mathbb B}(y,j_{y})} \left| B_{x}u_{k}(x,y) \right|.
	\end{equation}
	
	If the geometry of $({\mathbb B}(y,j_{y}),g^{M})$ is uniformly bounded, then $C_{i}(y)$ and $D_{i}(y)$ $(0\leq i \leq k)$
	are also uniformly bounded by construction of $u_{i}(x,y)$ in \cite[Th.\,2.26]{BerlineGetzlerVergne92}. 
	\par
	Let $\chi_{y}\in C^{\infty}(M)$ be a nonnegative function such that $\chi_{y}(x)=1$ on ${\mathbb B}(y,\frac{1}{2}j_{y})$, 
	$\chi_{y}(x)=0$ on $M\setminus{\mathbb B}(y,j_{y})$ and $|d\chi_{y}|\leq 4j^{-1}_{y}$. 
	We define
	$$
	\begin{aligned}
		G_{k}(t,x,y) 
		&:= 
		\chi_{y}(x)\left(\frac{\partial}{\partial t}+\square^{L}_{x}\right)F_{k}(t,x,y)
		\\
		&=
		(4\pi)^{-n}t^{k-n}\chi_{y}(x) \exp\left( -\frac{d(x,y)^{2}}{4t} \right) B_{x}u_{k}(x,y).
	\end{aligned}
	$$
	Then
	\begin{equation}
		\label{eqn:estimate:G:k}
		\left| G_{k}(t,x,y) \right| \leq t^{k-n} D_{k}(y) \exp\left( -\frac{d(x,y)^{2}}{4t} \right).
	\end{equation}
	Set
	$$
	\begin{aligned}
		H_{k}(t,x,y) 
		&:=
		\int^{t}_{0}d\tau\int_{M}K^{L}(t-\tau,x,z)G_{k}(\tau,z,y)\,dv(z)
		\\
		&=
		\int^{t}_{0}d\tau\int_{{\mathbb B}(y,j_{y})}K^{L}(t-\tau,x,z)G_{k}(\tau,z,y)\,dv(z).
	\end{aligned}
	$$
	Then $H_{k}(t,x,y)$ satisfies the heat equation 
	\begin{equation}
		\label{eqn:heat:eq:H:k}
		(\partial_{t}+\square^{L}_{x}) H_{k}(t,x,y) = G_{k}(t,x,y)  = \chi_{y}(x)\left(\partial_{t}+\square^{L}_{x}\right)F_{k}(t,x,y)
	\end{equation}
	with $\lim_{t\to 0} H_{k}(t,x,y) =0$. Since $\chi_{y}=1$ on ${\mathbb B}(y,j_{y}/2)$, we get
	\begin{equation}
		\label{eqn:heat:eq}
		(\partial_{t}+\square^{L}_{x})\{ F_{k}(t,x,y)-H_{k}(t,x,y) \} = 0
		\qquad
		(\forall\,x\in{\mathbb B}(y, j_{y}/2),\,\, t>0).
	\end{equation}
	Recall that $\kappa=\max_{x\in M}|Q^{L}(x)|$, where $Q^{L}=2\pi\Lambda c_{1}(L,h)$.

	\begin{lemma}
		\label{lemma:sup:Fk-Hk}
		For all $t\in (0,1]$, the following inequality holds:
		$$
		\begin{aligned}
			\,&
			\sup_{(0,t]\times{\mathbb B}(y,\frac{1}{2} j_{y})}\left| F_{k}(\cdot,\cdot,y)-H_{k}(\cdot,\cdot,y) \right|
			\\
			&\leq 
			e^{\kappa/2} \{ \sup_{(0,t]\times\partial{\mathbb B}(y,\frac{1}{2} j_{y})}|F_{k}(\cdot,\cdot,y)|
			+\sup_{(0,t]\times\partial{\mathbb B}(y,\frac{1}{2} j_{y})}|H_{k}(\cdot,\cdot,y)| \}.
		\end{aligned}
		$$
	\end{lemma}
	
	\begin{pf}
		Recall that $\nabla^{L}=\partial^{L}+\bar{\partial}^{L}$ is the Chern connection of $(L, h)$. Then for any $s \in A^{0}(L)$,
		$$
		\partial\bar{\partial}h(s,s) = h(\partial^{L}\bar{\partial}^{L}s, s) - h(\bar{\partial}^{L}s, \bar{\partial}^{L}s) + h({\partial}^{L}s, {\partial}^{L}s) 
		+ h(s, \bar{\partial}^{L}{\partial}^{L}s).
		$$ 
		This, together with $\square^{L} = (\bar{\partial}^{L})^{*}\bar{\partial}^{L} = -\sqrt{-1}\Lambda{\partial}^{L}\bar{\partial}^{L}$,
		$\overline{\square}^{L} = ({\partial}^{L})^{*}{\partial}^{L} =\sqrt{-1}\Lambda \bar{\partial}^{L}{\partial}^{L}$, and the Bochner-Kodaira-Nakano formula
		$\overline{\square}^{L} - {\square}^{L} = 2\pi\Lambda c_{1}(L,h)=Q^{L}$ on $A^{0}(L)$, yields that
		$$
		\square h(s,s) = h({\square}^{L}s, s) + h(s, {\square}^{L}s) + Q^{L} h(s,s) - |\nabla^{L}s|^{2}
		\qquad
		(s\in A^{0}(L)).
		$$
		By this equality, we get for any $s\in A^{0}({\mathbf R}_{>0}\times M, L)$,
		$$
		\left( \partial_{t} + \frac{1}{2}\Delta \right) h(s,s) 
		= 
		-| \nabla^{L} s |^{2} + h( (\square^{L}+\partial_{t})s, s) + h(s, (\square^{L}+\partial_{t})s) + Q^{L} h(s,s).
		$$
		Putting $s=F_{k}(\cdot,\cdot,y)-H_{k}(\cdot,\cdot,y)$ in this equality and using \eqref{eqn:heat:eq}, we get
		$$
		\begin{aligned}
			\left( \partial_{t} + \frac{1}{2}\Delta_{x} \right) \left| F_{k}(t,x,y)-H_{k}(t,x,y) \right|^{2}
			&\leq
			Q^{L}(x) \left| F_{k}(t,x,y)-H_{k}(t,x,y) \right|^{2} 
			\\
			&\leq \kappa \left| F_{k}(t,x,y)-H_{k}(t,x,y) \right|^{2}
		\end{aligned}
		$$
		for all $x\in{\mathbb B}(y,j_{y}/2)$ and $t>0$. Namely, on ${\mathbf R}_{>0} \times {\mathbb B}(y,j_{y}/2)$, we have
		$$
		\left( \partial_{t} + \frac{1}{2}\Delta_{x} \right) \left\{e^{-\kappa t}\left| F_{k}(t,x,y)-H_{k}(t,x,y) \right|^{2}\right\} \leq 0.
		$$
		From the weak maximum principle for subsolutions of the heat operator, it follows that for $t\in(0,1]$,
		$$
		\begin{aligned}
			\,&
			e^{-\kappa}\max_{[0,t]\times{\mathbb B}(y,j_{y}/2)} \left| F_{k}(\cdot,\cdot,y) - H_{k}(\cdot,\cdot,y) \right|^{2}
			\\
			&\leq
			\max_{(\tau,x)\in[0,t]\times{\mathbb B}(y,j_{y}/2)} \left\{ e^{-\kappa \tau}\left| F_{k}(\tau,x,y) - H_{k}(\tau,x,y) \right|^{2} \right\}
			\\
			&\leq
			\max_{(\tau,x)\in ([0,t]\times\partial{\mathbb B}(y,j_{y}/2))\cup(\{0\}\times{\mathbb B}(y,j_{y}/2))} 
			\left\{ e^{-\kappa \tau}\left| F_{k}(\tau,x,y) - H_{k}(\tau,x,y) \right|^{2} \right\}
			\\
			&\leq
			\max_{(\tau,x)\in[0,t]\times\partial{\mathbb B}(y,j_{y}/2)} \left| F_{k}(\tau,x,y) - H_{k}(\tau,x,y)  \right|^{2}
			\\
			&\leq
			\left( \max_{(\tau,x)\in[0,t]\times\partial{\mathbb B}(y,j_{y}/2)} \left| F_{k}(\tau,x,y) \right| + 
			\max_{(\tau,x)\in[0,t]\times\partial{\mathbb B}(y,j_{y}/2)} \left| H_{k}(\tau,x,y)  \right| 
			\right)^{2},
		\end{aligned}
		$$
		where we used $\kappa\geq0$ and $F_{k}(0,x,y)=H_{k}(0,x,y)=0$ for $x\in {\mathbb B}(y,j_{y}/2)$ to get the third inequality. 
		The result follows from this inequality.
	\end{pf}

	\begin{lemma}
		\label{lemma:sup:Fk}
		Set $B(n):=\sup_{x>0}x^{n}e^{-x/64}$. Then for all $t\in (0, 1]$,
		$$
		\sup_{(0,t]\times\partial{\mathbb B}(y,\frac{1}{2}j_{y})}|F_{k}(\cdot,\cdot,y)|
		\leq
		\widetilde{C}_{1}(y)\,\exp\left(-\frac{j_{y}^{2}}{64t}\right),
		$$
		where $\widetilde{C}_{1}(y) =(k+1)C_{1}e^{\kappa}B(n)j_{y}^{-2n} \max_{1\leq i \leq k} C_{i}(y)$ 
		with the constants $C_{i}(y)$ defined in \eqref{eqn:constants}.
	\end{lemma}
	
	\begin{pf}
		For $(s,x)\in (0,t]\times\partial{\mathbb B}(y,j_{y}/2)$, we get by Lemma~\ref{lemma:Gaussian:nabla:HK} (1)
		$$
		\begin{aligned}
			|F_{k}(s,x,y)|
			&\leq 
			C_{1}e^{\kappa}\,s^{-n}\exp\left(-\frac{d(x,y)^{2}}{8s}\right)\{1+C_{1}(y)s+ \cdots + C_{k}(y)s^{k}\}
			\\
			&\leq 
			(k+1)C_{1}e^{\kappa}\max_{1\leq i \leq k}C_{i}(y)\,s^{-n}\exp\left(-\frac{j_{y}^{2}}{32s}\right)
			\\
			&\leq
			(k+1)C_{1}e^{\kappa}B(n)\max_{1\leq i \leq k}C_{i}(y)\,j_{y}^{-2n}\exp\left(-\frac{j_{y}^{2}}{64t}\right).
		\end{aligned}
		$$
		This proves the result.
	\end{pf}

	In what follows, we assume $k>n+4$. For $y\in M$, we set
	$$
	E(y):=\sup_{\exp_{y}(\xi)\in{\mathbb B}(y,2j_{y})} \frac{ (\exp_{y})^{*}(dv)(\xi)}{i^{n^{2}}d\xi\wedge \overline{d\xi}} \geq 1.
	$$

	\begin{lemma}
		\label{lemma:sup:Hk}
		For $t\in (0, 1]$, one has
		$$
		\sup_{(0,1]\times\partial{\mathbb B}(y,\frac{1}{2}j_{y})} \left| H_{k}(\cdot,\cdot,y) \right|
		\leq
		\widetilde{C}_{2}(y)\,t^{k+1-n}\exp\left(-\frac{j_{y}^{2}}{256t}\right),
		$$
		where $\widetilde{C}_{2}(y) = (16\pi)^{n}(k-n+1)^{-1}C_{1}e^{\kappa}D_{k}(y) \sup_{x\in{\mathbb B}(y,\frac{1}{2}j_{y})}E(x)$
		with the constants $D_{k}(y)$ defined in \eqref{eqn:constants}.
	\end{lemma}
	
	\begin{pf}
		Let $(s,x)\in (0,t]\times\partial{\mathbb B}(y,\frac{1}{2}j_{y})$. Then we have
		\begin{equation}
			\label{eqn:estimate:H:1}
			\begin{aligned}
				\left| H_{k}(s,x,y) \right|
				&=
				\left| \int^{s}_{0}d\tau\int_{{\mathbb B}(y,j_{y})} K(s-\tau,x,z)G_{k}(\tau,z,y) dv(z) \right|
				\\
				&\leq
				\int^{s}_{0}d\tau\int_{{\mathbb B}(y,\frac{1}{4}j_{y})} \left| K(s-\tau,x,z) \right| \cdot \left| G_{k}(\tau,z,y) \right| dv(z)
				\\
				&\quad+
				\int^{s}_{0}d\tau\int_{{\mathbb B}(y,j_{y})\setminus{\mathbb B}(y,\frac{1}{4}j_{y})} \left| K(s-\tau,x,z) \right| \cdot \left| G_{k}(\tau,z,y)\right| dv(z).
			\end{aligned}
		\end{equation}
		Since $x\in\partial{\mathbb B}(y,\frac{1}{2}j_{y})$, we have $d(x,z)\geq d(z,y)$ and $d(x,z)\geq\frac{1}{4}j_{y}$ for
		$z\in{\mathbb B}(y,\frac{1}{4}j_{y})$. 
		By Lemma~\ref{lemma:Gaussian:nabla:HK} (1) and \eqref{eqn:estimate:G:k}, we get
		\begin{equation}
			\label{eqn:estimate:H:2}
			\begin{aligned}
				\,&
				\int^{s}_{0}d\tau\int_{{\mathbb B}(y,\frac{1}{4}j_{y})} \left| K(s-\tau,x,z)\right| \cdot \left| G_{k}(\tau,z,y)\right| dv(z)
				\\
				&\leq 
				\int^{s}_{0}d\tau\int_{{\mathbb B}(y,\frac{1}{4}j_{y})}
				C_{1}e^{\kappa}(s-\tau)^{-n}\exp\left(-\frac{d(x,z)^{2}}{8(s-\tau)}\right)
				\tau^{k-n} D_{k}(y) \exp\left( -\frac{d(z,y)^{2}}{4\tau} \right)
				dv(z)
				\\
				&\leq 
				C_{1}e^{\kappa}D_{k}(y)
				\int^{s}_{0} d\tau
				\int_{{\mathbb B}(y,\frac{1}{4}j_{y})}\left\{(s-\tau)^{-n}\exp\left(-\frac{j_{y}^{2}}{256(s-\tau)}\right)
				\exp\left(-\frac{d(y,z)^{2}}{16(s-\tau)}\right)\right.
				\\
				&\qquad\qquad\qquad\qquad\qquad\times\left.\tau^{k-n} \exp\left( -\frac{d(z,y)^{2}}{16\tau} \right)\right\}
				dv(z)
			\end{aligned}
		\end{equation}
		Write $z=\exp_{y}(\xi)$, where $\xi=(\xi_{1},\ldots,\xi_{2n})$ is the system of
		geodesic normal coordinates centered at $y$. Then $d(y,z)^{2}=\| \xi \|^{2} = \sum_{i}(\xi_{i})^{2}$.
		We set $dV(\xi):=i^{n^{2}}d\xi\wedge \overline{d\xi}$. By \eqref{eqn:estimate:H:2}, we get
		\begin{equation}
			\label{eqn:estimate:H:3}
			\begin{aligned}
				\,&
				\int^{s}_{0}d\tau\int_{{\mathbb B}(y,\frac{1}{4}j_{y})} \left| K(s-\tau,x,z) \right| \cdot \left| G_{k}(\tau,z,y)\right| dv(z)
				\\
				&\leq 
				C_{1}e^{\kappa}D_{k}(y)\,
				s^{-n}\exp\left(-\frac{j_{y}^{2}}{256s}\right)
				\int^{s}_{0}\tau^{k}d\tau
				\\
				&\qquad\times
				\int_{\| \xi \|\leq\frac{1}{4}j_{y}}
				\left\{\frac{\tau(s-\tau)}{s}\right\}^{-n} \exp\left(-\frac{s\|\xi\|^{2}}{16\tau(s-\tau)}\right)E(y)\,dV(\xi)
				\\
				&\leq 
				\frac{(16\pi)^{n}C_{1}e^{\kappa}D_{k}(y)E(y)}{k+1}\,s^{k+1-n}\exp\left(-\frac{j_{y}^{2}}{256s}\right).
			\end{aligned}
		\end{equation}
		Since $d(z,y)\geq\frac{1}{4}j_{y}$ and $d(x,z)\leq d(x,y)+d(y,z) \leq \frac{3}{2}j_{y}$ for $z\in{\mathbb B}(y,j_{y})\setminus{\mathbb B}(y,\frac{1}{4}j_{y})$,
		we get by Lemma~\ref{lemma:Gaussian:nabla:HK} (1) and \eqref{eqn:estimate:G:k}
		\begin{equation}
			\label{eqn:estimate:H:4}
			\begin{aligned}
				\,&
				\int^{s}_{0}d\tau\int_{{\mathbb B}(y,j_{y})\setminus{\mathbb B}(y,\frac{1}{4}j_{y})} \left| K(s-\tau,x,z) \right| \cdot \left| G_{k}(\tau,z,y)\right| dv(z)
				\\
				&\leq 
				\int^{s}_{0}d\tau\int_{{\mathbb B}(y,j_{y})\setminus{\mathbb B}(y,\frac{1}{4}j_{y})}
				C_{1}e^{\kappa}(s-\tau)^{-n}\exp\left(-\frac{d(x,z)^{2}}{8(s-\tau)}\right)
				\tau^{k-n} D_{k}(y) \exp\left( -\frac{d(z,y)^{2}}{4\tau} \right)
				dv(z)
				\\
				&\leq 
				C_{1}e^{\kappa}D_{k}(y)\,\exp\left(-\frac{j_{y}^{2}}{64s}\right)
				\int^{s}_{0}\tau^{k-n}d\tau\int_{{\mathbb B}(y,j_{y})\setminus{\mathbb B}(y,\frac{1}{4}j_{y})}
				(s-\tau)^{-n}\exp\left(-\frac{d(x,z)^{2}}{8(s-\tau)}\right) dv(z)
				\\
				&\leq 
				C_{1}e^{\kappa}D_{k}(y)\,\exp\left(-\frac{j_{y}^{2}}{64s}\right)
				\int^{s}_{0}\tau^{k-n}d\tau\int_{{\mathbf C}^{n}}
				(s-\tau)^{-n}\exp\left(-\frac{\| \xi \|^{2}}{8(s-\tau)}\right) E(x)\,dV(\xi)
				\\
				&\leq 
				\frac{(8\pi)^{n}C_{1}e^{\kappa}D_{k}(y)\sup_{x\in{\mathbb B}(y,j_{y}/2)}E(x)}{k-n+1}\,s^{k+1-n}\exp\left(-\frac{j_{y}^{2}}{64s}\right).
			\end{aligned}
		\end{equation}
		By \eqref{eqn:estimate:H:1}, \eqref{eqn:estimate:H:3}, \eqref{eqn:estimate:H:4}, we get
		$$
		\left| H_{k}(s,x,y) \right|
		\leq 
		\frac{(16\pi)^{n}C_{1}e^{\kappa}D_{k}(y)\sup_{x\in{\mathbb B}(y,j_{y}/2)}E(x)}{k-n+1}\,s^{k+1-n}\exp\left(-\frac{j_{y}^{2}}{256s}\right).
		$$
		The result follows from this inequality.
	\end{pf}

	\begin{proposition}
		\label{prop:Hk-Fk}
		For $t\in (0, 1]$, the following inequality holds:
		$$
		\sup_{(0,t]\times{\mathbb B}(y,\frac{1}{2}j_{y})} \left| H_{k}(\cdot,\cdot,y)-F_{k}(\cdot,\cdot,y) \right|
		\leq
		\widetilde{C}_{3}(y)\,t^{k+1-n} \exp\left(-\frac{j_{y}^{2}}{256t}\right),
		$$
		where
		$$
		\widetilde{C}_{3}(y)
		=
		(k+1)C_{1}e^{3\kappa/2}B(n)j_{y}^{-2n} \max_{1\leq i \leq k} C_{i}(y)
		+
		\frac{(16\pi)^{n}C_{1}e^{3\kappa/2}D_{k}(y)}{k-n+1} \sup_{x\in{\mathbb B}(y,\frac{1}{2}j_{y})}E(x).
		$$
	\end{proposition}
	
	\begin{pf}
		The result follows from Lemmas~\ref{lemma:sup:Fk-Hk}, \ref{lemma:sup:Fk}, \ref{lemma:sup:Hk}.
	\end{pf}

	Next, we estimate $H_{k}$ on the diagonal.

	\begin{proposition}
		\label{prop:Hk:diag}
		For all $t\in (0,1]$ and $y \in M$, 
		$$
		\left| H_{k}(t,y,y) \right| \leq \widetilde{C}_{4}(y)\,t^{k+1-n},
		$$
		where $\widetilde{C}_{4}(y) = (2\pi)^{n}C_{1}e^{\kappa}D_{k}(y)E(y)/(k+1)$.
	\end{proposition}
	
	\begin{pf}
		Since
		$$
		\left| H_{k}(t,y,y) \right| \leq \int^{t}_{0}d\tau\int_{{\mathbb B}(y,j_{y})} \left| K^{L}(t-\tau,y,z) \right| \cdot \left| G_{k}(\tau,z,y) \right|\,dv(z),
		$$
		we get by Lemma~\ref{lemma:Gaussian:nabla:HK} (1) and \eqref{eqn:estimate:G:k}
		$$
		\begin{aligned}
			\,&
			\left| H_{k}(t,y,y) \right|
			\\
			&\leq 
			\int^{t}_{0}d\tau \int_{{\mathbb B}(y,j_{y})}
			C_{1}e^{\kappa}(t-\tau)^{-n}\exp\left(-\frac{d(y,z)^{2}}{8(t-\tau)}\right)
			\tau^{k-n}D_{k}(y)\exp\left(-\frac{d(y,z)^{2}}{4\tau}\right)dv(z)
			\\ 
			&\leq 
			C_{1}e^{\kappa}D_{k}(y) \int^{t}_{0}\tau^{k-n}d\tau
			\int_{\|\xi\|<j_{y}}(t-\tau)^{-n} \exp\left\{ -\frac{1}{8}\| \xi \|^{2}\left( \frac{1}{\tau}+\frac{1}{t-\tau}\right)\right\} E(y)\,dV(\xi)
			\\
			&\leq 
			C_{1}e^{\kappa}D_{k}(y)E(y)\,t^{-n}\int^{t}_{0}\tau^{k}d\tau
			\int_{{\mathbf C}^{n}}\left\{\frac{\tau(t-\tau)}{t}\right\}^{-n}\exp\left(-\frac{t \| \xi \|^{2}}{8\tau(t-\tau)}\right) dV(\xi)
			\\
			&=
			\frac{(2\pi)^{n}C_{1}e^{\kappa}D_{k}(y)E(y)}{k+1}\,t^{k+1-n}.
		\end{aligned}
		$$
		This proves the result.
	\end{pf}

	\begin{theorem}
		\label{thm:asym:exp:H:K}
		Let $k>n+4$. For all $t\in (0, 1]$ and $y\in M$, the following inequality holds:
		$$
		\left| K^{L}(t,y,y)-p_{k}(t,y,y) \right| \leq \widetilde{D}_{k}(y)\,t^{k+1-n},
		$$
		where the constant $\widetilde{D}_{k}(y)$ is given by
		$$
		\begin{aligned}
			\widetilde{D}_{k}(y)
			&=
			(k+1)C_{1}e^{\kappa}B(n)j_{y}^{-2n} \max_{1\leq i \leq k} C_{i}(y)
			+
			\frac{(16\pi)^{n}C_{1}e^{3\kappa/2}D_{k}(y)}{k-n+1} \sup_{x\in{\mathbb B}(y,\frac{1}{2}j_{y})}E(x)
			\\
			&\quad
			+ \frac{(2\pi)^{n}C_{1}e^{\kappa}D_{k}(y)E(y)}{k+1}.
		\end{aligned}
		$$
	\end{theorem}
	
	\begin{pf} 
		Since
		$$
		\begin{aligned}
			\left| K^{L}(t,y,y)-p_{k}(t,y,y) \right| 
			&=
			\left| F_{k}(t,y,y) \right|
			\\
			&\leq
			\sup_{(0,t]\times{\mathbb B}(y,\frac{1}{2}j_{y})} \left| H_{k}(\cdot,\cdot,y) - F_{k}(\cdot,\cdot,y) \right| + \left| H_{k}(t,y,y) \right|,
		\end{aligned}
		$$
		the result follows from Propositions~\ref{prop:Hk-Fk} and \ref{prop:Hk:diag}.
		This completes the proof.
	\end{pf}

	\section{Partial analytic torsions and the ratio of analytic torsions}
	\label{sect:3}
	\par
	In the rest of this article, we maintain the notation from the introduction.
	Let $(L,h^{L})$ be a holomorphic Hermitian line bundle on $X$. We assume that $H := L^{-1}$ is ample and that the Chern form $c_{1}(L, h^{L})$ is semi-negative
	on $X$ and vanishes on a neighborhood of ${\rm Sing}\,X_{0}$ in $X$. The existence of such a Hermitian metric will be shown in Lemma~\ref{lemma:HM} below.
	In this section, we compare the analytic torsions $\tau(X_{s},{\mathcal O}_{X_{s}})$ and $\tau(X_{s},L_{s})$, where we set $L_{s} := L|_{X_{s}}$.
	\par
	Since $X$ admits an ample line bundle $H=L^{-1}$, by shrinking $S$ if necessary, there exists an embedding 
	$\iota \colon X \hookrightarrow S \times {\mathbb P}^{N}$ such that $f = {\rm pr}_{1}\circ \iota$. Let $g_{{\mathbb P}^{N}}$ be the Fubini-Study metric on
	$g_{{\mathbb P}^{N}}$ and let $g^{X}_{\rm FS}$ be the K\"ahler metric on $X$ defined as $g^{X}_{\rm FS} = \iota^{*}(ds\otimes d\bar{s} + g_{{\mathbb P}^{N}})$.
	Shrinking $S$ again if necessary, there exists a constant $\gamma>0$ such that
	$$
	\gamma^{-1} g^{X}_{\rm FS} \leq g^{X} \leq \gamma g^{X}_{\rm FS}.
	$$
	By this inequality, we have the following inequality for all $s \in S^{o}$:
	\begin{equation}
		\label{eqn:comparison:X:FS}
		\gamma^{-1} \iota_{s}^{*}g_{{\mathbb P}^{N}} \leq g_{s} \leq \gamma\, \iota_{s}^{*} g_{{\mathbb P}^{N}},
	\end{equation}
	where $\iota_{s} := \iota|_{X_{s}}$ and $g_{s} = g^{X}|_{X_{s}}$.

	\subsection{Analytic torsion}
	\par
	Let us recall the definition of analytic torsion for compact Riemann surfaces.
	Let $(M,h^{M})$ be a compact Riemann surface endowed with a K\"ahler metric.
	Let $(E,h^{E})$ be a holomorphic Hermitian vector bundle on $M$.
	Let $\square_{0,q}=(\bar{\partial}+\bar{\partial}^{*})^{2}$ be the Laplacian acting on $A^{0,q}_{M}(E)$.
	Let $\zeta_{0,q}(s)$ be the zeta function of $\square_{0,q}$:
	$$
	\zeta_{0,q}(s)
	:=
	\sum_{\lambda\in\sigma(\square_{0,q})\setminus\{0\}} \frac{ \dim E(\lambda,\square_{0,q}) }{ \lambda^{s} }
	=
	\frac{1}{\Gamma(s)}\int_{0}^{\infty} \{ {\rm Tr}\,e^{-t\square_{0,q}}-h^{0,q}(E)\} t^{s-1}\,dt,
	$$
	where $E(\lambda,\square_{0,q})$ is the eigenspace of $\square_{0,q}$ corresponding to the eigenvalue $\lambda$.
	\par
	The analytic torsion of $(M,E)$ with respect to the metrics $h_{M}$, $h_{E}$ is the real number
	$$
	\tau(M,E) := \exp\{ -\sum_{q\geq0}(-1)^{q}q\,\zeta'_{0,q}(0) \} = e^{\zeta'_{0,1}(0)} = e^{\zeta'_{0,0}(0)} .
	$$

	\subsection{Partial analytic torsion}
	\par
	Let $B( {\rm Sing}\,X_{0}, \delta) = \bigcup_{p\in {\rm Sing}\,X_{0}} B( p, \delta)$, where $B( p, \delta)$ is the open metric ball of radius $\delta>0$ 
	centered at $p\in {\rm Sing}\,X_{0}$. Let $v$ be a $C^{\infty}$ complex vector field of $X \setminus {\rm Crit}(f)$ such that $f_{*} v  = \partial/\partial s$.
	Let $0< \epsilon_{\delta} \ll \delta$ be a sufficiently small positive number. We denote by $\Delta(\epsilon_{\delta}) \subset {\mathbf C}$ the closed disc of radius $\epsilon_{\delta}$ centered at the origin.
	Integrating $v$, we have a $C^{\infty}$ map
	$$
	\Phi \colon \Delta(\epsilon_{\delta}) \times ( X_{0} \setminus B( {\rm Sing}\,X_{0}, \delta) ) \hookrightarrow X
	$$ 
	satisfying the following conditions (cf. Section~\ref{sect:An upper bound of the small eigenvalues} or \cite[Proof of Th.\,2.3]{Kodaira86}):
	\begin{itemize}
		\item[(1)] 
		$\Phi$ is a diffeomorphism from $\Delta(\epsilon_{\delta}) \times (X_{0} \setminus B( {\rm Sing}\,X_{0}, \delta) )$ to its image.
		\item[(2)] 
		$\Phi_{s} := \Phi(s, \cdot)$ sends $\{ s \} \times ( X_{0} \setminus B( {\rm Sing}\,X_{0}, \delta))$ to $X_{s}$.
		\item[(3)] 
		$\Phi_{0}=\Phi(0,\cdot) = {\rm id}_{X_{0}}|_{X_{0} \setminus B( {\rm Sing}\,X_{0}, \delta)}$.
		\item[(4)] 
		$\frac{1}{2} g_{0}|_{X_{0}\setminus B({\rm Sing}\,X_{0},\delta)} \leq \Phi_{s}^{*}g_{s} \leq 2g_{0}|_{X_{0}\setminus B({\rm Sing}\,X_{0},\delta)}$.
		\item[(5)] 
		$\{ \Phi_{s}^{*}g_{s} \}_{s\in\Delta(\epsilon_{\delta})}$ converges to $g_{0}|_{X_{0}\setminus B({\rm Sing}\,X_{0},\delta)}$ in the $C^{\infty}$-topology.
	\end{itemize}
	We define 
	$$
	\Omega_{\delta} := f^{-1}(\Delta(\epsilon_{\delta})) \setminus \Phi( \Delta(\epsilon_{\delta}) \times ( X_{0} \setminus B( {\rm Sing}\,X_{0}, \delta) ) ).
	$$
	We fix $0<\delta_{0}\ll1$ and we write $\Omega$ for $\Omega_{\delta_{0}}$. Shrinking $S$ if necessary, namely replacing $S$ with $\Delta(\epsilon_{\delta_{0}})$, 
	we can assume $\Phi_{s}$ is defined for all $s \in S$ and
	$$
	\Omega = X \setminus \Phi( S \times ( X_{0} \setminus B( {\rm Sing}\,X_{0}, \delta_{0}) ) ).
	$$
	Then $\Omega$ is an open neighborhood of ${\rm Crit}(f) = {\rm Sing}\,X_{0}$ in $X$.

	Since ${\rm Sing}\,X_{0}$ consists of isolated points of $X$, the following lemma is well-known. 
	For completeness, we give its proof.
	
	\begin{lemma}
		\label{lemma:HM}
		There exists a Hermitian metric $h^{L}$ on $L$ with semi-negative Chern form
		such that $(L,h^{L})$ is flat on a neighborhood of ${\rm Sing}\,X_{0}$.
	\end{lemma}
	
	\begin{pf}
		Let $c\in(0,1/2)$ be a small number. Then there exists a smooth convex increasing function $F_{c}\in C^{\infty}([0,1])$ such that
		$F_{c}(t)=0$ for $t\leq c$ and $F_{\delta}(t) = t + A_{c}$ for $t\geq2c$, where $A_{c}$ is a constant.
		For instance, if $G\in C^{\infty}([0,1])$ is a non-negative increasing function such that $G(t)=0$ for $t\leq c$ and $G(t)=1$ for $t \geq 2c$, 
		we define $F(t):=\int_{0}^{t} G(s)\,ds$. Then $F(t)$ is a desired convex increasing function with $A_{c} = -2c+\int_{c}^{2c} G(s)\,ds$.
		\par
		Let $h^{H}$ be a Hermitian metric on $H=L^{-1}$ with positive first Chern form. Let $p\in {\rm Sing}\,X_{0}$. 
		Let $\sigma$ be a local defining section of $H$ defined on a coordinate neighborhood $(U,z)$ centered at $p$ with $h^{H}(\sigma(z),\sigma(z))(p)=1$. 
		Set $\varphi(z):=-\log h^{H}(\sigma(z),\sigma(z))$. Since $i\partial\bar{\partial}\varphi>0$ is a K\"ahler form on $U$,
		we may assume by changing the local coordinates suitably that $\varphi(z) = \|z\|^{2} + O(\|z\|^{4})$ on $U$.
		If $c \in (0,1)$ is sufficiently small, we may assume that $0\leq\varphi(z) < c$ for $\|z\|<\sqrt{c/2}$ and that
		$\varphi(z)>2c$ for $\|z\|>\sqrt{3c}$. Under this condition, we set $\psi(z) := F_{c}(\varphi(z))$. 
		Then $\psi(z)=0$ for $\|z\|<\sqrt{c/2}$ and $\psi(z) = \varphi(z) + A_{c}$ for $\|z\|>\sqrt{3c}$.
		Since $F'_{c}\geq0$ and $F''_{c}\geq0$, we see that 
		$i\partial\bar{\partial}\psi = F'_{c}(\varphi)i\partial\bar{\partial}\varphi+F''_{c}(\varphi)i\partial\varphi\bar{\partial}\varphi$
		is a semi-positive $(1,1)$-form. Moreover, there are open subsets $W \subset \overline{W} \subset V \subset \overline{V} \subset U$
		such that $\psi=0$ on $W$ and $\psi = \varphi + A_{c}$ on $U\setminus V$.
		We define a Hermitian metric $\widetilde{h}_{H}$ on $H$ by $\widetilde{h}_{H}(\sigma,\sigma)(z): = \exp( -\psi(z) + A_{c} )$ on $U$
		and $\widetilde{h}_{H}=h_{H}$ on $X\setminus U$. Then $h^{L}:=(\widetilde{h}^{H})^{-1}$ is a Hermitian metric on $L$
		with the desired property.
	\end{pf}

	By Lemma~\ref{lemma:HM}, we can assume that $(L, h^{L})$ is a trivial holomorphic Hermitian line bundle on 
	$\overline{\Omega}$. In what follows, we fix the following isomorphism of holomorphic line bundles over $\overline{\Omega}$:
	\begin{equation}
		\label{eqn:trivial:L}
		(L, h^{L})|_{\overline{\Omega}} \cong ({\mathcal O}_{X}, h^{{\mathcal O}_{X}})|_{\overline{\Omega}}.
	\end{equation}
	
	\medskip
	Let $(F,h^{F})$ be a holomorphic Hermitian line bundle on $X$. Later, we consider the cases $(F,h^{F}) = (L,h^{L})$ and 
	$(F,h^{F}) = ({\mathcal O}_{X}, h^{{\mathcal O}_{X}})$, where $h^{{\mathcal O}_{X}}$ is the trivial metric on ${\mathcal O}_{X}$. Set $F_{s} := F|_{X_{s}}$. 
	For $s\in S^{o}$, let
	$$
	K^{F_{s}}(t,x,x) \sim \frac{a_{0}(x, F_{s})}{t} + a_{1}(x,F_{s}) + O( t )
	$$
	be the asymptotic expansion of the heat kernel of the Laplacain $\square^{F_{s}}$ as $t \to 0$. Then
	\begin{equation}
		\label{eqn:torsion}
		\begin{aligned}
			\log\tau(X_{s},F_{s})
			&=
			\int_{0}^{1} \frac{dt}{t} \int_{X_{s}} \{ K^{F_{s}}(t,x,x) - \frac{a_{0}(x, F_{s})}{t} - a_{1}(x,F_{s}) \} dv_{x}
			\\
			&\quad
			+\int_{1}^{\infty} \frac{dt}{t} \{ \int_{X_{s}} K^{F_{s}}(t,x,x) - h^{0}(F_{s}) \} dv_{x} 
			\\
			&\quad
			- \Gamma'(1) \{ \int_{X_{s}} a_{0}(x,F_{s}) dv_{x} - h^{0}(F_{s}) \}.
		\end{aligned}
	\end{equation}
	Since $g^{X}$ is K\"ahler, $\displaystyle  \int_{X_{s}}a_{0}(x,L_{s}) dv_{x} = \frac{{\rm Area}(X_{s})}{4\pi}$ is independent of $s\in S^{o}$.
	It is classical \cite[Th.\,4.8.16]{Gilkey84} that $\displaystyle  \int_{X_{s}}a_{1}(x,L_{s}) dv_{x}$ is a topological constant independent of $s \in S^{o}$.

	\medskip
	Define the partial analytic torsions of $(X_{s},F_{s})$ by
	\begin{equation}
		\label{eqn:pat:1}
		\log \tau^{\Omega}_{[0,1]}(X_{s},F_{s}) 
		:= 
		\int_{0}^{1} \frac{dt}{t} \int_{\Omega\cap X_{s}} \{ K^{F_{s}}(t,x,x) - \frac{a_{0}(x,F_{s})}{t} - a_{1}(x,F_{s}) \} dv_{x},
	\end{equation}
	\begin{equation}
		\label{eqn:pat:2}
		\log \tau^{X_{s}\setminus\Omega}_{[0,1]}(X_{s},F_{s}) 
		:= 
		\int_{0}^{1} \frac{dt}{t} \int_{X_{s}\setminus\Omega} \{ K^{F_{s}}(t,x,x) - \frac{a_{0}(x,F_{s})}{t} - a_{1}(x,F_{s}) \} dv_{x},
	\end{equation}
	\begin{equation}
		\label{eqn:pat:3}
		\log \tau_{[1,\infty]}(X_{s},F_{s}) 
		:=
		\int_{1}^{\infty} \frac{dt}{t} \{ \int_{X_{s}} K^{F_{s}}(t,x,x) dv_{x} - h^{0}(F_{s}) \}.
	\end{equation}
	Since $L^{-1}$ is an ample line bundle on $X$, we have $h^{0}(L_{s}) =0$ for all $s \in S$.
	Hence for $F=L$ or ${\mathcal O}_{X}$, there is a constant $C_{F}$ independent of $s\in S^{o}$ such that
	\begin{equation}
		\label{eqn:at:vs:pat}
		\log \tau(X_{s},F_{s}) 
		=
		\log \tau^{\Omega}_{[0,1]}(X_{s},F_{s}) 
		+ 
		\log \tau^{X_{s}\setminus\Omega}_{[0,1]}(X_{s},F_{s}) 
		+
		\log \tau_{[1,\infty]}(X_{s},F_{s}) 
		+
		C_{F}.
	\end{equation}

	\subsection{The parameter dependence of the eigenvalues and the heat kernels}
	In this subsection, in order to study the behavior of various partial analytic torsions, we prove the continuity of the eigenvalues of the Hodge-Kodaira Laplacian
	and the heat kernel with respect to the deformation parameter $s\in S^{o}$.
	We identify $(X_{s} \setminus \Omega_{\delta}, g_{s})$ as $(X_{0} \setminus \Omega_{\delta}, \Phi_{s}^{*}g_{s})$ via the diffeomorphism $\Phi_{s}$. 
	
	\begin{lemma}
		\label{lemma:SE}
		There exist a constant $\lambda>0$ independent of $s\in S^{o}$ such that for all $s\in S^{o}$
		$$
		(\square^{L_{s}}f,f) = \| \bar{\partial}f \|_{L^{2}}^{2} \geq \lambda \| f \|_{L^{2}}^{2}
		\qquad
		\forall\,f \in A_{X_{s}}^{0}(L_{s}).
		$$
	\end{lemma}
	
	\begin{pf}
		Since $L^{-1}$ is ample, there is a Hermitian metric $h^{\prime L}$ on $L$ such that $\omega':=-c_{1}(L,h^{\prime L})$ is a K\"ahler form on $X$.
		We write $\| \cdot \|'_{L^{2}}$ for the $L^{2}$-norm with respect to $\omega'$ and $h^{\prime L}$. By the Bochner-Kodaira-Nakano formula,
		we have
		$$
		(\| \bar{\partial}f \|'_{L^{2}})^{2} \geq (\| f \|'_{L^{2}})^{2}
		\qquad
		\forall\,f \in A_{X_{s}}^{0}(L_{s}).
		$$ 
		Recall that $\omega$ is the K\"ahler form of $g$. Since $S$ is compact, there is a constant $C_{0}>0$ such that 
		$C_{0}^{-1}h^{L} \leq h^{\prime L} \leq C_{0} h^{L}$ and $C_{0}^{-1}\omega \leq \omega' \leq C_{0}\omega$ on $X$. Then
		$$
		C_{0}\| \bar{\partial}f \|_{L^{2}}^{2} \geq (\| \bar{\partial}f \|'_{L^{2}})^{2} \geq (\| f \|'_{L^{2}})^{2} \geq C_{0}^{-2} \| f \|_{L^{2}}^{2}
		$$
		for any $f \in A_{X_{s}}^{0}(L_{s})$. We get the result by setting $\lambda=C_{0}^{-3}$.
	\end{pf}

	Let $\{ \phi_{k}^{(s)} \}_{k \in {\mathbf N}}$ be a complete orthonormal system of the Hilbert space of the $L^{2}$-sections of $L_{s}$ consisting of 
	the eigenfunctions of $\square^{L_{s}}$. Let $\lambda_{k}(s)$ be the eigenvalues of $\phi_{k}^{(s)}$. We assume that
	$0< \lambda_{1}(s) \leq \lambda_{2}(s) \leq \cdots$.

	\begin{lemma}
		\label{lemma:asym:EV}
		For all $s \in S^{o}$ and $k \geq 1$, the following inequality holds:
		$$
		\lambda_{k}(s) \geq C k,
		$$
		where $C:=\lambda e^{-\lambda}/\{C_{1}(1,\gamma)e^{\kappa}{\rm Vol}(X_{s})\}$ with $\lambda > 0$ being the same constant as in Lemma~\ref{lemma:SE}.
	\end{lemma}
	
	\begin{pf}
		For all $t\in(0,1]$ and $k\geq1$, it follows from Lemma~\ref{lemma:Gaussian:nabla:HK} (1) that
		$$
		\sum_{i=1}^{k}e^{-t\lambda_{i}(s)}
		\leq
		\sum_{i=1}^{\infty}e^{-t\lambda_{i}(s)}
		=
		{\rm Tr}\,e^{-t\square^{L_{s}}}
		=
		\int_{X_{s}} K^{L_{s}}(t, x,x) dv_{x}
		\leq
		\frac{C_{1}(1,\gamma)e^{\kappa}{\rm Vol}(X_{s})}{t}.
		$$
		Since $\lambda/\lambda_{k}(s)\leq1$ for $k\geq 1$ by Lemma~\ref{lemma:SE}, substituting $t:=\lambda/\lambda_{k}(s)$ in the above inequality, we get
		$$
		k\,e^{-\lambda}
		\leq
		\sum_{i=1}^{k}e^{-\frac{\lambda\lambda_{i}(s)}{\lambda_{k}(s)}}
		\leq
		\frac{C_{1}(1,\gamma)e^{\kappa}{\rm Vol}(X_{s})}{\lambda} \lambda_{k}(s)
		$$
		The result follows from this inequality.
	\end{pf}

	\begin{proposition}
		\label{prop:conv:ev}
		For all $k\geq 1$, $\lambda_{k}(s)$ extends to a continuous function on $S$. Namely, 
		$$
		\lim_{s \to 0} \lambda_{k}(s) = \lambda_{k}(0).
		$$
	\end{proposition}
	
	\begin{pf}
		By \eqref{eqn:comparison:X:FS} and Lemma~\ref{lemma:Gaussian:nabla:HK}, we have the uniform upper bound of the heat kernel of $(X_{s}, g_{s})$.
		Namely, there exists a constant $C>0$ independent of $s \in S^{o}$ such that for all $x, y \in X_{s}$ and $t\in (0,1]$, one has $k_{s}(t, x,y) \leq Ct^{-1}$.
		By \cite[Ths.\,2.1 and 2.16]{CarlenKusuokaStroock87}, there exists a constant $S>0$ independent of $s \in S$ such that for all $f \in C^{\infty}(X_{s})$,
		\begin{equation}
			\label{eqn:Sobolev}
			\| f \|_{L^{4}} \leq S ( \| df \|_{L^{2}} + \| f \|_{L^{2}} ).
		\end{equation}
		Let $K>0$ be a constant such that $-K \omega^{X} \leq c_{1}(L, h^{L}) \leq K \omega^{X}$ on $X$, where $\omega^{X}$ is the K\"ahler form of $(X,g^{X})$. 
		Let $\sigma \in A^{0}(X_{s},L_{s})$. Since $|d|\sigma|| \leq |\nabla^{L}\sigma|$, we have
		\begin{equation}
			\label{eqn:Kato:ineq}
			\| d |\sigma| \|_{L^{2}}^{2} \leq 
			\| \nabla^{L}\sigma \|_{L^{2}}^{2} = ( (\nabla^{L})^{*}\nabla^{L}\sigma, \sigma)_{L^{2}} \leq 2\| \bar{\partial}^{L}\sigma \|_{L^{2}}^{2} + K \|\sigma\|_{L^{2}}^{2},
		\end{equation}
		where the last inequality follows from the Bochner-Kodaira-Nakano formula. By \eqref{eqn:Sobolev}, \eqref{eqn:Kato:ineq}, there exists a constant $S'>0$
		independent of $s \in S^{o}$ such that 
		\begin{equation}
			\label{eqn:unif:Sobolev}
			\| \sigma \|_{L^{4}} \leq S' ( \| \bar{\partial}^{L}\sigma \|_{L^{2}} + \| \sigma \|_{L^{2}} )
		\end{equation}
		for all $\sigma \in A^{0}(X_{s}, L_{s})$. Since \cite[p.114 Conditions (C1), (C2)]{Yoshikawa97} can be verified by using Lemmas~\ref{lemma:H1norm} and
		\ref{lemma:cut:off} below and since we have the uniformity of the Sobolev constant by \eqref{eqn:unif:Sobolev},
		the result can be proved in the same way as \cite[Th.\,5.1]{Yoshikawa97}.
	\end{pf}

	\begin{proposition}
		\label{prop:conv:ef}
		Let $0< \widetilde{\lambda}_{1} < \widetilde{\lambda}_{2} < \cdots$ be the spectrum of $\square^{L_{0}}$.
		Let $\lambda_{k,1}(s) \leq \ldots \leq \lambda_{k,\mu_{k}}(s)$ be the eigenvalues of $\square^{L_{s}}$ converging to $\widetilde{\lambda}_{k}$ as $s \to 0$.
		Let ${\mathfrak K}$ be an arbitrary compact subset of $X_{0}\setminus{\rm Sing}\,X_{0}$. Then the following hold.
		\begin{itemize}
			\item[(1)] 
			For all $k\geq 1$, $\sum_{i=1}^{\mu_{k}}  \Phi_{s}^{*}|\phi_{k,i}^{(s)}|^{2}$ converges to $\sum_{i=1}^{\mu_{k}} |\phi_{k,i}^{(0)}|^{2}$ uniformly on ${\mathfrak K}$ 
			as $s\to 0$.
			\item[(2)]
			$K^{L_{s}}(t, \Phi_{s}(x), \Phi_{s}(x))$ converges to $K^{L_{0}}(t,x,x)$ uniformly on ${\mathfrak K}$ as $s \to 0$.
			\item[(3)]
			$K^{{\mathcal O}_{X_s}}(t, \Phi_{s}(x), \Phi_{s}(x))$ converges to $K^{{\mathcal O}_{X_0}}(t,x,x)$ uniformly on ${\mathfrak K}$ as $s \to 0$.
		\end{itemize}
	\end{proposition}
	
	\begin{pf}
		Since the proof of (3) is completely parallel to that of (2), we only prove (1) and (2). 
		Let $\delta>0$ be such that ${\mathfrak K} \subset X_{0}\setminus\Omega_{\delta}$.
		\par{\bf (1) }
		Let $\{ s_{n} \}_{n\in{\mathbf N}} \subset S^{o}$ be an arbitrary sequence with $\lim_{n\to\infty} s_{n} =0$. 
		By the same argument as in \cite[Prop.\,5.2]{Yoshikawa97}, there exist a subsequence $\{ s_{n(\nu)} \}_{\nu\in{\mathbf N}} \subset \Delta(\epsilon_{\delta})$ 
		and $L^{2}$ sections $\psi_{k,i}$ ($i=1,\ldots,\mu_{k}$) of $L_{0}$ such that $\{ \psi_{k,1}, \ldots, \psi_{k,\mu_{k}} \}$ is an orthonormal basis of the eigenspace 
		$E(\widetilde{\lambda}_{k}, \square^{L_{0}})$ and such that $\Phi_{s_{n(\nu)}}^{*}\phi_{k,i}^{(s_{n(\nu)})}$ converges to $\psi_{k,i}$ in $L^{2}({\mathfrak K}, dv_{0})$.
		Since $\sum_{i=1}^{\mu_{k}} \psi_{k,i}(x) \otimes \langle \cdot, \psi_{k,i}(y)\rangle$ is the integral kernel of the orthogonal projection operator from 
		$L^{2}(X_{0}, L_{0})$ to $E(\widetilde{\lambda}_{k}, \square^{L_{0}})$, we have $\sum_{i=1}^{\mu_{k}} |\psi_{k,i}|^{2} = \sum_{i=1}^{\mu_{k}} |\phi_{k,i}^{(0)}|^{2}$.
		Hence $\sum_{i=1}^{\mu_{k}} \Phi_{s_{n(\nu)}}^{*}|\phi_{k,i}^{(s_{n(\nu)})}|^{2}$ converges to $\sum_{i=1}^{\mu_{k}} |\phi_{k,i}^{(0)}|^{2}$ 
		in $L^{1}({\mathfrak K}, dv_{0})$. Since $\{ s_{n} \}_{n\in{\mathbf N}} \subset S^{o}$ is an arbitrary sequence, this implies that
		$\sum_{i=1}^{\mu_{k}} \Phi_{s}^{*}|\phi_{k,i}^{(s)}|^{2}$ $(s\in\Delta(\epsilon_{\delta}))$ converges to 
		$\sum_{i=1}^{\mu_{k}} |\phi_{k,i}^{(0)}|^{2}$ in $L^{1}({\mathfrak K}, dv_{0})$.
		\par
		Since $\phi_{k,i}^{(s)}$ is a normalized eigenform of $\square^{L_{s}}$ with uniformly bounded eigenvalue $\lambda_{k,i}(s)$ (cf. Proposition~\ref{prop:conv:ev})
		and since $\Phi_{s}^{*}g_{s}$ converges to $g_{0}$ in the $C^{\infty}$-topology on ${\mathfrak K}$, 
		we have $\| \nabla^{\ell} (\Phi_{s}^{*}\phi_{k,i}) \|_{L^{\infty}} \leq C_{k,\ell}$ for $i=1,\ldots,\mu_{k}$ by the elliptic regularity,
		where the constant $C_{k,\ell}>0$ is independent of $s \in \Delta(\delta)$. By Arzel\`a-Ascoli theorem, for any sequence 
		$\{ s_{n} \}_{n\in{\mathbf N}} \subset \Delta(\epsilon_{\delta})$ with $\lim_{n\to\infty} s_{n} =0$, there is a subsequence $\{ s_{n(\nu)} \}_{\nu\in{\mathbf N}}$
		such that $\sum_{i=1}^{\mu_{k}} \Phi_{s_{n(\nu)}}^{*}|\phi_{k,i}^{(s_{n(\nu)})}|^{2}$ converges to $\sum_{i=1}^{\mu_{k}} |\phi_{k,i}^{(0)}|^{2}$ in $C^{0}({\mathfrak K})$.
		Since the limit is independent of the choice of a subsequence, this implies the result.
		\par{\bf (2) }
		Recall that $K^{L_{s}}(t,x,x) = \sum_{m=1}^{\infty} e^{-t\lambda_{m}(s)} | \phi_{m}^{(s)}(x) |^{2}$ for all $t>0$ and $x \in X_{s}$, $s \in S$.
		Since $K^{L_{s}}(t,x,x) \leq C_{1}(1, \gamma)e^{\kappa} (t^{-1}+1)$ for all $0<t\leq 1$ by Lemma~\ref{lemma:Gaussian:nabla:HK} (1), 
		substituting $t=1/\lambda_{m}(s)$, we get
		$$
		e^{-1} | \phi_{m}^{(s)}(x) |^{2} \leq \sum_{j=1}^{\infty} e^{-\lambda_{j}(s)/\lambda_{m}(s)} | \phi_{j}^{(s)}(x) |^{2} \leq C_{1}(1, \gamma)e^{\kappa} (\lambda_{m}(s)+1).
		$$
		Hence
		$$
		\begin{aligned}
			\left| K^{L_{s}}(t,x,x) - \sum_{m = 1}^{N} e^{-t\lambda_{m}(s)} | \phi_{m}^{(s)}(x) |^{2} \right| 
			&= 
			\sum_{m = N+1}^{\infty} e^{-t\lambda_{m}(s)} | \phi_{m}^{(s)}(x) |^{2}
			\\
			&\leq
			C_{1}(1, \gamma)e^{\kappa+1} \sum_{k = N+1}^{\infty} e^{-t\lambda_{k}(s)} (\lambda_{k}(s) + 1).
		\end{aligned}
		$$
		Set $B' = \sup_{x\geq0}xe^{-x/2}$. This, together with Lemma~\ref{lemma:asym:EV}, yields that
		$$
		\begin{aligned}
			\left| K^{L_{s}}(t,x,x) - \sum_{m = 1}^{N} e^{-t\lambda_{m}(s)} | \phi_{m}^{(s)}(x) |^{2} \right| 
			&\leq
			2B'C_{1}(1, \gamma)e^{\kappa+1} \sum_{k = N+1}^{\infty} e^{-tCk/2}
			\\
			&=
			\frac{2B'C_{1}(1, \gamma)e^{\kappa+1}}{ 1 - e^{-Ct/2}} (e^{-Ct/2})^{N+1}.
		\end{aligned}
		$$ 
		Let $\epsilon > 0$ be an arbitrary given number. By this inequality, there exists $N =N(\epsilon, t) \in {\mathbf N}$ such that for all $s \in S$ and $x \in X_{s}$, 
		\begin{equation}
			\label{eqn:diff:HK:Four:ser}
			\left| K^{L_{s}}(t,x,x) - \sum_{m = 1}^{N} e^{-t\lambda_{m}(s)} |\phi_{i}^{(s)}(x)|^{2} \right| < \frac{\epsilon}{3}.
		\end{equation}
		We can assume that $N = \sum_{k=1}^{M} \mu_{k}$. Hence
		$$
		\sum_{m = 1}^{N} e^{-t\lambda_{m}(s)} |\phi_{i}^{(s)}(\Phi_{s}(x))|^{2} 
		=
		\sum_{k=1}^{M} \sum_{i=1}^{\mu_{k}} e^{-t\lambda_{k,i}(s)} |\phi_{k,i}^{(s)}(\Phi_{s}(x))|^{2}.
		$$
		Since for any $x \in {\mathfrak K}$
		\begin{equation}
			\label{eqn:diff:PHK:1}
			\begin{aligned}
				\,&
				\left| \sum_{m = 1}^{N} e^{-t\lambda_{m}(s)} |\phi_{i}^{(s)}(\Phi_{s}(x))|^{2} - \sum_{m = 1}^{N} e^{-t\lambda_{m}(0)} |\phi_{i}^{(s)}(x)|^{2} \right|
				\\
				&\leq
				\sum_{k=1}^{M} \sum_{i=1}^{\mu_{k}} | e^{-t\lambda_{k,i}(s)} - e^{-t\widetilde{\lambda}_{k}} | \cdot |\phi_{k,i}^{(s)}(\Phi_{s}(x))|^{2} 
				\\
				&\quad
				+ \sum_{k=1}^{M} e^{-t\widetilde{\lambda}_{k}} 
				\left| \sum_{i=1}^{\mu_{k}} |\phi_{k,i}^{(s)}(\Phi_{s}(x))|^{2} - \sum_{i=1}^{\mu_{k}} |\phi_{k,i}^{(0)}(x)|^{2} \right|
				\\
				&\leq
				C_{1}(1, \gamma)e^{\kappa+1}\sum_{k=1}^{M} (2\widetilde{\lambda}_{k}+1) \sum_{i=1}^{\mu_{k}} | e^{-t\lambda_{k,i}(s)} - e^{-t\widetilde{\lambda}_{k}} | 
				\\
				&\quad+ \sum_{k=1}^{M} e^{-t\widetilde{\lambda}_{k}} 
				\left\| \sum_{i=1}^{\mu_{k}} \Phi_{s}^{*}|\phi_{k,i}^{(s)}|^{2} - \sum_{i=1}^{\mu_{k}} |\phi_{k,i}^{(0)}|^{2} \right\|_{\mathfrak K},
			\end{aligned}
		\end{equation}
		it follows from Proposition~\ref{prop:conv:ev} and (1) of this proposition and \eqref{eqn:diff:PHK:1} that there exists $\delta=\delta(\epsilon,t)>0$ such that
		for all $x\in X_{s}$, $s \in \Delta(\delta)$,
		\begin{equation}
			\label{eqn:diff:PHK:2}
			\left| \sum_{m = 1}^{N} e^{-t\lambda_{m}(s)} |\phi_{i}^{(s)}(\Phi_{s}(x))|^{2} - \sum_{m = 1}^{N} e^{-t\lambda_{m}(0)} |\phi_{i}^{(s)}(x)|^{2} \right|
			<
			\frac{\epsilon}{3}.
		\end{equation}
		By \eqref{eqn:diff:HK:Four:ser} and \eqref{eqn:diff:PHK:2}, we get $\| \Phi_{s}^{*}K^{L_{s}}(t,\cdot,\cdot) - K^{L_{0}}(t,\cdot,\cdot) \|_{\mathfrak K} < \epsilon$
		for $s \in \Delta(\delta)$. This completes the proof.
	\end{pf}

	\subsection{Estimate for the partial analytic torsion I}
	Recall that $(L,h)|_{\overline{\Omega}}$ is a trivial holomorphic Hermitian line bundle.
	In this subsection, we study the ratio of the partial analytic torsions 
	$\log \tau^{\Omega}_{[0,1]}(X_{s},{\mathcal O}_{X_{s}}) - \log \tau^{\Omega}_{[0,1]}(X_{s}, L_{s})$ as $s \to 0$. 
	\par
	Recall that $\Omega_{\delta} = f^{-1}(\Delta(r_{\delta})) \setminus \Phi( \Delta(r_{\delta}) \times ( X_{0} \setminus B( {\rm Sing}\,X_{0}, \delta) ) )$
	and $\Omega = \Omega_{\delta_{0}}$, $S=\Delta(\delta_{0})$.
	Let $(V_{p},z)$ be a coordinate neighborhood of $p \in {\rm Sing}\,X_{0}$ in $X$. 
	On $V_{p} \cap \Omega_{4\delta_{0}}$, we define $\rho(z) = \| z -p \|^{2}$ and we extend $\rho$ to a smooth function on $X$ in such a way that 
	$\rho \geq 16\delta_{0}^{2}$ on $X \setminus \Omega_{4\delta_{0}}$. Then there exists a constant $A>0$ such that for all $s \in S$,
	\begin{equation}
		\label{eqn:dist:bdry}
		0<A_{s} =\frac{2}{{\rm dist}_{s}( \partial \Omega_{3\delta_{0}/2}\cap X_{s}, \partial\Omega_{2\delta_{0}}\cap X_{s} )} \leq A,
	\end{equation}
	where ${\rm dist}_{s}(\cdot, \cdot)$ is the distance with respect to the metric $g_{s}$.

	\begin{theorem}
		\label{thm:diff:pat}
		The following equality holds:
		$$
		\lim_{s\to 0} \{ \log \tau^{\Omega}_{[0,1]}(X_{s},{\mathcal O}_{X_{s}}) - \log \tau^{\Omega}_{[0,1]}(X_{s}, L_{s}) \} 
		= \log \tau^{\Omega}_{[0,1]}(X_{0},{\mathcal O}_{X_{0}}) - \log \tau^{\Omega}_{[0,1]}(X_{0}, L_{0}).
		$$
	\end{theorem}
	
	\begin{pf}
		Write $k_{s}(t,x,y)$ for the heat kernel of the Hodge-Kodaira Laplacian acting on the sections of trivial Hermitian line bundle ${\mathcal O}_{X_{s}}$ on $X_{s}$.
		Since $(L,h)|_{\Omega}$ is a trivial holomorphic line bundle, we have $a_{i}(x,L) = a_{i}(x, {\mathcal O}_{X_{s}})$ for any $x\in \Omega \cap X_{s}$. 
		We apply Theorem~\ref{thm:diff:HK} by setting $M= X_{s}$, $L=L_{s}$, $\Omega_{c} = X_{s} \cap \Omega_{c}$. 
		By \eqref{eqn:comparison:X:FS}, \eqref{eqn:dist:bdry}, we have the uniformity of the constants $\gamma$ and $A$ in Theorem~\ref{thm:diff:HK}
		with respect to $s \in S$. Namely, we can take $\gamma$ and $A$ independent of $s \in S$ in Theorem~\ref{thm:diff:HK}. 
		Then there exists a constant $D>0$ independent of $s \in S^{o}$ such that for all $x,y \in \Omega \cap X_{s}$, 
		$$
		\left| k_{s}(t, x, y)- K^{L_{s}}(t, x, y) \right| 
		\leq
		D\,\rho_{s}(x,y)^{-2(2n+1)}e^{-\frac{\rho_{s}(x,y)^{2}}{16t}},
		$$
		where $\rho_{s}(x,y) = \min\{ d_{s}(x, \partial\Omega_{3/2}\cap X_{s}), d_{s}(y, \partial\Omega_{3/2}\cap X_{s}) \}$, 
		$d_{s}(\cdot,\cdot)$ being the distance function on $(X_{s}, g_{s})$.
		Set $\rho = \min_{s\in S^{o}}\min_{x\in X_{s}\cap \Omega_{1}} d_{s}(x, \partial\Omega_{3/2}\cap X_{s})>0$. 
		Then for all $x,y \in \Omega\cap X_{s}$, $s \in S^{o}$, we have
		\begin{equation}
			\label{eqn:est:H:K:X:s}
			\left| k_{s}(t, x, y)- K^{L_{s}}(t, x, y) \right| 
			\leq
			D\,\rho^{-2(2n+1)}e^{-\frac{\rho^{2}}{16t}}.
		\end{equation}
		We remark that \eqref{eqn:est:H:K:X:s} holds also for $X_{0}$ with possibly different positive constants $D$, $\rho$. 
		By \eqref{eqn:est:H:K:X:s}, there exists a constant $C(\rho)>0$ depending only on $\rho>0$ such that 
		for all $0<\delta<\delta_{0}$ and $s\in\Delta(r_{\delta})$
		\begin{equation}
			\label{eqn:diff:HK:O:L:0}
			\int_{0}^{1} \frac{dt}{t} \int_{\Omega_{\delta}\cap X_{s}} | k_{s}(t,x,x) - K^{L_{s}}(t,x,x) |\,dv_{s}(x) \leq C(\rho){\rm Area}(\Omega_{\delta}\cap X_{s}).
		\end{equation}
		\par
		Let $\epsilon > 0$ be an arbitrary number. We take $0<\delta<\delta_{0}$ in such a way that ${\rm Area}(\Omega_{\delta}\cap X_{s}) < \epsilon/2C(\rho)$
		for all $s \in \Delta(r_{\delta})$. By \eqref{eqn:diff:HK:O:L:0}, we get
		\begin{equation}
			\label{eqn:diff:HK:O:L:1}
			\int_{0}^{1} \frac{dt}{t} \int_{\Omega_{\delta}\cap X_{s}} | k_{s}(t,x,x) - K^{L_{s}}(t,x,x) |\,dv_{s}(x) \leq \frac{\epsilon}{2}.
		\end{equation}
		By \eqref{eqn:est:H:K:X:s}, Proposition~\ref{prop:conv:ef} (2), (3) and Lebesgue's convergence theorem, there exists $r'>0$ such that for all $s \in \Delta(r')$, 
		\begin{equation}
			\label{eqn:diff:HK:O:L:2}
			\begin{aligned}
				\,&
				\left| 
				\int_{0}^{1} \frac{dt}{t} \left\{ \int_{X_{s}\cap(\Omega\setminus\Omega_{\delta})} \{ k_{s}(t,\cdot,\cdot) - K^{L_{s}}(t,\cdot,\cdot) \}\,dv_{s}
				-
				\int_{X_{0}\cap(\Omega\setminus\Omega_{\delta})} \{ k_{0}(t,\cdot,\cdot) - K^{L_{0}}(t,\cdot,\cdot) \}\,dv_{0} \right\}
				\right|
				\\
				&\quad< 
				\frac{\epsilon}{2}.
			\end{aligned}
		\end{equation}
		By \eqref{eqn:trivial:L}, \eqref{eqn:pat:1}, 
		\begin{equation}
			\label{eqn:diff:H:K:X:s}
			\begin{aligned}
				\,&
				\log \tau^{\Omega}_{[0,1]}(X_{s},{\mathcal O}_{X_{s}}) - \log \tau^{\Omega}_{[0,1]}(X_{s}, L_{s})
				\\
				&=
				\int_{0}^{1} \frac{dt}{t} \int_{\Omega_{\delta}\cap X_{s}} \{ k_{s}(t,x,x) - K^{L_{s}}(t,x,x) \}\,dv_{s}(x)
				\\
				&\quad+
				\int_{0}^{1} \frac{dt}{t} \int_{X_{s} \cap (\Omega\setminus\Omega_{\delta})} \{ k_{s}(t,x,x) - K^{L_{s}}(t,x,x) \}\,dv_{s}(x).
			\end{aligned}
		\end{equation}
		We deduce from \eqref{eqn:diff:HK:O:L:1}, \eqref{eqn:diff:HK:O:L:2}, \eqref{eqn:diff:H:K:X:s} that 
		$$
		| \{ \log \tau^{\Omega}_{[0,1]}(X_{s},{\mathcal O}_{X_{s}}) - \log \tau^{\Omega}_{[0,1]}(X_{s}, L_{s}) \}
		-
		\{ \log \tau^{\Omega}_{[0,1]}(X_{0},{\mathcal O}_{X_{0}}) - \log \tau^{\Omega}_{[0,1]}(X_{0}, L_{0}) \} | < \epsilon
		$$
		for all $s \in \Delta(r'')$, $r''=\min\{ r_{\delta}, r' \}$.
		This proves the result.
	\end{pf}

	\subsection{Estimate for the partial analytic torsion II}
	In this subsection, we study the asymptotic behavior of $\tau^{X_{s}\setminus\Omega}_{[0,1]}(X_{s},{\mathcal O}_{X_{s}})$ and 
	$\tau^{X_{s}\setminus\Omega}_{[0,1]}(X_{s}, L_{s})$ as $s \to 0$.

	\begin{theorem}
		\label{thm:bd:pat}
		The following equalities hold:
		$$
		\lim_{s\to0} \log \tau^{X_{s}\setminus\Omega}_{[0,1]}(X_{s},{\mathcal O}_{X_{s}}) = \log \tau^{X_{0}\setminus\Omega}_{[0,1]}(X_{0},{\mathcal O}_{X_{0}}), 
		$$
		$$
		\lim_{s\to0} \log \tau^{X_{s}\setminus\Omega}_{[0,1]}(X_{s}, L_{s}) = \log \tau^{X_{0}\setminus\Omega}_{[0,1]}(X_{0}, L_{0}).
		$$
	\end{theorem}
	
	\begin{pf}
		We only prove the second equality, since the proof of the first one is similar. 
		We regard $(X_{s} \setminus \Omega, g_{s})$ as $(X_{0} \setminus \Omega, \Phi_{s}^{*}g_{s})$ via the diffeomorphism $\Phi_{s}$. 
		Since $\{\Phi_{s}^{*}g_{s}\}_{s\in S}$ is a family of Riemannian metrics
		on $X_{0} \setminus \Omega$ depending smoothly in $s$, by shrinking $S$ if necessary, there exists $j>0$ such that $j_{x} \geq j$ for all 
		$x \in X_{s} \setminus \Omega$, $s \in S$ and such that on the ball ${\mathbb B}(x, 3j)$ endowed with the geodesic normal coordinates centered at $x$, 
		the metric tensor and its higher derivatives up to order $k(>5)$ are uniformly bounded for all $x \in X_{s} \setminus \Omega$, $s \in S$. 
		By the formula for the constsnt $\widetilde{D}_{k}(y)$ in Theorem~\ref{thm:asym:exp:H:K} and this uniformity, 
		there exists a constant $\widetilde{D}_{k}>0$ such that $\widetilde{D}_{k}(y) \leq \widetilde{D}_{k}$ for all $y \in X_{s} \setminus \Omega$, $s\in S$.
		Namely, we have
		\begin{equation}
			\label{eqn:asym:exp:H:K:X:s}
			\left| K^{L_{s}}(t,y,y)- \left( \frac{a_{0}(y,L_{s})}{t} + a_{1}(y, L_{s}) + \cdots + a_{k}(y, L_{s})t^{k-1} \right) \right| \leq \widetilde{D}_{k}\,t^{k}
		\end{equation}
		for all $y \in X_{s} \setminus \Omega$, $s\in S$ and $t\in (0,1]$.
		\par
		By Proposition~\ref{prop:conv:ef} (2), $\Phi_{s}^{*}K^{L_{s}}(t,y,y)$ converges to $K^{L_{0}}(t,y,y)$ uniformly on $X_{0}\setminus \Omega$ as $s \to 0$. 
		Since $\Phi_{s}^{*}g_{s}$ converges to $g_{0}$ in the $C^{\infty}$-topology of $X_{0}\setminus \Omega$ as $s \to 0$, we see that $a_{i}(\Phi_{s}(y), L_{s})$
		converges to $a_{i}(y, L_{0})$ uniformly on $X_{0}\setminus \Omega$ as $s \to 0$.
		Hence, by \eqref{eqn:asym:exp:H:K:X:s} and Lebesgue's convergence theorem applied to the integral
		$$
		\int_{0}^{1} \frac{dt}{t} \int_{X_{s}\setminus\Omega}
		\left\{ K^{L_{s}}(t,y,y)- \left( \frac{a_{0}(y,L_{s})}{t} + a_{1}(y, L_{s}) + \cdots + a_{k}(y, L_{s})t^{k-1} \right) \right\} dv_{s}(y),
		$$
		we get
		$$
		\lim_{s \to 0} \log \tau^{X_{s}\setminus\Omega}_{[0,1]}(X_{s}, L_{s}) =
		\int_{0}^{1} \frac{dt}{t} \int_{X_{0}\setminus\Omega} \left\{ K^{L_{0}}(t,y,y)- \left( \frac{a_{0}(y,L_{0})}{t} + a_{1}(y, L_{0}) \right) \right\} dv_{0}(y).
		$$
		This completes the proof.
	\end{pf}

	\subsection{Small eigenvalues}
	
	Recall that $k_{s}(t,x,y)$ is the heat kernel of the Laplacain acting on the functions on $X_{s}$. 
	Since $L^{-1}$ is an ample line bundle on $X$, $H^{0}(X_{s}, L_{s})=0$ for all $s\in S^{o}$.
	The partial analytic torsions $\tau_{[1,\infty]}(X_{s},{\mathcal O}_{X_{s}})$ and $\tau_{[1,\infty]}(X_{s},L_{s})$ are given respectively by
	$$
	\log \tau_{[1,\infty]}(X_{s},{\mathcal O}_{X_{s}}) 
	=
	\int_{1}^{\infty} \frac{dt}{t} \{\int_{X_{s}} k(t,x,x) dv_{s}(x) - 1 \},
	$$
	$$
	\log \tau_{[1,\infty]}(X_{s},L_{s}) 
	=
	\int_{1}^{\infty} \frac{dt}{t} \int_{X_{s}} K^{L_{s}}(t,x,x) dv_{s}(x).
	$$
	Recall that
	$$
	N = \dim H^{0}(X_{0}\setminus {\rm Sing}\,X_{0}, {\mathbf C}) = \#\{ \hbox{irreducible components of } X_{0}\}. 
	$$

	\begin{theorem}
		\label{thm:SEF}
		The function $\lambda_{k}$ on $S^{o}$ extends to a continuous function on $S$. In particular, $\lambda_{k}(s)\to0$ as $s\to0$ for $k\leq N-1$. 
		Moreover, there exists $\lambda>0$ such that for all $k\geq N$,
		$$
		\lambda_{k}(s) \geq \lambda.
		$$
	\end{theorem}
	
	\begin{pf}
		See \cite[Th.\,A]{JiWentworth92} and \cite[Main Th.]{Yoshikawa97}.
	\end{pf}

	\begin{theorem}
		\label{thm:pat:smev}
		As $s \to 0$,
		$$
		\log \tau_{[1,\infty]}(X_{s},{\mathcal O}_{X_{s}}) = -\log\{\prod_{i=1}^{N-1} \lambda_{i}(s)\} + \log \tau_{[1,\infty]}(X_{0},{\mathcal O}_{X_{0}}) + c + o(1),
		$$
		where $\displaystyle c = (N-1) \{ \int_{1}^{\infty} e^{-t} \frac{dt}{t} + \int_{0}^{1} (e^{-t}-1) \frac{dt}{t} \}$.
	\end{theorem}
	
	\begin{pf}
		Following Chen-Li \cite[Th.\,1]{ChenLi81}, we derive a lower bound of $\lambda_{k}(s)$.
		By \cite[Cor.\,4.2]{Yoshikawa97}, there is a constant $A>0$ such that 
		$$
		\| f \|_{L^{4}} \leq A( \| df \|_{L^{2}} + \| f \|_{L^{2}} ) 
		\qquad
		\forall\,f\in C^{\infty}(X_{s}), \quad s\in S^{o}.
		$$
		By \cite[Ths.\,2.1 and 2.16]{CarlenKusuokaStroock87}, there exists a constant $C>0$ such that $k(t,x,y) \leq C\,t^{-2}$ for all $t\in(0,1]$, $x,y\in X_{s}$, $s\in S^{o}$. 
		Hence for all $t\in(0,1]$ and $k\geq1$,
		$$
		\sum_{i=1}^{k}e^{-t\lambda_{i}(s)}
		\leq
		\sum_{i=1}^{\infty}e^{-t\lambda_{i}(s)}
		=
		{\rm Tr}\,e^{-t\square_{s}}
		\leq
		C{\rm Vol}(X_{s}) \,t^{-2}.
		$$
		Let $k\geq N$. Since $\lambda/\lambda_{i}(s)\leq1$ for $i\geq N$ by Theorem~\ref{thm:SEF}, 
		substituting $t:=\lambda/\lambda_{k}(s)$
		in the above inequality and using $\lambda_{i}(s)/\lambda_{k}(s)\leq1$ for $i\leq k$, we get
		$$
		\left( k-(N-1) \right)\,e^{-\lambda}
		\leq
		\sum_{i=N}^{k}e^{-\frac{\lambda\lambda_{i}(s)}{\lambda_{k}(s)}}
		\leq
		C{\rm Vol}(X_{s})\,\left(\frac{\lambda}{\lambda_{k}(s)}\right)^{-2}
		=
		C{\rm Vol}(X_{s})\,\left(\frac{\lambda_{k}(s)}{\lambda}\right)^{2}.
		$$
		We set $B := \lambda\sqrt{k e^{-\lambda}/\{C{\rm Vol}(X_{s})\}}$. Then we get for all $k\geq N$ and $s\in S^{o}$
		$$
		\lambda_{k}(s) \geq B \sqrt{k-(N-1)}.
		$$
		Since $\sum_{i=N}^{\infty}e^{-tB\sqrt{k-(N-1)}}/t$ is an integrable function on $[1, \infty)$ dominating the function 
		$\sum_{i=N}^{\infty}e^{-t\lambda_{i}(s)}/t$, we get
		$$
		\log \tau_{[1,\infty]}(X_{s},{\mathcal O}_{s}) - \sum_{i=1}^{N-1} \int_{\lambda_{i}(s)}^{\infty} e^{-t} \frac{dt}{t}
		=
		\int_{1}^{\infty} \sum_{i=N}^{\infty}e^{-t\lambda_{i}(s)} \frac{dt}{t} 
		=
		\int_{1}^{\infty} \sum_{i=N}^{\infty}e^{-t\lambda_{i}(0)} \frac{dt}{t} + o(1)
		$$
		as $s\to0$.
		This, together with
		$$
		\sum_{i=1}^{N-1} \int_{\lambda_{i}(s)}^{\infty} e^{-t} \frac{dt}{t} 
		= 
		-\sum_{i=1}^{N-1} \log \lambda_{i}(s) + \sum_{i=1}^{N-1} \int_{1}^{\infty} e^{-t} \frac{dt}{t} + \sum_{i=1}^{N-1} \int_{0}^{1} (e^{-t}-1) \frac{dt}{t} + o(1)
		$$
		implies the result.
	\end{pf}

	\begin{theorem}
		\label{thm:pat:bd2}
		The following equality holds as $s \to 0$:
		$$
		\log \tau_{[1,\infty]}(X_{s},L_{s}) = \log \tau_{[1,\infty]}(X_{0},L_{0}) + o(1)
		\quad
		(s\to0).
		$$
	\end{theorem}
	
	\begin{pf}
		Let $\Delta_{L_{s}}=(\nabla^{L_{s}})^{*}\nabla^{L_{s}}$ be the Bochner Laplacian acting on $A^{0}_{X_{s}}(L_{s})$, where 
		$\nabla^{L_{s}}= \partial_{L_{s}}+\bar{\partial}$ is the Chern connection of $(L_{s},h_{L_{s}})$. 
		Let $R_{L_{s}}$ be the curvature of $(L_{s},h_{L_{s}})$.
		Set $Q_{s}:=i\Lambda_{s} R_{L}|_{X_{s}}$. Since $(L,h_{L})$ is a semi-negative line bundle, $Q_{s}\leq 0$ on $X_{s}$.
		Since $\partial_{L_{s}}^{*}\partial_{L_{s}} - \bar{\partial}^{*}\bar{\partial} = i\Lambda R_{L_{s}}$ 
		by the Bochner-Kodaira-Nakano formula, we have $\Delta_{L_{s}}=2\square_{L_{s}}+Q_{s}$. 
		Since $2\square_{L_{s}}=\Delta_{L_{s}}-Q_{s}$ and $-Q_{s}\geq0$, we get by \cite{HessSchraderUhlenbrock80}
		$$
		| K^{L_{s}}(2t,x,y) | \leq k(t,x,y)
		\qquad
		(\forall\,t>0,\quad\forall\,x,y\in X_{s}).
		$$
		\par
		By \cite[Cor.\,4.2]{Yoshikawa97}, there is a constant $A>0$ such that 
		$\| f \|_{L^{4}} \leq A( \| df \|_{L^{2}} + \| f \|_{L^{2}} )$ for all $f\in C^{\infty}(X_{s})$ and $s\in S^{o}$.
		Then there exists a constant $C>0$ such that $k(t,x,y) \leq C\,t^{-2}$ for all $t\in(0,1]$, $x,y\in X_{s}$, $s\in S^{o}$. 
		Hence for any $t\in(0,1]$,
		$$
		{\rm Tr}\,e^{-t\square_{L_{s}}}
		\leq
		C{\rm Area}(X_{s}) \,t^{-2}.
		$$
		Let $0<\lambda_{1}^{L}(s) \leq \lambda_{2}^{L}(s) \leq\cdots$ be the eigenvalues of $\square_{L_{s}}$. For any $m \geq 1$,
		$$
		\sum_{i=1}^{m}e^{-t\lambda_{i}^{L}(s)}
		\leq
		\sum_{i=1}^{\infty}e^{-t\lambda_{i}^{L}(s)}
		=
		{\rm Tr}\,e^{-t\square_{L_{s}}}
		\leq
		C{\rm Area}(X_{s}) \,t^{-2}.
		$$
		Since $\lambda/\lambda_{m}^{L}(s)\leq1$ by Lemma~\ref{lemma:SE}, substituting $t:=\lambda/\lambda_{m}^{L}(s)$
		in the above inequality and using $\lambda_{i}^{L}(s)/\lambda_{m}^{L}(s)\leq1$ for $i\leq m$, we get
		$$
		m\,e^{-\lambda}
		\leq
		\sum_{i=1}^{m}e^{-\frac{\lambda\lambda_{i}^{L}(s)}{\lambda_{m}^{L}(s)}}
		\leq
		C{\rm Area}(X_{s})\,\left(\frac{\lambda}{\lambda_{m}^{L}(s)}\right)^{-2}
		=
		C{\rm Area}(X_{s})\,\left(\frac{\lambda_{m}^{L}(s)}{\lambda}\right)^{2}
		$$
		We set $B := \lambda\sqrt{k e^{-\lambda}/\{C{\rm Area}(X_{s})\}}$. Then we get for all $m\geq1$ and $s\in S^{o}$
		$$
		\lambda_{m}^{L}(s) \geq B m^{1/2}.
		$$
		Since $\sum_{m=1}^{\infty}e^{-tB\sqrt{m}}/t \in L^{1}([1,\infty))$ dominates $\sum_{m=1}^{\infty}e^{-t\lambda_{m}^{L}(s)}/t$, we get
		$$
		\log \tau_{[1,\infty]}(X_{s},L_{s}) 
		=
		\int_{1}^{\infty} \sum_{m=1}^{\infty}e^{-t\lambda_{m}^{L}(s)} \frac{dt}{t} 
		=
		\int_{1}^{\infty} \sum_{m=1}^{\infty}e^{-t\lambda_{m}^{L}(0)} \frac{dt}{t} + o(1)
		$$
		as $s \to 0$. This completes the proof.
	\end{pf}

	\begin{theorem}
		\label{thm:at:bd1}
		The following equality holds as $s \to 0$:
		$$
		\log \tau(X_{s},{\mathcal O}_{X_{s}}) - \log \tau(X_{s},L_{s}) 
		= 
		-\log\{\prod_{i=1}^{N-1} \lambda_{i}(s)\}  + c + o(1)
		$$
		with $c = \log \tau(X_{0},{\mathcal O}_{X_{0}}) - \log \tau(X_{0},L_{0}) + (N-1) ( \int_{1}^{\infty} e^{-t} \frac{dt}{t} + \int_{0}^{1} (e^{-t}-1) \frac{dt}{t} )$.
	\end{theorem}
	
	\begin{pf}
		The result follows from Theorems~\ref{thm:diff:pat}, \ref{thm:bd:pat}, \ref{thm:pat:smev}, \ref{thm:pat:bd2}.
	\end{pf}

	\section{Quillen metrics and the ratio of analytic torsions}
	\label{sect:4}

	We maintain the notation from Section~\ref{sect:3}.
	In this section, we give another expression of 
	$\log ( \tau(X_{s},{\mathcal O}_{X_{s}}) / \tau(X_{s},L_{s}))$ as $s \to 0$ in terms of certain period integrals to prove Theorem~\ref{Main:Thm:2}.
	To this end, we use the notion of Quillen metrics, for which we refer the reader to \cite{BGS88}, \cite{BismutBost90}. 
	In Sections~\ref{sect:4} and \ref{sect:5}, we assume that the line bundle $H=L^{-1}$ is defined by an effective divisor on $X$ whose support is disjoint from ${\rm Sing}\,X_{0}$ and intersects each irreducible component of $X_{0}$ transversally at a single point.
	Consequently, $\deg H_{s} = -\deg L_{s} = N$ for $s\in S$.
	See Subsection~\ref{sect:4:2} for the construction of $H$. We note that this assumption on $H$ can always be satisfied in our situation since our family is over the unit disc with a reduced central fiber. (If the base were a compact Riemann surface, such an ample line bundle would not exist in general.)

	\subsection{Semi-stable reduction}
	To give an expression of $\log ( \tau(X_{s},{\mathcal O}_{X_{s}}) / \tau(X_{s},L_{s}) )$ in terms of certain period integrals,
	we consider a semistable reduction of the family $f \colon (X, X_{0}) \to (S,0)$, which consists of the following commutative diagram:
	$$
	\begin{CD}
		(Y,Y_{0}=\psi^{-1}(0))@>F>> (X,X_{0})
		\\
		@V \widetilde{f} VV  @V f VV
		\\
		(T,0) @> \mu >> (S,0).
	\end{CD}
	$$
	Here $(T,0)$ is another unit disc of ${\mathbf C}$, $\mu \colon (T,0) \to (S,0)$ is given by $\mu(t) = t^{\nu}$ for some $\nu \in {\mathbf N}$, 
	$Y$ is a smooth complex surface such that $Y\setminus Y_{0} \cong X \times_{S\setminus\{0\}} (T\setminus\{0\})$ is the family induced from 
	$f \colon X\setminus X_{0} \to S\setminus\{0\}$ by $\mu$, $Y_{0}$ is a {\em reduced} normal crossing divisor of $Y$, and $F \colon Y \to X$ is the composition
	of the projection $X \times_{T}S \to X$ and a holomorphic map $Y \to X \times_{T}S$, which is a sequence of blowing-ups. 
	In this section, contrary to the preceding sections, $t$ is a holomorphic coordinate of $T$ centered at $0$.
	We set $Y_{t} := \widetilde{f}^{-1}(t)$. Then $Y_{t} \cong X_{\mu(t)} = X_{t^{\nu}}$ for $t\not=0$. 
	Recall that $X_{0} = C_{0}+C_{1}+\cdots+C_{N-1}$ is the irreducible decomposition of $X_{0}$. Then we have
	$$
	Y_{0} = \widetilde{C}_{0}+\cdots+\widetilde{C}_{N-1}+E_{1}+\cdots+E_{m},
	$$ 
	where $F(\widetilde{C}_{i})=C_{i}$ and $F(E_{j})$ is a singular point of $X_{0}$. 
	Since $Y$ is obtained from $X \times_{T}S$ by a sequence of blowing-ups, $Y$ is K\"ahler. 
	\par
	We consider the following two determinants of the cohomology:
	$$
	\lambda({\mathcal O}_{Y}) = \det R\widetilde{f}_{*}{\mathcal O}_{Y} = \widetilde{f}_{*}{\mathcal O}_{Y} \otimes (\det R^{1}\widetilde{f}_{*}{\mathcal O}_{Y} )^{\lor}
	= \widetilde{f}_{*}{\mathcal O}_{Y} \otimes \det \widetilde{f}_{*} K_{Y/T},
	$$
	$$
	\lambda(F^{*}L) = \det R\widetilde{f}_{*}(F^{*}L) = (\det R^{1}\widetilde{f}_{*}F^{*}L)^{\lor} = \det \widetilde{f}_{*} K_{Y/T}(F^{*}L^{-1}),
	$$
	where $K_{Y/T} := K_{Y}\otimes \widetilde{f}^{*}K_{T}^{-1}$ is the relative canonical bundle of the family $\widetilde{f} \colon Y \to T$. 
	Recall that $H=L^{-1}$.

	\begin{lemma}
		\label{lemma:dim:coh}
		For all $t\in T$, one has $h^{0}(Y_{t}, K_{Y_{t}}(F^{*}H_{t})) = g+N-1$.
	\end{lemma}
	
	\begin{pf}
		Since the family $\widetilde{f} \colon Y \to T$ is flat, it suffices by Riemann-Roch theorem to prove that $H^{1}(Y_{t}, K_{Y_{t}}(F^{*}H_{t}))=0$ for all $t \in T$.
		For $t\not=0$, $H^{1}(Y_{t}, K_{Y_{t}}(F^{*}H_{t}))^{\lor} = H^{0}(Y_{t}, F^{*}L_{t}) =0$ because $L$ is a negative line bundle on $X$.
		Let $\sigma \in H^{0}(Y_{0}, F^{*}L)$. Then $\sigma|_{C_{i}} = 0$ since $(F^{*}L)|_{C_{i}}$ is a negative line bundle on $C_{i}$.
		Hence, if $E_{j}\cap C_{i} \not= \emptyset$, then $\sigma|_{E_{j}}$ has zeros. Since $(F^{*}L)|_{E_{j}}$ is a trivial line bundle on $E_{j}$,
		this implies $\sigma|_{E_{j}} = 0$ if $E_{j} \cap C_{i} \not= \emptyset$ for some $C_{i}$. In the same way, if $E_{j} \cap E_{k} \not= \emptyset$ and
		$\sigma|_{E_{j}} = 0$, then $\sigma|_{E_{k}} = 0$. Since $Y_{0}$ is connected, we conclude $\sigma=0$. 
		This proves that $H^{1}(Y_{0}, K_{Y_{0}}(F^{*}H_{0}))=0$. 
	\end{pf}

	\subsection{The $L^{2}$-metric on the determinant of the cohomology}
	\label{sect:4:2}
	\par
	Let $g^{Y}$ be a K\"ahler metric on $Y$. We also consider the degenerate K\"ahler metric $F^{*}g^{X}$ on $Y$, 
	which is a genuine K\"ahler metric on $Y \setminus Y_{0}$. 
	Then $\lambda(F^{*}L)|_{T^{o}}$ is endowed with the $L^{2}$-metric $\| \cdot \|_{L^{2},\lambda(F^{*}L)}$ (resp. $\| \cdot \|'_{L^{2},\lambda(F^{*}L)}$) 
	with respect to $g^{Y}$, $F^{*}h^{L}$ (resp. $F^{*}g^{X}$, $F^{*}h^{L}$). 
	Similarly, $\lambda({\mathcal O}_{Y})|_{T^{o}}$ is endowed with the $L^{2}$-metric $\| \cdot \|_{L^{2},\lambda({\mathcal O}_{Y})}$ 
	(resp. $\| \cdot \|'_{L^{2},\lambda({\mathcal O}_{Y})}$) with respect to $g^{Y}$, $F^{*}h^{L}$ (resp. $F^{*}g^{X}$, $F^{*}h^{L}$). 
	\par
	Let $\omega_{1}, \ldots, \omega_{g}$ be a basis of $\widetilde{f}_{*}K_{Y/T}$ as a free ${\mathcal O}_{T}$-module near $0$.
	Then $ \lambda({\mathcal O}_{Y}) =  \widetilde{f}_{*}{\mathcal O}_{Y} \otimes \det \widetilde{f}_{*} K_{Y/T}$ is generated by
	$$
	\sigma := 1 \otimes (\omega_{1} \wedge \ldots \wedge \omega_{g}).
	$$
	We set $\omega_{i}(t) := \omega_{i}|_{Y_{t}}$. 
	Let $A_{1} = {\rm Area}(X_{s}, g_{s})$ for $s\not=0$ and $A_{2} = {\rm Area}(Y_{t}, g^{Y}|_{Y_{t}})$ for $t\not=0$. By definition of the $L^{2}$-metrics, we have
	\begin{equation}
		\label{eqn:L2:met:1}
		\| \sigma \|_{L^{2},\lambda({\mathcal O}_{Y})}^{2} (t) = A_{2} \det\left( \frac{i}{2} \int_{Y_{t}} \omega_{l}\wedge \overline{\omega}_{m} \right)_{1\leq l,m \leq g},
	\end{equation}
	\begin{equation}
		\label{eqn:L2:met:2}
		\| \sigma \|_{L^{2},\lambda({\mathcal O}_{Y})}^{\prime\,2} (t) 
		= A_{1} \det\left( \frac{i}{2} \int_{Y_{t}} \omega_{l}\wedge \overline{\omega}_{m} \right)_{1\leq l,m \leq g}.
	\end{equation}
	
	\par
	Let $p_{i} \colon S \to X$ $(0\leq i \leq N-1)$ be a section such that $p_{i}(s) \not= p_{j}(s)$ for $i\not=j$ and such that 
	$p_{i}(0) \in C_{i} \setminus {\rm Sing}\,X_{0}$ for all $i$.
	Then $\sum_{i=0}^{N-1} p_{i}$ is an ample divisor of $X$, which does not meet ${\rm Sing}\,X_{0}$. We define 
	$$
	H := {\mathcal O}_{X}(\sum_{i=0}^{N-1} p_{i} ).
	$$
	Since $F \colon Y \setminus F^{-1}({\rm Sing}\,X_{0}) \to X \setminus {\rm Sing}\,X_{0}$ is an isomorphism, $F^{-1}\circ p_{i}$ is a divisor on $Y$. 
	We set $\widetilde{p}_{i} := F^{-1}\circ p_{i}$. Then 
	$$
	F^{*}H = {\mathcal O}_{Y}( \sum_{i=0}^{N-1} \widetilde{p}_{i}).
	$$
	Since $F^{*}H_{t}=F^{*}H|_{Y_{t}} = {\mathcal O}_{Y_{t}}(\sum_{i=0}^{N-1}\widetilde{p}_{i}(t))$ with 
	$\widetilde{p}_{i}(0) \in \widetilde{C}_{i}\setminus F^{-1}({\rm Sing}\,X_{0})$, 
	an element of $H^{0}(Y_{t}, K_{Y_{t}}(F^{*}H_{t}))$ is viewed as a meromorphic Abelian differential with at most logarithmic poles on 
	$\sum_{i=0}^{N-1} \widetilde{p}_{i}(t)$. In particular, $H^{0}(Y_{t}, K_{Y_{t}}) \subset H^{0}(Y_{t}, K_{Y_{t}}(H_{t}))$. 
	Since ${\mathcal O}_{Y}(K_{Y}) \subset {\mathcal O}_{Y}(K_{Y}(F^{*}H))$, $\omega_{1}, \ldots, \omega_{g}$ are local sections of $\widetilde{f}_{*}K_{Y/T}(F^{*}H)$.
	Let $\omega_{g+1}, \ldots, \omega_{g+N-1}$ be local sections of $\widetilde{f}_{*}K_{Y/T}(F^{*}H)$ near $0\in T$ such that
	$\{ \omega_{1}, \ldots, \omega_{g+N-1} \}$ is a basis of $\widetilde{f}_{*}K_{Y/T}(F^{*}H)$ as a free ${\mathcal O}_{T}$-module near $0$.
	Shrinking $T$ if necessary, we can assume that $\omega_{i} \in H^{0}(Y, K_{Y/T})$ $(1\leq i \leq g)$ and
	$\omega_{j} \in H^{0}(Y, K_{Y/T}(F^{*}H))$ $(1 \leq j \leq g+N-1)$.  
	By Lemma~\ref{lemma:dim:coh} and Grauert's base change theorem,
	$\{ \omega_{1}(t), \ldots, \omega_{g+N-1}(t) \}$ is a basis of $H^{0}(Y_{t}, K_{Y_{t}}(F^{*}H))$ with
	$\omega_{i}(t) \in H^{0}(Y_{t}, K_{Y_{t}})$ $(1\leq i \leq g)$ and $\omega_{j}(t) \in H^{0}(Y_{t}, K_{Y_{t}}(F^{*}H))$ $(g+1 \leq j \leq g+N-1)$ for all $t \in T$.
	Since $H^{0}(Y_{t}, K_{Y_{t}}(\widetilde{p}_{i}(t) + \widetilde{p}_{j}(t) )) \not= H^{0}(Y_{t}, K_{Y_{t}})$ for all $t \in T$ by Riemann-Roch,
	we can choose $\omega_{g+i}(t)$ $(1 \leq i \leq N-1)$ in such a way that 
	the only poles of $\omega_{g+i}(t)$ are $\widetilde{p}_{0}(t)$ and $\widetilde{p}_{i}(t)$ for all $t \in T$. We set
	$$
	\widetilde{\sigma} := \omega_{1} \wedge\cdots\wedge \omega_{g+N-1}.
	$$
	Let $h^{H}$ be the Hermitian metric on $H=L^{-1}$ induced from $h^{L}$.
	By definition of the $L^{2}$-metrics, we have
	\begin{equation}
		\label{eqn:L2:tilde:sigma}
		\| \widetilde{\sigma}(t) \|_{L^{2},\lambda(F^{*}L)}^{2} 
		=
		\| \widetilde{\sigma}(t) \|_{L^{2},\lambda(F^{*}L)}^{\prime 2} 
		= 
		\det\left( \frac{i}{2} \int_{Y_{t}} F^{*}h^{H}( \omega_{l}(t) \wedge \overline{\omega}_{m}(t) ) \right)_{1\leq l,m \leq g+N-1}
	\end{equation}
	since the $L^{2}$-metric on $H^{0}(Y_{t}, K_{Y_{t}}(H))$ is independent of the choice of a K\"ahler metric on $Y_{t}$.

	\subsection{The ratio of analytic torsions via Quillen metrics}

	By Bismut-Gillet-Soul\'e \cite{BGS88},
	the line bundle $\lambda({\mathcal O}_{Y})|_{T^{o}}$ is endowed with the Quillen metric $\| \cdot \|_{Q,\lambda({\mathcal O}_{Y})}$ 
	(resp. $\| \cdot \|'_{Q,\lambda({\mathcal O}_{Y})}$) with respect to $g^{Y}$ (resp. $F^{*}g^{X}$). 
	Similarly, $\lambda(F^{*}L)|_{T^{o}}$ is endowed with the Quillen metric $\| \cdot \|_{Q,\lambda(F^{*}L)}$ 
	(resp. $\| \cdot \|'_{Q,\lambda(F^{*}L)}$ ) with respect to $g^{Y}$, $F^{*}h^{L}$ (resp. $F^{*}g^{X}$, $F^{*}h^{L}$). 
	Recall that $(L,h^{L})$ is isomorphic to a trivial Hermitian line bundle on a neighborhood $U$ of ${\rm Sing}\,X_{0}$. 
	Hence $(F^{*}L, F^{*}h^{L})$ is a trivial Hermitian line bundle on the neighborhood $F^{-1}(U)$ of $F^{-1}({\rm Sing}\,X_{0})$.

	\begin{proposition}
		\label{prop:diff:Quillen}
		There exists a constant $\gamma_{1} \in {\mathbf R}$ such that as $t \to 0$,
		$$
		\log\left( \frac{\| \cdot \|_{Q,\lambda({\mathcal O}_{Y})}}{\| \cdot \|'_{Q,\lambda({\mathcal O}_{Y})}} \right)^{2}(t)
		-
		\log\left( \frac{\| \cdot \|_{Q,\lambda(F^{*}L)}}{\| \cdot \|'_{Q,\lambda(F^{*}L)}} \right)^{2}(t)
		=
		\gamma_{1} + o(1).
		$$
	\end{proposition}
	
	\begin{pf}
		Let $\widetilde{\rm Td}(TY/T; g^{Y}, F^{*}g^{X})$ be the Bott-Chern secondary class such that
		$$
		-dd^{c} \widetilde{\rm Td}(TY/T; g^{Y}, F^{*}g^{X}) = {\rm Td}(TY/T, g^{Y}) - {\rm Td}(TY/T, F^{*}g^{X}).
		$$
		By the anomaly formula for Quillen metrics \cite[Th.\,0.2]{BGS88}, we have
		\begin{equation}
			\label{eqn:anomaly:1}
			\log\left( \frac{\| \cdot \|_{Q,\lambda({\mathcal O}_{Y})}}{\| \cdot \|'_{Q,\lambda({\mathcal O}_{Y})}} \right)^{2}(t) 
			= 
			\int_{Y_{t}} \widetilde{\rm Td}(TY_{t}; g^{Y}|_{Y_{t}}, F^{*}g^{X}|_{Y_{t}}),
		\end{equation}
		\begin{equation}
			\label{eqn:anomaly:2}
			\log\left( \frac{\| \cdot \|_{Q,\lambda(F^{*}L)}}{\| \cdot \|'_{Q,\lambda(F^{*}L)}} \right)^{2}(t) 
			= 
			\int_{Y_{t}} \widetilde{\rm Td}(TY_{t}; g^{Y}|_{Y_{t}}, F^{*}g^{X}|_{Y_{t}}) {\rm ch}(F^{*}L, F^{*}h^{L})|_{Y_{t}}.
		\end{equation}
		By the triviality of $(F^{*}L, F^{*}h^{L})$ on $F^{-1}(U)$, we get ${\rm ch}(F^{*}L, F^{*}h^{L}) = 1$ on $F^{-1}(U)$. 
		By \eqref{eqn:anomaly:1}, \eqref{eqn:anomaly:2}, we get for $t \not= 0$,
		\begin{equation}
			\label{eqn:diff:Quillen}
			\begin{aligned}
				\,&
				\log\left( \frac{\| \cdot \|_{Q,\lambda({\mathcal O}_{Y})}}{\| \cdot \|'_{Q,\lambda({\mathcal O}_{Y})}} \right)^{2}(t)
				-
				\log\left( \frac{\| \cdot \|_{Q,\lambda(F^{*}L)}}{\| \cdot \|'_{Q,\lambda(F^{*}L)}} \right)^{2}(t)
				\\
				&=
				\int_{Y_{t}\setminus F^{-1}(U)} \widetilde{\rm Td}(TY_{t}; g^{Y}|_{Y_{t}}, F^{*}g^{X}|_{Y_{t}}) \left(1 - {\rm ch}(F^{*}L, F^{*}h^{L})|_{Y_{t}} \right).
			\end{aligned}
		\end{equation}
		After shrinking $T$ if necessary,  $\widetilde{f} \colon Y\setminus F^{-1}(U) \to T$ is a trivial family of compact smooth manifolds with boundary.
		Hence the right hand side of \eqref{eqn:diff:Quillen} extends to a smooth function on $T$. This completes the proof.
	\end{pf}

	Let $\widetilde{g}^{Y}$ be a Hermitian metric on $TY/T|_{Y\setminus{\rm Sing}\,Y_{0}}$ such that for every ${\mathfrak p} \in {\rm Sing}\,Y_{0}$, one has
	$$
	\widetilde{g}^{Y}|_{U_{{\mathfrak p}\cap Y_{t}}} = \left.\frac{dz\cdot d\bar{z}}{|z|^{2}} \right|_{Y_{t}} =  \left.\frac{dw\cdot d\bar{w}}{|w|^{2}} \right|_{Y_{t}} 
	$$
	on a coordinate neighborhood $(U_{\mathfrak p}, (z,w))$ centered at ${\mathfrak p}$, where $\widetilde{f}(z,w) = zw$ on $U_{\mathfrak p}$. 
	Let $\| \cdot \|''_{Q,\lambda({\mathcal O}_{Y})}$ be the Quillen metric on $\lambda({\mathcal O}_{Y})|_{T^{o}}$  
	with respect to $\widetilde{g}^{Y}$. 
	Similarly, let $\| \cdot \|''_{Q,\lambda(F^{*}L)}$ be the Quillen metric on $\lambda(F^{*}L)|_{T^{o}}$ with respect to $\widetilde{g}^{Y}$, $F^{*}h^{L}$.

	\begin{proposition}
		\label{prop:diff:Quillen:BB}
		There exists a constant $\gamma_{2} \in {\mathbf R}$ such that as $t \to 0$,
		$$
		\log\left( \frac{\| \cdot \|_{Q,\lambda({\mathcal O}_{Y})}}{\| \cdot \|''_{Q,\lambda({\mathcal O}_{Y})}} \right)^{2}(t)
		-
		\log\left( \frac{\| \cdot \|_{Q,\lambda(F^{*}L)}}{\| \cdot \|''_{Q,\lambda(F^{*}L)}} \right)^{2}(t)
		=
		\gamma_{2} + o(1).
		$$
	\end{proposition}
	
	\begin{pf}
		Replacing $F^{*}g^{X}$ with $\widetilde{g}^{Y}$, we can prove the assertion in the same way as in Proposition~\ref{prop:diff:Quillen}.
	\end{pf}

	\begin{theorem}
		\label{thm:ratio:torsion:1}
		There exists a constant $\gamma \in {\mathbf R}$ such that as $t\to 0$,
		$$
		\begin{aligned}
			\,&
			\log \frac{ \tau(X_{\mu(t)},{\mathcal O}_{X_{\mu(t)}}) }{ \tau(X_{\mu(t)}, L_{\mu(t)}) } =
			\log \frac{\tau(Y_{t}, {\mathcal O}_{Y_{t}};F^{*}g^{X}) }{ \tau(Y_{t}, F^{*}L_{t};F^{*}g^{X},F^{*}h^{L}) }
			\\
			&=
			\log 
			\left[
			\frac{ \det\left(  \int_{Y_{t}} h_{F^{*}H}(\omega_{i}(t) \wedge \overline{\omega_{j}(t)}) \right)_{1\leq i,j \leq g+N-1} }
			{ \det\left(  \int_{Y_{t}} \omega_{i}(t) \wedge \overline{\omega_{j}(t)} \right)_{1\leq i,j \leq g} }
			\right]
			+ \gamma + o(1).
		\end{aligned}
		$$
	\end{theorem}
	
	\begin{pf}
		
		Since $\tau(X_{\mu(t)}, {\mathcal O}_{X_{\mu(t)}}) = \tau(Y_{t}, {\mathcal O}_{Y_{t}};F^{*}g^{X})$
		and 
		$$
		\tau(X_{\mu(t)}, L_{\mu(t)}) = \tau(Y_{t}, F^{*}L_{t};F^{*}g^{X},F^{*}h^{L}), 
		$$
		it suffices to prove the second equality. 
		By the definition of Quillen metrics and \eqref{eqn:L2:met:2}, \eqref{eqn:L2:tilde:sigma}, we have
		\begin{equation}
			\label{eqn:ratio:Quillen:0}
			\begin{aligned}
				\log\left( \frac{\| \sigma \|'_{Q,\lambda({\mathcal O}_{Y})}}{\| \widetilde{\sigma} \|'_{Q,\lambda(F^{*}L)}} \right)^{2}(t)
				&=
				\log \frac{\tau(Y_{t}, {\mathcal O}_{Y_{t}};F^{*}g^{X})}{\tau(Y_{t}, F^{*}L_{t};F^{*}g^{X},F^{*}h^{L})}
				+
				\log\left( \frac{\| \sigma \|'_{L^{2},\lambda({\mathcal O}_{Y})}}{\| \widetilde{\sigma} \|'_{L^{2},\lambda(F^{*}L)}} \right)^{2}(t)
				\\
				&=
				\log \frac{\tau(Y_{t}, {\mathcal O}_{Y_{t}};F^{*}g^{X})}{\tau(Y_{t}, F^{*}L_{t};F^{*}g^{X},F^{*}h^{L})}
				+
				\log \det\left( \displaystyle \int_{Y_{t}} \omega_{i}(t) \wedge \overline{\omega_{j}(t)} \right)_{1\leq i,j \leq g}
				\\
				&\quad - \det\left( \displaystyle \int_{Y_{t}} h_{F^{*}H}(\omega_{i}(t) \wedge \overline{\omega_{j}(t)}) \right)_{1\leq i,j \leq g+N-1} - A_{1}.
			\end{aligned}
		\end{equation}
		By Bismut-Bost \cite[Th.\,2.2]{BismutBost90}, there exist $\gamma_{3}, \gamma_{4} \in {\mathbf R}$ such that as $t \to 0$,
		\begin{equation}
			\label{eqn:sing:Quillen:1}
			\log (\| \sigma \|''_{Q,\lambda({\mathcal O}_{Y})})^{2}(t) = \frac{ \# {\rm Sing}\,Y_{0}}{12} \log |t|^{2} + \gamma_{3} + o(1),
		\end{equation}
		\begin{equation}
			\label{eqn:sing:Quillen:2}
			\log (\| \widetilde{\sigma} \|''_{Q,\lambda(F^{*}L)})^{2}(t) = \frac{ \# {\rm Sing}\,Y_{0}}{12} \log |t|^{2} + \gamma_{4} + o(1).
		\end{equation}
		By Propositions~\ref{prop:diff:Quillen} and \ref{prop:diff:Quillen:BB} and \eqref{eqn:sing:Quillen:1}, \eqref{eqn:sing:Quillen:2}, as $t \to 0$, we get
		\begin{equation}
			\label{eqn:ratio:Quillen:2}
			\begin{aligned}
				\log\left( \frac{\| \sigma \|'_{Q,\lambda({\mathcal O}_{Y})}}{\| \widetilde{\sigma} \|'_{Q,\lambda(F^{*}L)}} \right)^{2}(t)
				&=
				\log\left( \frac{\| \sigma \|''_{Q,\lambda({\mathcal O}_{Y})}}{\| \widetilde{\sigma} \|''_{Q,\lambda(F^{*}L)}} \right)^{2}(t) - \gamma_{1} + \gamma_{2} + o(1)
				\\
				&=
				\gamma_{3} - \gamma_{4} - \gamma_{1} + \gamma_{2} + o(1).
			\end{aligned}
		\end{equation}
		Comparing \eqref{eqn:ratio:Quillen:0} and \eqref{eqn:ratio:Quillen:2}, we get the result.
		This completes the proof.
	\end{pf}

	Let $C^{*}(T)$ and $C^{*}(T^{o})$ be the abelian groups of nowhere vanishing real valued continuous functions on $T$ and $T^{o} = T \setminus \{0\}$, respectively,
	where the group structure is given by the point wise multiplication of functions. The equality in $C^{*}(T^{o})/C^{*}(T)$ is denoted by $\equiv$. 
	By Theorems~\ref{thm:at:bd1} and \ref{thm:ratio:torsion:1}, we have the following:

	\begin{corollary}
		\label{cor:McKean:Singer:mult}
		The following identity holds in $C^{*}(T^{o})/C^{*}(T)$:
		$$
		\prod_{i=1}^{N-1} \lambda_{i}(\mu(t))^{-1} 
		\equiv
		\frac{ \tau(X_{\mu(t)},{\mathcal O}_{X_{\mu(t)}}) }{ \tau(X_{\mu(t)}, L_{\mu(t)}) }
		\equiv
		\frac{ \det\left( \int_{Y_{t}} h_{F^{*}H}(\omega_{i}(t) \wedge \overline{\omega_{j}(t)}) \right) }{ \det\left( \int_{Y_{t}} \omega_{i}(t) \wedge \overline{\omega_{j}(t)} \right) }.
		$$
	\end{corollary}

	\section{Asymptotic behavior of the determinants of the period integrals}
	\label{sect:5}
	\par
	In this section, we determine the asymptotic behavior of $\tau(X_{s},{\mathcal O}_{X_{s}}) / \tau(X_{s}, L)$ as $s \to 0$.
	To do this, in view of Theorem~\ref{thm:ratio:torsion:1}, we determine the singularity of the $L^{2}$-metrics $\| \cdot \|_{L^{2},\lambda(F^{*}L)}$, 
	$\| \cdot \|'_{L^{2},\lambda(F^{*}L)}$, $\| \cdot \|_{L^{2},\lambda({\mathcal O}_{Y})}$, $\| \cdot \|'_{L^{2},\lambda({\mathcal O}_{Y})}$.
	Throughout this section, we keep the notation of Section~\ref{sect:4}.

	\subsection{Determinants of the period integrals}
	\par
	Let $\nu \colon \widetilde{Y}_{0} \to Y_{0}$ be the normalization. Let $k\in {\mathbf N}$ be such that
	$$
	\nu^{*}\omega_{1}(0),\ldots,\nu^{*}\omega_{k}(0) \in H^{0}(\widetilde{Y}_{0}, K_{\widetilde{Y}_{0}}),
	\quad
	\nu^{*}\omega_{k+1}(0),\ldots,\nu^{*}\omega_{g}(0) \not\in H^{0}(\widetilde{Y}_{0}, K_{\widetilde{Y}_{0}}).
	$$

	\begin{proposition}
		\label{prop:Barlet}
		The following hold.
		\begin{itemize}
			\item[(1)]
			There exist constants $a_{ij}$, $b_{ij}$ $(1 \leq i,j \leq g+N-1)$ such that
			$$
			\sqrt{-1} \int_{Y_{t}} F^{*}h^{H}( \omega_{i}(t) \wedge \overline{\omega}_{j}(t) ) 
			=
			a_{ij} \log |t|^{-2} + b_{ij} + o(1)
			\qquad (t\to 0).
			$$
			\item[(2)]
			$a_{ij}=0$ if $1 \leq i \leq k$ or $1 \leq j \leq k$.
			\item[(3)]
			The Hermitian matrices $(b_{ij})_{1\leq i,j \leq k}$ and $(a_{ij})_{k+1 \leq i,j \leq g+N-1}$ are positive-definite. 
		\end{itemize}
	\end{proposition}
	
	\begin{pf}
		Let ${\mathfrak p} \in {\rm Sing}\,Y_{0}$. Let $(x,y)$ be a system of coordinates centered at ${\mathfrak p}$ defined on $U \subset Y$ such that $g(x,y) =xy$. 
		Near ${\mathfrak p}$, we can express
		$$
		\omega_{i}(t)(x,y) = \alpha_{i}(x,y) \left.\frac{dx}{x} \otimes {\mathbf e}\right|_{Y_{t}\cap U},
		$$
		where $\alpha_{i}(x,y) \in {\mathcal O}(U)$ and ${\mathbf e} \in \Gamma(U, F^{*}H)$ is a holomorphic frame of $F^{*}H$ on $U$.
		Since $h^{L}$ is flat near $F({\mathfrak p})$, we can assume that $F^{*}h^{H}({\mathbf e}, {\mathbf e})=1$ on $U$.
		Rescaling the coordinates if necessary, we may assume that $\overline{\varDelta}_{\mathfrak p}^{2} \subset U$,
		where $\overline{\varDelta}_{\mathfrak p}^{2}$ is the closed unit polydisc centered at ${\mathfrak p}$.
		Then, as $t \to 0$, 
		\begin{equation}
			\label{eqn:fib:integral:1}
			\int_{Y_{t}} \rho F^{*}h^{H}( \omega_{i}(t) \wedge \overline{\omega}_{j}(t) ) 
			= 
			\sum_{{\mathfrak p} \in {\rm Sing}\,Y_{0}} \int_{Y_{t}\cap \overline{\varDelta}_{\mathfrak p}^{2}} 
			\alpha_{i}(x,y)\overline{\alpha_{j}(x,y)} \left.\frac{dx\wedge d\bar{x}}{|x|^{2}} \right|_{Y_{t}\cap \overline{\varDelta}_{\mathfrak p}^{2}} + c + o(1),
		\end{equation}
		where 
		$c = \int_{Y_{0}\setminus\bigcup_{{\mathfrak p}\in{\rm Sing}\,Y_{0}}{\varDelta}_{\mathfrak p}^{2}} F^{*}h^{H}( \omega_{i}(0) \wedge \overline{\omega}_{j}(0))$.
		Since $Y_{t}\cap\overline{\varDelta}^{2}_{\mathfrak p} \cong \{ |t| \leq |x| \leq 1 \}$ is an annulus,
		making use of the Taylor series expansion of the holomorphic functions $\alpha_{i}, \alpha_{j} \in {\mathcal O}(\overline{\varDelta}^{2})$, we infer 
		in the same way as in \cite[Prop.\,13.5]{BismutBost90} (see also \cite[p.140, proof of Lemma 2, cas 1 et 2]{Barlet82}, \cite[Lemma 3.4]{Takayama21} with
		$d=d'=1$, $q=0$, $w=dx\wedge d\bar{x}$) that as $t \to 0$,
		\begin{equation}
			\label{eqn:fib:integral:2}
			\begin{aligned}
				\,&
				\sqrt{-1} \int_{|t| \leq |x| \leq 1} \alpha_{i}(x,t/x)\overline{\alpha_{j}(x, t/x)} \frac{dx\wedge d\bar{x}}{|x|^{2}}
				\\
				&=
				4\pi \alpha_{i}(0,0)\overline{\alpha_{j}(0,0)} \log |t|^{-2} 
				+
				\sqrt{-1} \int_{\varDelta}\{ \alpha_{i}(x,0)\overline{\alpha_{j}(x,0)} - \alpha_{i}(0,0)\overline{\alpha_{j}(0,0)} \} \frac{dx\wedge d\bar{x}}{|x|^{2}}
				\\
				&\quad+
				\sqrt{-1} \int_{\varDelta}\{ \alpha_{i}(0,y)\overline{\alpha_{j}(0,y)} - \alpha_{i}(0,0)\overline{\alpha_{j}(0,0)} \} \frac{dy\wedge d\bar{y}}{|y|^{2}} + O\left( |t|\log|t| \right).
			\end{aligned}
		\end{equation}
		Since $\alpha_{i}(0,0)\overline{\alpha_{j}(0,0)} = {\rm Res}_{\mathfrak p}\omega_{i}(0) \overline{{\rm Res}_{\mathfrak p}\omega_{j}(0)}$, 
		we deduce from \eqref{eqn:fib:integral:1}, \eqref{eqn:fib:integral:2} that
		\begin{equation}
			\label{eqn:fib:integral:3}
			\sqrt{-1}\int_{Y_{t}} F^{*}h^{H}( \omega_{i}(t) \wedge \overline{\omega}_{j}(t) ) 
			= 
			4\pi ({\rm Res}_{\mathfrak p}\omega_{i}(0) \overline{{\rm Res}_{\mathfrak p}\omega_{j}(0)}) \log |t|^{-2} + b_{ij} + o(1),
		\end{equation}
		where $b_{ij}$ is a constant.
		This proves (1).
		Since $\nu^{*}\omega_{i}(0)$ is a regular $1$-form on $\widetilde{Y}_{0}$, we have ${\rm Res}_{\mathfrak p}\omega_{i}(0) =0$ for all
		${\mathfrak p} \in {\rm Sing}\,Y_{0}$ and $1 \leq i \leq k$.
		By \eqref{eqn:fib:integral:3}, we get $a_{ij}=0$ if $1 \leq i \leq k$ or $1 \leq j \leq k$. This proves (2).
		\par
		Let ${\mathbf c} = (c_{k+1}, \ldots, c_{g+N-1})\in {\mathbf C}^{N+g-1}$ be such that $\| {\mathbf c} \|^{2} = \sum_{i} |c_{i}|^{2} = 1$.
		Set $\varphi(t) := \sum_{i=k+1}^{g+N-1} c_{i} \omega_{i}(t)$. We have $H_{0} = {\mathcal O}_{X_{0}}(\sum_{i} p_{i} )$ with $p_{i} \in C_{i}\setminus{\rm Sing}\,X_{0}$. 
		Since $F \colon Y_{0} \setminus F^{-1}({\rm Sing}\,X_{0}) \to X_{0} \setminus {\rm Sing}\,X_{0}$ is an isomorphism, there exists a unique 
		$\widetilde{p}_{i} \in \widetilde{C}_{i}\setminus F^{-1}({\rm Sing}\,X_{0}) \subset Y_{0}\setminus{\rm Sing}\,Y_{0}$ with $F(\widetilde{p}_{i}) = p_{i}$.
		By \eqref{eqn:fib:integral:3} and (1), we have
		\begin{equation}
			\label{eqn:fib:integral:4}
			\| \varphi(t) \|_{L^{2}}^{2} 
			= 
			4\pi (\sum_{{\mathfrak q}\in {\rm Sing}\,Y_{0}} \sum_{j} |{\rm Res}_{\mathfrak q} \varphi(0)|_{\widetilde{C}_{j}}|^{2}) \log |t|^{-2} + \gamma + o(1)
			\quad
			(t\to 0).
		\end{equation}
		\par{\em (Case 1) }
		Suppose $(c_{g+1},\ldots,c_{g+N-1}) \not= (0,\ldots,0)$. Then there exist $i_{0} \in \{ 1, \ldots, N-1 \}$ with ${\rm Res}_{\widetilde{p}_{i_{0}}}(\varphi(0)) \not=0$.
		Indeed, if ${\rm Res}_{\widetilde{p}_{i}}(\varphi(0)) =0$ for $1\leq i \leq N-1$, then $\varphi(0) \in H^{0}(Y_{0}, K_{Y_{0}})$. 
		Hence we can express $\varphi(0) = \sum_{j=1}^{g} d_{j} \omega_{j}(0)$. Namely, 
		$\sum_{j=1}^{g} d_{j} \omega_{j}(0) - \sum_{i=k+1}^{g+N-1} c_{i} \omega_{i}(0) =0$.
		Since $\omega_{1}(0), \ldots, \omega_{g+N-1}(0)$ is a basis of $H^{0}(Y_{t}, K_{Y_{t}}(F^{*}H))$, we get $c_{g+1}=\cdots=c_{g+N-1}=0$. 
		This contradicts $(c_{g+1},\ldots,c_{g+N-1}) \not= (0,\ldots,0)$. Let ${\rm Res}_{\widetilde{p}_{i_{0}}}(\varphi(0)) \not=0$ in what follows. 
		\par
		If ${\rm Res}_{\mathfrak q}\varphi(0)|_{\widetilde{C}_{i_{0}}}= 0$ for all ${\mathfrak q} \in{\rm Sing}\,Y_{0}\cap\widetilde{C}_{i_{0}}$,
		since $\varphi(0)|_{\widetilde{C}_{i_{0}}}$ is a logarithmic $1$-form on $\widetilde{C}_{i_{0}}$, which is holomorphic on 
		$\widetilde{C}_{i_{0}} \setminus ({\rm Sing}\,Y_{0} \cup \{ \widetilde{p}_{i_{0}} \})$, the residue theorem 
		implies ${\rm Res}_{\widetilde{p}_{i_{0}}} \varphi(0)|_{\widetilde{C}_{i_{0}}} =0$. This contradicts ${\rm Res}_{\widetilde{p}_{i_{0}}}(\varphi(0)) \not=0$. 
		Hence ${\rm Res}_{\mathfrak q}\varphi(0)|_{\widetilde{C}_{i_{0}}} \not= 0$ for some ${\mathfrak q} \in {\rm Sing}\,Y_{0} \cap \widetilde{C}_{i_{0}}$. 
		By \eqref{eqn:fib:integral:4}, there exists $\alpha>0$ with
		$$
		\| \varphi(t) \|_{L^{2}}^{2} = \alpha \log |t|^{-2} + \gamma + o(1)
		\qquad
		(t \to 0).
		$$
		\par{\em (Case 2) }
		Suppose $(c_{g+1},\ldots,c_{g+N-1}) = (0,\ldots,0)$ and $(c_{k+1},\ldots,c_{g}) \not= (0,\ldots,0)$. Then $\varphi(t) \in H^{0}(Y_{t}, K_{Y_{t}})$. 
		There exist $j_{0} \in \{1, \ldots, N-1\}$ and ${\mathfrak q} \in {\rm Sing}\, Y_{0} \cap \widetilde{C}_{j_{0}}$ such that 
		${\rm Res}_{\mathfrak q} \varphi(0)|_{\widetilde{C}_{j_{0}}} \not= 0$.
		Indeed, if ${\rm Res}_{\mathfrak q} \varphi(0)|_{\widetilde{C}_{j}} = 0$ for all $1\leq j \leq N-1$ and ${\mathfrak q} \in {\rm Sing}\, Y_{0}\cap \widetilde{C}_{j}$,
		then $\nu^{*}\varphi(0) \in H^{0}(\widetilde{Y}_{0}, K_{\widetilde{Y}_{0}})$. Hence, we can express $\varphi(0) = \sum_{i=1}^{k} d_{i} \omega_{i}(0)$.
		Since $\{ \omega_{1}(0), \ldots, \omega_{g}(0) \}$ is a basis of $H^{0}(Y_{0}, K_{Y_{0}})$, we get a contradiction $c_{k+1}=\cdots=c_{g}=0$.
		Since ${\rm Res}_{\mathfrak q} \varphi(0)|_{\widetilde{C}_{j_{0}}} \not= 0$ for some $j$ and ${\mathfrak q} \in {\rm Sing}\,Y_{0}\cap \widetilde{C}_{j}$,
		we deduce from \eqref{eqn:fib:integral:4} that there exists $\alpha'>0$ with
		$$
		\| \varphi(t) \|_{L^{2}}^{2} = \alpha' \log |t|^{-2} + \gamma + o(1)
		\qquad
		(t \to 0).
		$$
		\par
		Set $A=(a_{ij})_{k+1\leq i,j \leq g+N-1}$. Now we prove that $A$ is positive-definite. By \eqref{eqn:fib:integral:4}, $A$ is positive-semidefinite.
		Let ${\mathbf c} = (c_{k+1},\ldots,c_{g+N-1})\not={\mathbf 0}$ be such that $\sum_{i,j=k+1}^{g+N-1} a_{ij} c_{i} \overline{c}_{j} =0$.
		Set $\varphi := \sum_{i=k+1}^{g+N-1} c_{i}\omega_{i}$. By definition of $a_{ij}$ and (1), this implies 
		$$
		\| \varphi(t) \|_{L^{2}}^{2} = (\sum_{i,j=k+1}^{g+N-1} a_{ij} c_{i} \overline{c}_{j} ) \log |t|^{-2} + \gamma + o(1) = \gamma + o(1)
		\qquad
		(t \to 0).
		$$
		By Cases 1 and 2, we obtain ${\mathbf c} = {\mathbf 0}$. This proves that $A$ is positive-definite.
		\par
		Next, we prove that $B:= (b_{ij})_{1\leq i,j \leq k}$ is positive-definite. Since
		$$
		b_{ij} = \sqrt{-1} \int_{\widetilde{Y}_{0}} \nu^{*}\omega_{i}(0) \wedge \nu^{*}\overline{\omega}_{j}(0)
		\qquad
		(1\leq i, j \leq k)
		$$
		and $\nu^{*}\omega_{1}(0), \ldots, \nu^{*}\omega_{k}(0)$ are linearly independent vectors of $H^{0}(\widetilde{Y}_{0}, K_{\widetilde{Y}_{0}})$,
		$B$ is positive definite. This completes the proof of (3).
	\end{pf}

	Recall that $\widetilde{\sigma} = \omega_{1}\wedge\cdots\wedge\omega_{g+N-1}$ is a nowhere vanishing holomorphic section of 
	$\lambda(F^{*}L) \cong \widetilde{f}_{*}(K_{Y/T}(F^{*}H))$ near $0 \in T$.

	\begin{proposition}
		\label{prop:period:1}
		There exists $\gamma_{0} \in {\mathbf R}$ such that as $t \to 0$,
		$$
		\begin{aligned}
			\log \| \widetilde{\sigma} \|_{L^{2},\lambda(F^{*}L)}^{2} (t) 
			&= 
			\log \| \widetilde{\sigma} \|_{L^{2},\lambda(F^{*}L)}^{\prime\, 2} (t) 
			\\
			&= 
			(g+N-1-k) \log\log( |t|^{-1}) + \gamma_{0} + o(1).
		\end{aligned}
		$$
	\end{proposition}
	
	\begin{pf}
		The first equality follows from the fact that the $L^{2}$-metric on $\widetilde{f}_{*}K_{Y/T}(F^{*}H)$ is independent of the choice of a Hermitian metric on
		the relative tangent bundle $TY/T$ which is fiberwise K\"ahler. Since
		$$
		\begin{aligned}
			\| \widetilde{\sigma} \|_{L^{2},\lambda(F^{*}L)}^{2} (t) 
			&= 
			\det\left( \frac{i}{2} \int_{Y_{t}} F^{*}h^{H}( \omega_{l}(t) \wedge \overline{\omega}_{m}(t) ) \right)_{1\leq l,m \leq g+N-1}
			\\
			&=
			\det B\cdot \det A \cdot ( \log |t|^{-2} )^{g+N-1-k} + O(( \log |t|^{-2} )^{g+N-2-k})
		\end{aligned}
		$$
		by Proposition~\ref{prop:Barlet} (1) and since $A$ and $B$ are positive-definite by Proposition~\ref{prop:Barlet} (3), we get the second equality.  
		Notice that $\log\log |t|^{-2} = \log\log |t|^{-1} + \log 2$.
	\end{pf}

	Recall that $\sigma = 1 \otimes(\omega_{1}\wedge\cdots\wedge\omega_{g})$ is a nowhere vanishing holomorphic section of $\lambda({\mathcal O}_{Y})$
	near $0 \in T$.

	\begin{proposition}
		\label{prop:period:2}
		There exist $\gamma_{1}, \gamma_{2} \in {\mathbf R}$ such that as $t \to 0$,
		$$
		\log \| \sigma \|_{L^{2},\lambda({\mathcal O}_{Y})}^{2} (t) = (g-k) \log\log (|t|^{-1}) + \gamma_{1} + o(1),
		$$
		$$
		\log \| \sigma \|_{L^{2},\lambda({\mathcal O}_{Y})}^{\prime\, 2} (t) = (g-k) \log\log (|t|^{-1}) + \gamma_{2} + o(1).
		$$
	\end{proposition}
	
	\begin{pf}
		Since $Y_{0}$ has at most ordinary double points, the monodromy of $\widetilde{f} \colon Y \to T$ around $t=0$ is unipotent. 
		By \cite[Th.\,C]{EFM21}, there exists a constant $c$ such that as $t \to 0$,
		\begin{equation}
			\label{eqn:fib:integral:5}
			\log \det\left( \frac{i}{2} \int_{Y_{t}} \omega_{l} \wedge \overline{\omega}_{m} \right)_{1\leq l,m \leq g} = (g-k) \log\log (|t|^{-1}) + c + o(1).
		\end{equation}
		By \eqref{eqn:L2:met:1}, \eqref{eqn:L2:met:2}, \eqref{eqn:fib:integral:5}, we get the result.
	\end{pf}

	\begin{remark}
		It is possible to prove \eqref{eqn:fib:integral:5} in the same way as the proof of Proposition~\ref{prop:Barlet}. Since the proof is parallel,
		we leave the details to the reader.
	\end{remark}

	\subsection{Asymptotic behavior of the ratio of analytic torsions}

	\begin{theorem}
		\label{thm:ratio:torsion:2}
		There exists a constant $\gamma \in {\mathbf R}$ such that as $s\to 0$,
		$$
		\log \frac{ \tau(X_{s}, {\mathcal O}_{X_{s}}) }{ \tau(X_{s}, L_{s}) } = -(N-1) \log\log (|s|^{-1}) + \gamma + o(1).
		$$
	\end{theorem}
	
	\begin{pf}
		Since $\mu(t) = t^{\deg\mu}$, the result follows from Theorem~\ref{thm:ratio:torsion:1} and Propositions~\ref{prop:period:1} and \ref{prop:period:2}.
	\end{pf}

	\subsection{Proof of Theorem~\ref{Main:Thm:2}}
	The result follows from Corollary~\ref{cor:McKean:Singer:mult} and Theorem~\ref{thm:ratio:torsion:2}.
	\qed

	\section{An upper bound of the small eigenvalues}
	\label{sect:An upper bound of the small eigenvalues}

	In this section, we give an upper bound of the small eigenvalues.

	\begin{proposition}
		\label{prop:lb:small:ev}
		There exist constants $K(i)>0$ $(1\leq i \leq N-1)$ such that
		$$
		\lambda_{i}(s) \leq  \frac{K(i)}{ \log (|s|^{-1}) }
		\qquad
		(s\in S^{o}).
		$$
	\end{proposition}

	\subsection{Some intermediary results}
	\par
	For every $p \in {\rm Sing}\,X_{0}$, we fix a system of coordinates $\zeta=(\zeta_{1}, \zeta_{2})$ centered at $p$. We denote by $\| \cdot \|$ the norm
	with respect to the Euclidean metric $\sum_{i}d\zeta_{i}d\bar{\zeta}_{i}$.

	\begin{lemma}
		\label{lemma:grad:f}
		There exists an integer $\nu \in {\mathbf N}$ and a constant $K_{0}>0$ such that the following inequality holds on a neighborhood of each  $p\in{\rm Sing}\,X_{0}$
		$$
		\| df(\zeta) \|^{2} \geq K_{0} \| \zeta \|^{2\nu},
		$$
		where $\|\zeta\|^{2} = |\zeta_{1}|^{2} + |\zeta_{2}|^{2}$.
	\end{lemma}
	
	\begin{pf}
		Since $f(z)$ has an isolated critical point at $z=0$, there exists $\nu \in {\mathbf N}$
		such that the Jacobi ideal $(\frac{\partial f}{\partial \zeta_{1}}, \frac{\partial f}{\partial \zeta_{2}})$ generates ${\mathfrak m}_{0}^{\nu}$, where ${\mathfrak m}_{0}$
		is the maximal ideal of ${\mathcal O}_{X,p}$. Hence there exist $g_{ij} \in {\mathbf C}\{ \zeta_{1}, \zeta_{2} \}$ such that
		$\zeta_{i}^{\nu} = \sum_{j=1}^{2} g_{ij} \frac{\partial f}{\partial \zeta_{j}}$ $(i=1,2)$ on a small neighborhood $U \subset X$ of $p$. 
		Then $\sum_{i} |\zeta_{i}|^{2\nu} \leq (\sum_{i,j} |g_{ij}|^{2}) (\sum_{k} |\frac{\partial f}{\partial \zeta_{k}}|^{2})$. The result follows easily from this inequality.
	\end{pf}

	Define a smooth vector field $\Theta$ of type $(1,0)$ on $X \setminus {\rm Sing}\,X_{0}$ by
	$$
	\Theta := \frac{g^{T^{*}X}( \cdot, df)}{\| df \|^{2}}.
	$$
	Then $\langle f_{*}\Theta, dt \rangle = \langle \Theta, df \rangle =1$. Since $\Theta$ is of type $(1,0)$, we get $f_{*}\Theta = \partial/\partial t$. 
	We define real vector fields $U, V$ on $X \setminus {\rm Sing}\,X_{0}$ by 
	$$
	U - iV := 2\Theta.
	$$
	Set $u:= {\rm Re}\,t$, $v:= {\rm Im}\,t$. Then we have
	\begin{equation}
		\label{eqn:v:f}
		f_{*}U = \frac{\partial}{\partial u},
		\qquad
		f_{*}V = \frac{\partial}{\partial v}.
	\end{equation}
	Let $p \in {\rm Sing}\,X_{0}$. Let $B(p,1)=\{\zeta\in{\mathbf C}^{2};\, \|\zeta\|<1\} \subset X$ be the unit ball centered at $p$.
	By Lemma~\ref{lemma:grad:f}, there exists a constant $C>0$ such that for all $\zeta \in B(p,1)\setminus\{0\}$, 
	$$
	\| U(\zeta) \| + \| V(\zeta) \| \leq C \| \zeta \|^{-2\nu},
	\qquad
	\| \nabla U(\zeta) \| + \| \nabla V(\zeta) \| \leq C \| \zeta \|^{-4\nu}.
	$$
	For $0<r\ll1$, we set $M_{r} := Cr^{-2\nu}$, $N_{r} := Cr^{-4\nu}$. On $B(p,1) \setminus B(p,r)$, we have
	$$
	\|U(\zeta)\| + \|V(\zeta)\| \leq M_{r},
	\qquad
	\| \nabla U(\zeta) \| + \| \nabla V(\zeta) \| \leq N_{r}.
	$$  
	\par
	Let $0<\delta <\min\{ r/M_{r}, 1/(2N_{r}) \} = r^{4\nu}/(2C)$. 
	For $z \in X_{0} \setminus \bigcup_{p\in {\rm Sing}\,X_{0}} B(p, 2r)$ and $\theta \in [0, 2\pi]$, 
	let $\Phi^{\theta}(\eta, z) \in C^{\infty}( [-\delta,\delta], X)$ be the unique solution of the ordinary differential equation
	\begin{equation}
		\label{eqn:ODE}
		\begin{cases}
			\displaystyle \frac{d}{d\eta}\Phi^{\theta}(\eta,z) = \cos\theta \cdot U_{ \Phi^{\theta}(\eta, z) } + \sin\theta \cdot V_{ \Phi^{\theta}(\eta,z) }
			\quad
			(-\delta \leq \eta \leq \delta),
			\\
			\Phi^{\theta}(0,z) = z \in X_{0} \setminus \bigcup_{p\in {\rm Sing}\,X_{0}} B(p, r).
		\end{cases}
	\end{equation}
	Since
	$\frac{d}{d\eta}f( \Phi^{\theta}(\eta,z) ) = f_{*}( \frac{d\Phi^{\theta}}{d\eta}(\eta,z) ) =
	\cos\theta \cdot \left(\frac{\partial}{\partial u} \right)_{f( \Phi^{\theta}(\eta, z) )} + \sin\theta \cdot \left( \frac{\partial}{\partial v} \right)_{f( \Phi^{\theta}(\eta,z) )}$
	by \eqref{eqn:v:f}, we have $f( \Phi^{\theta}(\eta, z) ) = \eta e^{i\theta}$. Hence $\Phi^{\theta}(\eta, z) \in X_{\eta e^{i\theta}}$. 
	
	Since $\| U(\zeta) \| + \| V(\zeta) \| \leq M_{r}$ on $B(p,1)\setminus B(p,r)$, we have 
	\begin{equation}
		\label{eqn:diffeo}
		\left\| \Phi^{\theta}(\eta, z) - \Phi^{\theta}(0,z) \right\| \leq M_{r} |\eta|
	\end{equation}
	for all $(\eta,\theta, z) \in [-\delta, \delta] \times [0,2\pi] \times \{X_{0}\cap(B(p,1/2)\setminus B(p, 2r))\}$.
	Then $\| \Phi^{\theta}(\eta, z) \| \geq \| \Phi^{\theta}(0,z) \| - M_{r}\delta \geq r$ and $\| \Phi^{\theta}(\eta, z) \| \leq \| \Phi^{\theta}(0,z) \| + M_{r}\delta \leq \frac{1}{2} +r <1$.
	Hence $\Phi^{\theta}(\eta, z) \in X_{\eta e^{i\theta}} \setminus \bigcup_{p\in {\rm Sing}\,X_{0}} B(p, r)$ for 
	$(\eta,\theta,z) \in [0, \delta] \times [0,2\pi] \times (X_{0}\setminus \bigcup_{p\in {\rm Sing}\,X_{0}} B(p, 2r))$. 
	Similarly, by fixing a system of local coordinates on a neighborhood of $X_{0}$, we may assume that \eqref{eqn:diffeo} holds on 
	$X_{0} \setminus \bigcup_{p \in {\rm Sing}\,X_{0}} B(p,r)$.
	\par
	Define 
	$$
	\Phi^{\theta}_{\eta} \colon X_{0}\setminus \bigcup_{p\in {\rm Sing}\,X_{0}} B(p, 2r)  \ni z \to \Phi^{\theta}_{\eta}(z) := \Phi^{\theta}(\eta, z) \in 
	X_{\eta e^{i\theta}} \setminus \bigcup_{p\in {\rm Sing}\,X_{0}} B(p, r).
	$$
	By the uniqueness of the solution of \eqref{eqn:ODE}, $\Phi_{\eta}^{\theta}$ is a diffeomorphism from $X_{0} \setminus \bigcup_{p \in {\rm Sing}\,X_{0}} B(p,2r)$ 
	to $\Phi_{\eta}^{\theta}( X_{0} \setminus \bigcup_{p \in {\rm Sing}\,X_{0}} B(p,2r) )$ for $\eta \in (-\delta, \delta)$. 
	Let $(\Phi^{\theta}_{\eta})_{*,z} \in {\rm Hom}(T^{\mathbf R}_{z}X_{0}, T^{\mathbf R}_{\Phi(\eta, z)}X)$ be the differential of the map $\Phi^{\theta}_{\eta}$ at $z$. 
	Identifying $T^{\mathbf R}_{\Phi^{\theta}(\eta,z)}X$ with ${\mathbf R}^{4}$, we get $(\Phi^{\theta}_{\eta})_{*,z} \in {\rm Hom}(T^{\mathbf R}_{z}X_{0}, {\mathbf R}^{4})$.
	Hence $(\Phi^{\theta}_{\eta})_{*,z} - (\Phi^{\theta}_{0})_{*,z} \in {\rm Hom}(T^{\mathbf R}_{z}X_{0}, {\mathbf R}^{4})$. 
	In the next lemma, the norm $\| (\Phi^{\theta}_{\eta})_{*,z} - (\Phi^{\theta}_{0})_{*,z} \|$ is the one with respect to the Euclidean metric on ${\mathbf C}^{2}={\mathbf R}^{4}$. 
	Let us consider the case $z \in B(p, 1)$.
	Since $\Phi^{\theta}_{0}(z) = z$ is the identity map, $(\Phi^{\theta}_{0})_{*,z}$ is the inclusion map $T^{\mathbf R}_{z}X_{0} \hookrightarrow T^{\mathbf R}_{z}B(p,1)$.
	Since the metric on $X_{0}$ is induced from the metric $g^{X}$ on $X$ and $g^{X}$ is quasi-isometric to the Euclidean metric on $B(p,1)$, 
	this implies the existence of  a constant $K>0$ with
	\begin{equation}
		\label{eqn:est:diffeo}
		\| (\Phi^{\theta}_{0})_{*,z} \| \leq K
		\qquad
		(z \in X_{0}\cap(B(p,1)\setminus\{0\})).
	\end{equation}
	Similarly, replacing $K$ with another constant if necessary, we may assume that \eqref{eqn:est:diffeo} holds on $X_{0} \setminus {\rm Sing}\,X_{0}$.

	\begin{lemma}
		\label{lemma:diffeo}
		There exists a constant $K_{1}>0$ such that 
		$$
		\left\| (\Phi^{\theta}_{\eta})_{*,z} - (\Phi^{\theta}_{0})_{*,z} \right\| \leq K_{1}N_{r} \eta
		$$
		for all $(\eta, \theta,z) \in [0, \delta] \times [0,2\pi] \times ( X_{0} \setminus \bigcup_{p\in{\rm Sing}\,X_{0}} B(2r) ))$. 
		In particular, for all $(\eta, \theta, z) \in [0, r^{8\nu}] \times [0,2\pi] \times ( X_{0} \setminus \bigcup_{p\in {\rm Sing}\,X_{0}} B(p,2r) ))$, one has
		$$
		\left\| (\Phi^{\theta}_{\eta})_{*,z} - (\Phi^{\theta}_{0})_{*,z} \right\| \leq K_{1}C \eta^{\frac{1}{2}}.
		$$
	\end{lemma}
	
	\begin{pf}
		It suffices to prove the assertion when $z \in B(p,1) \setminus B(p, 2r)$ and $(\eta, \theta) \in [0, \delta] \times [0, 2\pi]$. 
		Set $\Xi^{\theta}(\eta, z) := (\Phi^{\theta}_{\eta})_{*,z}$ and $W^{\theta}(\zeta) := \cos\theta \cdot U(\zeta) + \sin\theta \cdot V(\zeta)$. 
		Since $\partial_{\eta}\Phi^{\theta}(\eta,z) = W^{\theta}(\Phi^{\theta}(\eta, z))$, we have
		$$
		\frac{d}{d\eta}\Xi^{\theta}(\eta,z) = (\nabla_{\zeta}W^{\theta})(\Phi^{\theta}(\eta,z)) \cdot \Xi^{\theta}(\eta,z).
		$$
		Here, when we express $W^{\theta}(\zeta) = \sum_{i}W^{\theta}_{i}(\zeta)(\frac{\partial}{\partial x_{i}})_{\zeta}$ with $(x_{1},x_{2},x_{3},x_{4})$ 
		being the real coordinates of ${\mathbf C}^{2}$, $\nabla_{\zeta}W^{\theta}(\zeta)$ denotes the Jacobian matrix 
		$( \frac{\partial W^{\theta}_{i}}{\partial x_{j}}(\zeta) )$.
		\par
		Set $\psi^{\theta}(\eta) := \| (\Phi^{\theta}_{\eta})_{*,z} - (\Phi^{\theta}_{0})_{*,z} \| = \| \Xi^{\theta}(\eta, z) - \Xi^{\theta}(0,z) \|$. Then we get
		$$
		\psi^{\theta}(\eta) 
		= 
		\left\| \int_{0}^{\eta} \frac{d}{d\sigma}\Xi^{\theta}(\sigma,z) d\sigma \right\| 
		= 
		\left\| \int_{0}^{\eta} (\nabla_{\zeta}W^{\theta})(\Phi^{\theta}(\sigma,z)) \cdot \Xi^{\theta}(\sigma,z) d\sigma \right\| 
		\leq 
		N_{r} \int_{0}^{\eta} \| \Xi^{\theta}(\sigma,z) \| d\sigma,
		$$
		where we used $\| \nabla_{\zeta} W^{\theta} \| \leq N_{r}$ on $B(p,1)\setminus B(p,r)$ to get the last inequality. Hence we have
		$$
		\psi^{\theta}(\eta) \leq N_{r} \int_{0}^{\eta} \psi^{\theta}(\sigma)\,d\sigma + N_{r} \| (\Phi^{\theta}_{0})_{*,z} \| \cdot \eta 
		\leq N_{r} \int_{0}^{\eta} \psi^{\theta}(\sigma)\,d\sigma + N_{r}K \eta
		$$
		for all $\eta \in [0,\delta]$. By Gronwall's lemma, we get $\psi^{\theta}(\eta) \leq K(e^{N_{r} \eta} - 1)$. Since $0 \leq \eta \leq \delta \leq \frac{1}{2N_{r}}$, this implies 
		$\psi^{\theta}(\eta) \leq e^{\frac{1}{2}}KN_{r} \eta$. 
	\end{pf}

	\begin{lemma}
		\label{lemma:H1norm}
		There exists a constant $K_{2}>0$ such that 
		$$
		\left\| (\Phi^{\theta}_{\eta})^{*} g_{\eta e^{i\theta}} - (\Phi^{\theta}_{0})^{*} g_{0} \right\|_{X_{0}\setminus\bigcup_{p\in {\rm Sing}\,X_{0}}B(p,2r)} 
		\leq 
		K_{2} \eta^{\frac{1}{2}}
		$$
		for all $(\eta , \theta) \in [0, r^{8\nu}] \times [0, 2\pi]$. In particular, if $0 \leq \eta \ll1$, then 
		$$
		\left\| (\Phi^{\theta}_{\eta})^{*} g_{\eta e^{i\theta}} - (\Phi^{\theta}_{0})^{*} g_{0} \right\|_{X_{0}\setminus\bigcup_{p\in {\rm Sing}\,X_{0}}B(p,2\eta^{\frac{1}{8\nu}})} 
		\leq 
		K_{2} \eta^{\frac{1}{2}}.
		$$
	\end{lemma}
	
	\begin{pf}
		Since $g_{\eta e^{i\theta}}=g^{X}|_{X_{\eta e^{i\theta}}}$ and hence $(\Phi^{\theta}_{\eta e^{i\theta}})^{*}g_{\eta e^{i\theta}} = (\Phi_{\eta e^{i\theta}})^{*}g^{X}$, 
		the first inequality follows from \eqref{eqn:diffeo} and Lemma~\ref{lemma:diffeo}.
		The second inequality follows from the first one by setting $r = \eta^{\frac{1}{8\nu}}$. 
	\end{pf}

	For $\eta \in [0, \delta]$ and $\theta \in [0, 2\pi]$, set $\Psi_{\eta}^{\theta} := (\Phi_{\eta}^{\theta})^{-1}$. Then $\Psi_{\eta}^{\theta}$ is a diffeomorphism from 
	$\Phi_{\eta}^{\theta}(X_{0}\setminus \bigcup_{p\in{\rm Sing}\,X_{0}} B(p,2r))$ to $X_{0}\setminus \bigcup_{p\in{\rm Sing}\,X_{0}} B(p,2r)$.

	\begin{lemma}
		\label{lemma:est:L2}
		For any $\chi, \chi' \in C^{\infty}_{0}(X_{0}\setminus\bigcup_{p\in {\rm Sing}\,X_{0}}B(p, 2\eta^{\frac{1}{8\nu}}))$, the following inequalities hold:
		$$
		\left| \left( (\Psi^{\theta}_{\eta})^{*}\chi, (\Psi^{\theta}_{\eta})^{*}\chi' \right)_{L^{2}(X_{\eta e^{i\theta}})} - \left( \chi, \chi' \right)_{L^{2}(X_{0})} \right|
		\leq
		K_{3} \eta^{\frac{1}{2}} \| \chi \|_{L^{2}(X_{0})} \| \chi' \|_{L^{2}(X_{0})}, 
		\leqno{(1)}
		$$
		$$
		\left| \left\| d (\Psi_{\eta}^{\theta})^{*}\chi \right\|_{L^{2}(X_{\eta e^{i\theta}})}^{2} - \left\| d\chi \right\|_{L^{2}(X_{0})}^{2} \right|
		\leq 
		K_{3} \eta^{\frac{1}{2}} \| d \chi \|_{L^{2}(X_{0})}^{2},
		\leqno{(2)}
		$$
		where $K_{3}>0$ is a constant independent of $\chi$, $\chi'$ and $\eta$, $\theta$.
	\end{lemma}
	
	\begin{pf}
		By Lemma~\ref{lemma:H1norm}, there exists a constant $K'_{3}>0$ independent of $\eta$ and $\theta$ with
		$$
		\left\| \frac{(\Psi^{\theta}_{\eta})^{*}dv_{\eta e^{i\theta}}}{dv_{0}} - 1 \right\|_{X_{0}\setminus\bigcup_{p\in {\rm Sing}\,X_{0}}B(p, 2\eta^{\frac{1}{8\nu}})} 
		\leq 
		K'_{3} \eta^{\frac{1}{2}}.
		$$ 
		This, together with
		$$
		\left( (\Psi^{\theta}_{\eta})^{*}\chi, (\Psi^{\theta}_{\eta})^{*}\chi' \right)_{L^{2}(X_{\eta e^{i\theta}})} 
		=
		\int_{X_{0}\setminus\bigcup_{p\in {\rm Sing}\,X_{0}}B(p, 2\eta^{\frac{1}{8\nu}})} \chi(z)\overline{\chi'(z)} \, (\Phi^{\theta}_{\eta})^{*}dv_{\eta e^{i\theta}}
		$$
		and the Cauchy-Schwarz inequality, yields (1). 
		\par
		Let $*_{\eta}$ be the Hodge star operator with respect to $(\Phi^{\theta}_{\eta})^{*}g_{\eta e^{i\theta}}$ 
		acting on the $1$-forms on $X_{0}\setminus\bigcup_{p\in {\rm Sing}\,X_{0}}B(p, 2\eta^{\frac{1}{8\nu}})$. 
		By the second inequality of Lemma~\ref{lemma:H1norm}, there exists a constant $K''_{3}>0$ such that
		$\| *_{\eta} - *_{0} \|_{X_{0}\setminus\bigcup_{p\in {\rm Sing}\,X_{0}}B(p, 2\eta^{\frac{1}{8\nu}})} \leq K''_{3} \eta^{\frac{1}{2}}$.
		This, together with
		$$
		\left\| d (\Psi^{\theta}_{\eta})^{*}\chi \right\|_{L^{2}(X_{\eta e^{i\theta}})}^{2} 
		=
		\int_{X_{0}\setminus\bigcup_{p\in {\rm Sing}\,X_{0}}B(p, 2\eta^{\frac{1}{8\nu}})} d\chi(z) \wedge *_{\eta} (\overline{d\chi(z)}),
		$$
		yields (2). This completes the proof.
	\end{pf}

	Recall that $X_{0} = C_{1}+\cdots+C_{N}$ is the irreducible decomposition. 
	For $p \in {\rm Sing}\,X_{0} \cap C_{i}$, we fix a system of local coordinates $(U_{p},\zeta)$,  $\zeta = (\zeta_{1}, \zeta_{2})$, of $X$ centered at $p$. 
	On $U_{p}$, we define $r_{p}(z) := \| \zeta(z) \| = \sqrt{ |\zeta_{1}(z)|^{2} + |\zeta_{2}(z)|^{2}}$.

	\begin{lemma}
		\label{lemma:cut:off}
		For every $0<\epsilon \ll 1$, there exists $\chi_{\epsilon}^{(i)} \in C^{\infty}_{0}( C_{i} \setminus {\rm Sing}\,X_{0} )$ $(1 \leq i \leq N)$ with the following properties:
		\begin{itemize}
			\item[(1)]
			$0 \leq \chi^{(i)}_{\epsilon} \leq 1$. On $C_{i} \setminus \bigcup_{p \in {\rm Sing}\,X_{0} \cap C_{i}} U_{p}$, one has $\chi^{(i)}_{\epsilon} = 1$.
			\item[(2)]
			For any $p \in {\rm Sing}\,X_{0} \cap C_{i}$, one has $\chi_{\epsilon}^{(i)}(z) = 0$ if $r_{p}(z) \leq \frac{1}{2}\epsilon$ and 
			$\chi_{\epsilon}^{(i)}(z) = 1$ if $r_{p}(z) \geq 2\sqrt{\epsilon}$.
			\item[(3)]
			$\| d \chi^{(i)}_{\epsilon} \|_{L^{2}}^{2} \leq K/(\log \epsilon^{-1})$, where $K_{4}>0$ is a constant independent of $\epsilon$.
		\end{itemize}
	\end{lemma}
	
	\begin{pf}
		For $0< \epsilon \ll 1$, we define $B(p,\epsilon) := \{ z \in C_{i} ;\, r_{p}(z) < \epsilon \}$. We set
		$$
		\psi^{(i)}_{\epsilon}(z) 
		:=
		\begin{cases}
			\begin{array}{ll}
				0 & (z \in \overline{B(p,\epsilon)} ) 
				\\
				\displaystyle \frac{2}{\log \epsilon^{-1}} \int_{\epsilon}^{r_{p}(z)} \frac{d\rho}{\rho} & (z \in \overline{B(p,\sqrt{\epsilon})}\setminus B(p,\epsilon))
				\\
				1 & (z \in C_{i} \setminus \bigcup_{p \in C_{i}\cap {\rm Sing}\,X_{0}} B(p, \sqrt{\epsilon}) ).
			\end{array}
		\end{cases} 
		$$
		By definition, we have (1), (2) and ${\rm Supp}(d \psi^{(i)}_{\epsilon} ) \subset \overline{B(p, \sqrt{\epsilon})} \setminus B(p, \epsilon)$.
		Let ${\rm dist}_{C_{i}}(\cdot, \cdot)$ be the distance function on $C_{i}$ with respect to $g^{X}|_{C_{i}}$. 
		Since $|r_{p}(z) - r_{p}(w)| \leq {\rm dist}_{{\mathbf C}^{2}}( \zeta(z), \zeta(w) ) \leq {\rm dist}_{C_{i}}(z,w)$, $r_{p}(\cdot)$ is a Lipschitz function on $U_{p}$
		with Lipschitz constant $1$. 
		Hence we get
		$$
		\left| d \psi^{(i)}_{\epsilon}(z) \right|
		\leq
		\begin{cases}
			\begin{array}{ll}
				\displaystyle \frac{2}{\log \epsilon^{-1}} \frac{1}{r_{p}(z)} & (z \in \overline{B(p,\sqrt{\epsilon})}\setminus B(p,\epsilon),\,\, p\in C_{i}\cap {\rm Sing}\,X_{0})
				\\
				0 & (\hbox{otherwise}).
			\end{array}
		\end{cases}
		$$
		By \cite[Lemma 3.4]{Yoshikawa97}, there exists a constant $K_{4}>0$ independent of $0 < \epsilon\ll1$ such that 
		$$
		\int_{B(p,1)} \left| d \psi^{(i)}_{\epsilon} \right|^{2} dv_{C_{i}}  
		\leq
		\frac{4}{(\log \epsilon^{-1})^{2}} \int_{\epsilon \leq r_{p}(z) \leq \sqrt{\epsilon}} \frac{dv_{C_{i}} }{r_{p}(z)^{2}}
		\leq
		\frac{K_{4}}{\log \epsilon^{-1}}.
		$$
		By an argument using Friedrichs mollifier, we can find a function $\chi_{\epsilon}^{(i)} \in C^{\infty}_{0}( C_{i} \setminus {\rm Sing}\,X_{0} )$ with (1), (2), (3).
		This completes the proof.
	\end{pf}

	\subsection{\bf Proof of Proposition~\ref{prop:lb:small:ev}}
	\label{sec:6.2}
	\par
	By the mini-max principle, we have
	$$
	\lambda_{k}(s) 
	= 
	\min_{\substack{V \subset C^{\infty}(X_{s}) \\ \dim V =k}} \max_{\substack{\varphi \in V \\ \|\varphi\| = 1}} (\square \varphi, \varphi)_{L^{2}(X_{s})}
	=
	\min_{\substack{V \subset C^{\infty}(X_{s}) \\ \dim V =k}} \max_{\substack{\varphi \in V \\ \|\varphi\| = 1}} \| d \varphi \|_{L^{2}(X_{s})}^{2}.
	$$
	It suffices to prove that for all $s \in S^{o}$ with $0 <|s| \ll 1$, there exists an orthogonal system of functions 
	$\{ \varphi_{1}(s), \ldots, \varphi_{N}(s) \} \subset C^{\infty}(X_{s})$ with
	\begin{equation}
		\label{eqn:app:ONF:1}
		\| \varphi_{i}(s) \|_{L^{2}} = 1+o(1),
		\qquad
		\| d\varphi_{i}(s) \|_{L^{2}}^{2} \leq \frac{K}{\log (|s|^{-1}) },
	\end{equation}
	where $K>0$ is a constant independent of $s \in S^{o}$. Let $\nu \in {\mathbf N}$ be the same integer as in Lemma~\ref{lemma:est:L2}.
	Let $\chi_{\epsilon}^{(i)}$ be the function as in Lemma~\ref{lemma:cut:off}. Extending $\chi_{\epsilon}^{(i)}$ by zero on $C_{j}$ $(j\not=i)$,
	we regard $\chi_{\epsilon}^{(i)} \in C^{\infty}_{0}(X_{0} \setminus {\rm Sing}\,X_{0})$ with compact support in 
	$X_{0} \setminus \bigcup_{p\in {\rm Sing}\,X_{0}} B(p, \frac{1}{2}\epsilon)$.
	We set $\epsilon(s) := 2|s|^{\frac{1}{8\nu}}$. For $s=|s|e^{i\theta}$, we define
	$$
	\varphi_{i}(s) := (\Psi_{|s|}^{\theta})^{*} ( \chi_{\epsilon(s)}^{(i)} ) / \sqrt{{\rm Area}(C_{i})} \in C^{\infty}(X_{s}).
	$$
	Since ${\rm Supp}(\chi_{\epsilon(s)}^{(i)}) \cap {\rm Supp}(\chi_{\epsilon(s)}^{(j)}) = \emptyset$ for $i \not=j$, it is obvious that
	$\{ \varphi_{1}(s), \ldots, \varphi_{N}(s) \}$ is an orthogonal system of smooth functions on $X_{s}$.
	By Lemma~\ref{lemma:est:L2} (1) and Lemma~\ref{lemma:cut:off} (1), (2), we get
	\begin{equation}
		\label{eqn:app:ONF:2}
		\| \varphi_{i}(s) \|_{L^{2}(X_{s})}^{2} 
		= 
		\| \chi_{\epsilon(s)}^{(i)} \|_{L^{2}(C_{i})}^{2} / {\rm Area}(C_{i}) + O(|s|^{\frac{1}{8\nu}})
		=
		1 + O(|s|^{\frac{1}{8\nu}}).
	\end{equation}
	By Lemma~\ref{lemma:est:L2} (2) and Lemma~\ref{lemma:cut:off} (3), we get
	\begin{equation}
		\label{eqn:app:ONF:3}
		\| d \varphi_{i}(s) \|_{L^{2}(X_{s})}^{2}
		\leq 
		\frac{ \| d \chi_{\epsilon(s)}^{(i)} \|_{L^{2}(X_{0})}^{2} }{{\rm Area}(C_{i})} + C_{3} |s|^{\frac{1}{2}} \| d \chi_{\epsilon(s)}^{(i)} \|_{L^{2}(X_{0})}^{2} 
		\leq
		\frac{K'}{\log (|s|^{-1}) },
	\end{equation}
	where $K'>0$ is a constant independent of $s\in S^{o}$. We deduce \eqref{eqn:app:ONF:1} from \eqref{eqn:app:ONF:2}, \eqref{eqn:app:ONF:3}.
	This completes the proof of Proposition~\ref{prop:lb:small:ev}.
	\qed

	\section{Proof of Theorem~\ref{Main:Thm:1}}
	\label{sect:7}
	We keep the notation in Introduction.
	Since
	$$
	\lambda_{N-1}(s)^{N-1} \geq \prod_{i=1}^{N-1}  \lambda_{i}(s) = \frac{c + o(1)}{ (\log (|s|^{-1}) )^{N-1}}
	$$
	by Theorem~\ref{Main:Thm:2}, we get 
	\begin{equation}
		\label{eqn:ev:lb:N-1}
		\lambda_{N-1}(s) \geq \frac{c^{1/(N-1)} + o(1)}{ \log (|s|^{-1}) }.
	\end{equation}
	Combining \eqref{eqn:ev:lb:N-1} with Proposition~\ref{prop:lb:small:ev} for $i=N-1$, we get
	\begin{equation}
		\label{eqn:ev:N-1}
		\frac{c^{1/(N-1)}}{ \log (|s|^{-1}) } \leq \lambda_{N-1}(s) \leq \frac{K(N-1)}{ \log (|s|^{-1})}.
	\end{equation}
	This proves the assertion for $i=N-1$. By \eqref{eqn:ev:N-1} and Theorem~\ref{Main:Thm:2}, there exist constants $K',K''>0$ such that
	for all $s \in S^{o}$,
	\begin{equation}
		\label{eqn:ev:prod:N-2}
		\frac{K'}{ (\log (|s|^{-1}) )^{N-2}} \leq \prod_{i=1}^{N-2}  \lambda_{i}(s) \leq \frac{K''}{ (\log (|s|^{-1}) )^{N-2}}.
	\end{equation}
	Then we get
	\begin{equation}
		\label{eqn:ev:lb:N-2}
		\lambda_{N-2}(s)^{N-2} \geq \prod_{i=1}^{N-2}  \lambda_{i}(s) \geq \frac{K'}{ (\log (|s|^{-1}) )^{N-2}}.
	\end{equation}
	Namely, we have
	\begin{equation}
		\label{eqn:ev:lb:N-1:2}
		\lambda_{N-2}(s) \geq \frac{(K')^{1/(N-2)}}{ \log (|s|^{-1}) }.
	\end{equation}
	This, together with Proposition~\ref{prop:lb:small:ev} for $i=N-2$, yields the assertion for $i=N-2$. Inductively, we obtain the assertion for all $1\leq i \leq N-1$.
	This completes the proof.
	\qed

	\section{Examples}
	\label{sect:8}
	In this section, we discuss some illustrating examples concerning small eigenvalues of Laplacian for degenerating families of Riemann surfaces.

	\begin{example}
		\label{ex:Fermat}
		Let $d \in {\mathbf Z}_{>0}$. For $s \in {\mathbf C}$, we define a plane curve $X_{s} \subset {\mathbf P}^{2}$ by
		$$
		X_{s} = \{ (x:y:z) \in {\mathbf P}^{2} ;\, x^{d} + y^{d} + s z^{d} = 0 \}.
		$$
		Then $X_{s}$ $(s\not=0)$ is isomorphic to the Fermat curve $X_{1}$. When $d \geq 4$, since $X_{s}$ endowed with the hyperbolic metric $g_{s}^{\rm hyp}$
		of the Gauss curvature $-1$ is isometric to the hyperbolic curve $(X_{1}, g^{\rm hyp}_{1})$, the $k$-th eigenvalue $\lambda_{k}^{\rm hyp}(s)$ of the hyperbolic 
		Laplacian of $X_{s}$ is a constant function on ${\mathbf C}^{*}$: 
		\begin{equation}
			\label{eqn:e:v:Fremat:hyp}
			\lambda_{k}^{\rm hyp}(s) = \lambda_{k}^{\rm hyp}(1) 
			\qquad
			(s\not=0).
		\end{equation}
		On the other hand, let $g_{s} = g^{\rm FS}|_{X_{s}}$ be the restriction of the Fubini-Study metric
		of ${\mathbf P}^{2}$ to $X_{s}$ and let $\lambda_{k}(s)$ be the $k$-th eigenvalue of the Laplacian of $(X_{s}, g_{s})$. 
		Since $X_{0}$ is the union of $d$ lines of ${\mathbf P}^{2}$, the eigenvalues of the Laplacain of $(X_{0}, g_{0})$ are the $d$-copies of the eigenvalues 
		of the Laplacian of the round sphere $S^{2}$. 
		\par
		Let us see that the estimate for $\lambda_{1}(s)$ deduced from \eqref{eqn:e:v:Fremat:hyp} and Lemma~\ref{lemma:minimax} is of type \eqref{eqn:Gromov}.
		Suppose that $s \in {\mathbf R}_{>0}$ and define $\varphi_{s}(x : y : z) := (x : y : s^{1/d}z)$. Since $\varphi_{s} \in {\rm Aut}({\mathbf P}^{2})$ is 
		such that $\varphi_{s}(X_{s}) = X_{1}$, we have $g_{s}^{\rm hyp} = \varphi_{s}^{*} g_{1}^{\rm hyp}$.
		Since there are constants $K_{0}, K_{1}>0$ with $K_{0}g_{1} \leq g_{1}^{\rm hyp} \leq K_{1} g_{1}$, we have
		$K_{0}\varphi_{s}^{*}g_{1} \leq g_{s}^{\rm hyp} \leq K_{1} \varphi_{s}^{*}g_{1}$. By Lemma~\ref{lemma:minimax}, we get
		\begin{equation}
			\label{eqn:mini:max}
			\frac{\lambda_{1}(s)}{\lambda_{1}^{\rm hyp}} 
			\geq 
			\min_{X_{s}} \frac{g_{s}^{\rm hyp}}{g_{s}} 
			\geq 
			K_{0} \min_{X_{s}} \frac{\varphi_{s}^{*}g_{1}}{g_{s}}
			=
			K_{0} \min_{X_{s}} \frac{\varphi_{s}^{*}g^{\rm FS}|_{X_{s}}}{g^{\rm FS}|_{X_{s}}}
			\geq
			K_{0} \min_{{\mathbf P}^{2}} \frac{\varphi_{s}^{*}g^{\rm FS}}{g^{\rm FS}}.
		\end{equation}
		Since the last term of \eqref{eqn:mini:max} is bounded from below by $C|s|^{\alpha}$ for some positive constants $C$, $\alpha$, 
		\eqref{eqn:mini:max} yields an estimate for $\lambda_{1}(s)$ of type \eqref{eqn:Gromov}. 
		In this example, it seems difficult to obtain the genuine behavior of $\lambda_{1}(s)$ as in Theorem~\ref{Main:Thm:1} 
		by means of the mini-max principle like Lemma~\ref{lemma:minimax}.
	\end{example}

	\begin{example}
		\label{ex:Bismut:Bost}
		Let $f \colon X \to S$ be a degeneration of Riemann surfaces of genus $g>1$ such that $X_{0}$ is a stable curve with $N>1$ irreducible components. 
		Hence the singularities of $X_{0}$ consist of ordinary double points. Following Bismut-Bost \cite{BismutBost90}, we fix a Hermitian metric 
		on the relative canonical bundle of $f$. 
		Namely, let $h_{X/S}$ be a smooth Hermitian metric on the relative canonical bundle $K_{X/S} = K_{X} \otimes f^{*}K_{S}^{-1}$. 
		Then $h_{X/S}$ induces a Hermitian metric on $TX/S|_{X\setminus{\rm Sing}\,X_{0}}$, the relative tangent bundle restricted to the regular locus of $f$.
		This Hermitian metric on $TX/S|_{X\setminus{\rm Sing}\,X_{0}}$ is still denoted by $h_{X/S}$. Let $p \in {\rm Sing}\,X_{0}$ be an arbitrary  singular point of $X_{0}$. 
		We have a local coordinates $(U_{p},(z,w))$ of $X$ centered at $p$ such that $f(z,w) = zw$ on $U_{p}$. 
		Since $K_{X/S}|_{U_{p}} \cong {\mathcal O}_{U_{p}}\cdot(dz/z) = {\mathcal O}_{U_{p}}\cdot(dw/w)$ in the canonical way, under this isomorphism, 
		there exists a positive smooth function $a_{p}(x,y)>0$ defined on $U_{p}$ such that $h_{X/S}(dz/z, dz/z) = h_{X/S}(dw/w, dw/w) = a_{p}(x,y)$ on $U_{p}$. 
		We set $h_{s} := h_{X/S}|_{X_{s}}$. Then each connected component of $U_{p}\cap(X_{0}\setminus\{p\})$ endowed with $h_{0}$ is quasi-isometric to the cylinder
		$S^{1}\times(0,\infty)$ endowed with the metric $d\theta^{2} + dr^{2}$. Contrary to the case of hyperbolic metric or the metric induced from $X$, 
		even though the Riemann surfaces $X_{s}$ are pinched along closed simple curves, the length of the corresponding geodesics is uniformly positive in this case.
		\par
		Let $g^{\rm hyp}_{s}$ be the hyperbolic metric on $X_{s}$ with constant Gauss curvature $-1$. 
		Then there exist constants $C_{0},C_{1}>0$ such that on $U_{p}\cap X_{s}$, $p \in {\rm Sing}\,X_{0}$,
		\begin{equation}
			\label{eqn:cylinder}
			C_{0}\frac{dz d\bar{z}}{|z|^{2}} \leq h_{s}|_{U_{p}\cap X_{s}} \leq C_{1}\frac{dz d\bar{z}}{|z|^{2}},
		\end{equation}
		\begin{equation}
			\label{eqn:hyp}
			C_{0}\frac{dz d\bar{z}}{|z|^{2}(\log |z|)^{2}} \leq g_{s}^{\rm hyp}|_{U_{p}\cap X_{s}} \leq C_{1}\frac{dz d\bar{z}}{|z|^{2}(\log |z|)^{2}}.
		\end{equation}
		\par
		Let $\lambda_{k}^{\rm cyl}(s)>0$ be the $k$-th nonzero eigenvalue of the Laplacian of $X_{s}$ with respect to $h_{s}$. 
		Let $\lambda_{1}^{\rm hyp}(s)>0$ be the first nonzero eigenvalue of the hyperbolic Laplacian of $X_{s}$. 
		By \eqref{eqn:cylinder}, \eqref{eqn:hyp} and Lemma~\ref{lemma:minimax}, we get
		$$
		\frac{\lambda_{1}^{\rm cyl}(s)}{\lambda_{1}^{\rm hyp}(s)} 
		\geq 
		\min_{X_{s}} \frac{g_{s}^{\rm hyp}}{h_{s}} 
		\geq 
		\frac{C_{0}}{C_{1}} \min_{p\in{\rm Sing}\,X_{0},\,z\in U_{p}\cap X_{s}} \frac{1}{(\log |z|)^{2}} = \frac{C_{0}}{4C_{1}} \frac{1}{(\log (|s|^{-1}) )^{2}}.
		$$
		This, together with Theorem~\ref{eqn:small:ev:stable:c}, yields the following lower bound
		\begin{equation}
			\label{eqn:lower:bd:cyl}
			\lambda_{1}^{\rm cyl}(s) \geq \frac{K_{0}}{(\log (|s|^{-1}) )^{3}},
		\end{equation}
		where $K_{0}>0$ is a constant independent of $s \in \varDelta^{*}$. To obtain an upper bound for $\lambda_{k}^{\rm cyl}(s)$, consider the orthogonal system of
		smooth functions $\{ \varphi_{1}(s), \ldots, \varphi_{N}(s) \} \subset C^{\infty}(X_{s})$ constructed in Section~\ref{sec:6.2}. 
		With respect to the metric $h_{s}$ on $X_{s}$, we have
		$$
		\| \varphi_{k}(s) \|_{L^{2}, h_{s}}^{2} = 1 + O\left( |s|^{\frac{1}{8\nu}}\log ( |s|^{-1}) \right).
		$$
		Since ${\rm Area}( {\rm supp}\,\varphi_{k}(s) ) = O\left( \log (|s|^{-1}) \right)$, in the same way as in \eqref{eqn:app:ONF:3}, we get
		$$
		\| d\varphi_{k}(s) \|_{L^{2}, h_{s}}^{2} = O\left( ( \log (|s|^{-1}) )^{2} \right).
		$$
		By the mini-max principle, there exists a positive constant $K_{1}>0$ such that for all $s \in \varDelta^{*}$ and $1 \leq k \leq N-1$,
		\begin{equation}
			\label{eqn:upper:bd:cyl}
			\lambda_{k}^{\rm cyl}(s) \leq \frac{K_{1}}{(\log (|s|^{-1}) )^{2}}.
		\end{equation}
		By \eqref{eqn:lower:bd:cyl}, \eqref{eqn:upper:bd:cyl}, we conclude the following for the asymptotic behavior of the first $N-1$ eigenvalues of $(X_{s}, h_{s})$:
		\begin{equation}
			\label{eqn:asymp:bd:cyl}
			\frac{K_{0}}{(\log (|s|^{-1}) )^{3}} \leq \lambda_{1}^{\rm cyl}(s) \leq \cdots \leq \lambda_{N-1}^{\rm cyl}(s) \leq \frac{K_{1}}{(\log (|s|^{-1}) )^{2}}.
		\end{equation}
		Comparing \eqref{eqn:asymp:bd:cyl} with \eqref{eqn:small:ev:stable:c} and Theorem~\ref{Main:Thm:1} , 
		we infer that the behavior of the first $N-1$ eigenvalues of $(X_{s}, h_{s})$ differs from those of $(X_{s}, g^{\rm hyp}_{s})$ or $(X_{s}, g^{\rm ind}_{s})$,
		where $g^{\rm ind}_{s}$ is the K\"ahler metric on $X_{s}$ induced from the K\"ahler metric on the ambient space $X$.
	\end{example}

	\begin{remark}
		In Example~\ref{ex:Bismut:Bost}, since the area of $(X_{s}, h_{s})$ grows like ${\rm Const.}\log (|s|^{-1})$ as $s\to 0$ and the length of simple curves of $X_{s}$ 
		converging to $p \in {\rm Sing}\,X_{0}$ is uniformly bounded from below by a positive constant, it is very likely that the Cheeger constant 
		$h(X_{s})$ of $(X_{s}, h_{s})$ satisfies the inequality $h(X_{s}) \geq C/\log(|s|^{-1})$, where $C>0$ is a positive constant independent of $s$. 
		Assuming this estimate for $h(X_{s})$, we would deduce from \cite{Cheeger71} the following better estimate for the small eigenvalues of $(X_{s}, h_{s})$
		\begin{equation}
			\label{eqn:asymp:bd:cyl:2}
			\frac{K'_{0}}{(\log (|s|^{-1}) )^{2}} \leq \lambda_{1}^{\rm cyl}(s) \leq \cdots \leq \lambda_{N-1}^{\rm cyl}(s) \leq \frac{K'_{1}}{(\log ( |s|^{-1}) )^{2}},
		\end{equation}
		where $K'_{0}$, $K'_{1}$ are constants independent of $s \in S^{*}$.
	\end{remark}

	\section{Problems and conjectures}
	\label{sect:9}
	In this section, we propose some problems and conjectures related to the main results of this paper.

	\begin{problem}
		In Theorems~\ref{Main:Thm:1} and \ref{Main:Thm:2}, we assume that $X_{0}$ is reduced. Namely, $f$ has only isolated critical points.
		If $X_{0}$ is not reduced or equivalently $f$ has non-isolated critical points, do the statements of Theorems~\ref{Main:Thm:1} and \ref{Main:Thm:2} remain valid?
		It is also interesting to ask if these theorems remain valid when the total space $X$ admits singularities.
	\end{problem}

	\begin{problem}
		If $\lambda_{k}(s)$ is a small eigenvalue, then does the limit $\lim_{s\to 0} \lambda_{k}(s)\log (|s|^{-1})$ exist?
		When $X_{0}$ consists of two irreducible components, we have an affirmative answer. How about the case when $X_{0}$ consists of more than two components?
		If the answer is affirmative, can one give a geometric expression of the limit? When ${\rm Sing}\,X_{0}$ consists of a unique node, Ji-Wentworth give a conjectural
		expression of the limit \cite[Remark~5.10]{JiWentworth92}. 
	\end{problem}

	\begin{problem}
		Theorem~\ref{Main:Thm:1} gives an exact magnitude of the speed of convergence of the $k$-th eigenvalue function on $S$ for $k<N$.
		When $k \geq N$, does the estimate
		$$
		| \lambda_{k}(s) - \lambda_{k}(0) | \leq \frac{C_{k}}{ \log (|s|^{-1}) }
		$$
		hold for $0<|s|\ll 1$? Here $C_{k}>0$ is a constant.
	\end{problem}

	\begin{problem}
		Assume that $X_{s}$ is a hyperbolic Riemann surface endowed with the hyperbolic metric and $X_{0}$ is a stable curve.
		In \cite{Burger90}, Burger constructed a metric graph structure on the dual graph of $X_{0}$ by making use of certain geometric data of $X_{s}$ 
		such as the length of short geodesics and proved that the small eigenvalue $\lambda_{k}(s)$ is asymptotic to the $k$-th eigenvalue of 
		the Laplacian of this metric graph. Does the theorem of Burger hold true in the situation of Theorem~\ref{Main:Thm:1}?
		In general, can one construct a finite metric graph depending on the geometry of $X_{s}$ and $X_{0}$ whose eigenvalues are asymptotic to 
		the small eigenvalues $\lambda_{k}(s)$ of $X_{s}$? 
	\end{problem}

	\begin{conjecture}
		\label{conj:high:dim}
		It is natural to seek for a generalization of Theorems~\ref{Main:Thm:1} and \ref{Main:Thm:2} in higher dimensions.
		Let $f\colon X\to S$ be a one-parameter degeneration of compact K\"ahler manifolds of dimension $n$ 
		such that $f$ has only isolated critical points.
		Let $\{0=\cdots=0<\lambda_{1}(s) \leq \lambda_{2}(s) \leq\cdots\leq \lambda_{k}(s) \leq \cdots \}$ be the eigenvalues
		of the Laplacian $\square^{n,0}_{s}$ acting on $(n,0)$-forms on $X_{s}=f^{-1}(s)$ with respect to the metric induced from the K\"ahler metric on $X$. 
		Then we conjecture the following (1)-(4):
		\begin{enumerate}
			\item
			For all $k\in{\mathbf N}$, $\lambda_{k}(s)$ extends to a continuous function on $S$.
			\item
			Set $N:= \dim \ker \square^{n-1,0}_{0} - \dim \square^{n-1,0}_{s}$ $(s\not=0)$, where $\square^{n-1,0}_{0}$ is the Friedrichs extension of the Hodge-Kodaira
			Laplacian acting on the smooth $(n-1,0)$-forms with compact support on $X_{0}\setminus{\rm Sing}\,X_{0}$. Then $N<\infty$. 
			Moreover, $\lim_{s\to0}\lambda_{k}(s) =0$ for $k \leq N$ and $\lim_{s\to0}\lambda_{k}(s) >0$ for $k>N$.
			\item
			There exist constants $\nu\in{\mathbf N}$, $c\in{\mathbf R}_{>0}$ such that
			$$
			\prod_{k=1}^{N} \lambda_{k}(s) = \frac{c+o(1)}{(\log (|s|^{-1}) )^{\nu}}
			\qquad
			(s\to0).
			$$
			\item
			$\nu=N$. For $1\leq k \leq N$, there exist constants $K_{k},K'_{k}>0$ such that
			$$
			\frac{K_{k}}{\log (|s|^{-1}) } \leq \lambda_{k}(s) \leq \frac{K'_{k}}{\log (|s|^{-1}) }.
			$$ 
		\end{enumerate}
	\end{conjecture}

	\begin{conjecture}
		In the situation of Conjecture~\ref{conj:high:dim}, we conjecture that an analogue of \eqref{eqn:McKean:Singer:mult} holds and yields Conjecture~\ref{conj:high:dim} (3). 
		Let $H$ be an ample line bundle on $X$ endowed with a Hermitian metric of semi-positive curvature. 
		Let $\widetilde{f} \colon Y \to T$ be a semi-stable reduction of $f \colon X \to S$ associated to a ramified covering $\mu \colon T \to S$.
		Let $F \colon Y \to X$ be the holomorphic map of the total spaces. 
		Let $\{ \varphi_{1}, \ldots, \varphi_{m_{1}} \}$, $m_{1}=h^{0}(K_{Y_{t}}(F^{*}H))$ be a basis of $\widetilde{f}_{*}K_{Y/T}(F^{*}H)$ as a free ${\mathcal O}_{T}$-module.
		Similarly, let $\{ \omega_{1}, \ldots, \omega_{m_{2}} \}$, $m_{2} = h^{n,0}(Y_{t})$ be a basis of $\widetilde{f}_{*}K_{Y/T}$ as a free ${\mathcal O}_{T}$-module.
		Then the following identity holds in $C^{*}(T^{o})/C^{*}(T)$:
		$$
		\prod_{k=1}^{N} \lambda_{k}(\mu(t))^{-1}
		\equiv
		\frac{ \tau(X_{\mu(t)},K_{X_{\mu(t)}}) }{ \tau(X_{\mu(t)},K_{X_{\mu(t)}}(H_{\mu(t)})) } 
		\equiv
		\frac{\displaystyle \det\left( \int_{Y_{t}} F^{*}h_{H}(\varphi_{\alpha}(t) \wedge \overline{\varphi_{\beta}(t)} ) \right) }
		{\displaystyle \det\left( \int_{Y_{t}} \omega_{i}(t) \wedge \overline{\omega_{j}(t)} \right) }.
		$$
	\end{conjecture}

	\begin{problem}
		Can one extend the Schoen-Wolpert-Yau theorem \cite{SWY80} or Burger's theorem \cite{Burger90} in higher dimensions? 
		Consider the situation of Conjecture~\ref{conj:high:dim} and suppose that $X_{s}$ is endowed with a K\"ahler-Einstein metric of negative scalar curvature. 
		Then what can one say about the asymptotic behavior of the eigenvalues of the Laplacian $\square^{0,q}_{s}$ (or more general $\square^{p,q}_{s}$)? 
		Possibly, the answer will heavily depend on how bad the singularity of $X_{0}$ is. We conjecture that, as $s \to 0$, $\square^{0,n}_{s}$ has small eigenvalues 
		only when $X_{0}$ has non-canonical singularities. If $X_{0}$ has non-canonical singularities and if this conjecture is true, 
		by replacing the length of short geodesics with the volumes of the vanishing cycles of $X_{s}$, does the analogue of the Schoen-Wolpert-Yau theorem \cite{SWY80} 
		hold for the small eigenvalues of $\square^{0,n}_{s}$? 
	\end{problem}

	\begin{conjecture}
		For degenerations of hyperbolic Riemann surfaces to stable curves, the asymptotic behavior of the product of the small eigenvalues of
		the Laplacian was determined by Grotowski-Huntley-Jorgenson \cite[Th.\,1, Cor.\,2]{GHJ01} in terms of the length of short geodesics. 
		It is natural to seek for its counterpart in the following setting. Let $B$ be a polydisc of dimension $m \geq 2$. 
		Let $X$ be a complex manifold of dimension $m+1$ endowed with a positive line bundle $H$. 
		Let $f \colon X \to B$ be a proper surjective holomorphic map of relative dimension one with connected fibers. 
		Let $\Sigma = \Sigma_{f}$ be the critical locus of $f$ and let $\Delta = f(\Sigma)$ be the discriminant locus of $f$. 
		We assume that $0\in\Delta$, that $f$ is flat, and that $f$ induces a finite map from $\Sigma_{f}$ to $\Delta$. 
		We set $B^{o} := B \setminus \Delta$, $X^{o} := X \setminus f^{-1}( \Delta)$, and $f^{o} := f|_{X^{o}}$. 
		Then $f^{o} \colon X^{o} \to B^{o}$ is a family of compact Riemann surfaces.
		We set $X_{b} := f^{-1}(b)$ for $b \in B$. Then for $b \in \Delta$, $X_{b}$ is a singular projective curve with reduced structure.  
		Let $N$ be the number of the irreducible components of $X_{0}$. We assume that the line bundle $H$ is defined by an effective divisor on $X$ whose support is disjoint from ${\rm Sing}\,X_{0}$ and intersects each irreducible component of $X_{0}$ transversally at a single point.
		In particular, $\deg H|_{X_{b}}=N$ for $b \in B$.
		Fix a K\"ahler metric $h^{X}$ on $X$. By shrinking $B$ if necessary, since $\Sigma\cap X_{0}$ consists of isolated points,
		we can construct a Hermitian metric $h^{H}$ on $H$ with semi-positive curvature and vanishing Chern form near $\Sigma$ in the same way as in Lemma~\ref{lemma:HM}. For $b \in B^{o}$, let $\tau(X_{b}, K_{X_{b}})$ (resp. $\tau(X_{b}, K_{X_{b}}(H_{b}))$) be the analytic torsion of $K_{X_{b}}$ with respect to $h^{X}|_{X_{b}}$ (resp. $K_{X_{b}}(H_{b})$ with respect to $h^{X}|_{X_{b}}$, 
		$h^{H}|_{X_{b}}$).
		Let $0<\lambda_{1}(b) \leq \lambda_{2}(b) \leq \cdots$ be the eigenvalues of the Laplacian of $(X_{b}, h^{X}|_{X_{b}})$.
		Let $\omega_{1}, \ldots, \omega_{g+N-1}$ be a free basis of $f_{*}K_{X/B}(H)$ near $0\in B$ such that 
		$\omega_{1}, \ldots, \omega_{g}$ is a free basis of $f_{*}K_{X/B}$ near $0 \in B$. 
		By considering a semi-stable reduction of $f \colon X \to B$, we conjecture that the following generalizations of Theorems~\ref{thm:at:bd1} and \ref{thm:ratio:torsion:1} hold:
		\begin{itemize}
			\item[(1)] 
			In $C^{*}(B^{o})/C^{*}(B)$, one has
			$$
			\prod_{0< \lambda_{k}(b)<1} \lambda_{k}(b)^{-1} \equiv \frac{\tau(X_{b}, K_{X_{b}})}{\tau(X_{b}, K_{X_{b}}(H_{b}))}.
			$$
			\item[(2)]
			There exists a locally bounded function $\psi$ on $B$ such that on $B^{o}$, one has
			$$
			\log \frac{ \tau(X_{b},K_{X_{b}}) }{ \tau(X_{b}, K_{X_b}(H_{b})) } 
			=
			\log 
			\left[
			\frac{ \det\left(  \int_{X_{b}} h^{H}(\omega_{i}(b) \wedge \overline{\omega_{j}(b)}) \right)_{1\leq i,j \leq g+N-1} }
			{ \det\left(  \int_{X_{b}} \omega_{i}(b) \wedge \overline{\omega_{j}(b)} \right)_{1\leq i,j \leq g} }
			\right]
			+
			\psi(b).
			$$
			Moreover, there exists an alteration $\mu \colon B' \to B$ satisfying the following two conditions: (i) $f \colon X \to B$ admits a semi-stable reduction over $B'$, and (ii) $\mu^{*}\psi$ is a continuous function on $B'$. 
		\end{itemize}
		Since $\mu \colon B' \to B$ can contain some exceptional divisors in general, it seems unlikely that one can take $\psi \in C^{0}(B)$
		except for the case where $f \colon X \to B$ is already semi-stable.
	\end{conjecture}

	\section{Appendix}
	\label{sect:10}
	We keep the notation in Introduction. Then $X_{s} = f^{-1}(s)$ is endowed with the K\"ahler metric $g_{s} = g_{X}|_{X_{s}}$ induced from the K\"ahler metric 
	on the total space $X$. Let $K_{s}$ be the Gauss curvature of $(X_{s},g_{s})$. 
	In contrast to the hyperbolic metrics, we have the following:

	\begin{lemma}
		The minimum of $K_{s}$ diverges to $-\infty$ as $s \to 0$.
	\end{lemma}

	\begin{pf}
		Let us prove the assertion by contradiction. 
		By the Gauss-Codazzi equation, $K_{s}$ is uniformly bounded from above. Suppose that $\min_{X_{s}} K_{s}$ is bounded from below as $s \to 0$. 
		Then there exist constants $C_{0}, C_{1}$ such that $C_{0} \leq K_{s} \leq C_{1}$ for all $s\in S^{o}$. Let $B( p, r)$ be the open metric ball of radius $r>0$ 
		centered at $p\in {\rm Sing}\,X_{0}$. Then ${\rm Area}(X_{s}\cap B(p,r)) \leq C_{2}r^{2}$ for all $0<r\ll1$ and $s \in S^{o}$ sufficiently close to $0$ 
		with some constant $C_{2}>0$. Let $dv_{s}$ be the volume form of $(X_{s},g_{s})$. Since $g_{s}$ converges to $g_{0}$ on every compact subset of 
		$X_{0}\setminus{\rm Sing}\,X_{0}$, the assumption $C_{0} \leq K_{s} \leq C_{1}$ implies that 
		$$
		\int_{X_{0}\setminus\bigcup_{p\in{\rm Sing}\,X_{0}}B(p,\epsilon)} K_{0} dv_{0} 
		= \int_{X_{s}\setminus\bigcup_{p\in{\rm Sing}\,X_{0}}B(p,\epsilon)} K_{s} dv_{s} + O(\epsilon^{2}) = 2\pi\chi(X_{s})+O(\epsilon^{2}).
		$$ 
		Hence 
		\begin{equation}
			\label{eqn:GB}
			\int_{X_{0}\setminus{\rm Sing}\,X_{0}}K_{0}dv_{0} := \lim_{\epsilon\to0} \int_{X_{0}\setminus\bigcup_{p\in{\rm Sing}\,X_{0}}B(p,\epsilon)} K_{0} dv_{0} = 2\pi\chi(X_{s}).
		\end{equation}
		\par
		Let $\nu \colon \widetilde{X}_{0} \to X_{0}$ be the normalization. 
		By \cite[(4.12)]{BruningLesch96}, for every $q\in \nu^{-1}({\rm Sing}\,X_{0})$, there exists a positive integer $N_{q} \in {\mathbf Z}_{>0}$ such that
		\begin{equation}
			\label{eqn:BL}
			\frac{1}{2\pi}\int_{X_{0}\setminus{\rm Sing}\,X_{0}}K_{0}dv_{0} = \chi(\widetilde{X}_{0}) + \sum_{q\in \nu^{-1}({\rm Sing}\,X_{0})} (N_{q}-1).
		\end{equation}
		Let $a(X_{0}) = h^{1}({\mathcal O}_{X_{0}})$ be the arithmetic genus of $X_{0}$. By \cite[Chap.\,II, Sect.\,11]{BHPV04}, we have
		$a(X_{0}) = g(\widetilde{X}_{0}) + \sum_{p \in {\rm Sing}\,X_{0}} \delta_{p}$ with 
		$\delta_{p} = \dim_{\mathbf C} (\nu_{*}{\mathcal O}_{\widetilde{X}_{0}} / {\mathcal O}_{X_{0}})_{p} \geq 1$. 
		Since $s\mapsto h^{1}({\mathcal O}_{X_{s}})$ is a constant function on $S$ and hence $a(X_{0})=g(X_{s})$ ($s\not=0$), we get 
		\begin{equation}
			\label{eqn:Euler:char}
			\chi(\widetilde{X}_{0}) = 2( 1 - g(\widetilde{X}_{0}) ) = \chi(X_{s}) + 2\sum_{p \in {\rm Sing}\,X_{0}} \delta_{p}.
		\end{equation}
		By \eqref{eqn:BL}, \eqref{eqn:Euler:char},
		$$
		\frac{1}{2\pi}\int_{X_{0}\setminus{\rm Sing}\,X_{0}}K_{0}dv_{0} = \chi(X_{s}) + \sum_{p \in {\rm Sing}\,X_{0}} \{ 2\delta_{p} + \sum_{q\in{\nu}^{-1}(p)} (N_{q}-1) \}.
		$$
		Since $2\delta_{p} + \sum_{q\in{\nu}^{-1}(p)} (N_{q}-1) >0$ for $p\in {\rm Sing}\,X_{0}$, this contradicts \eqref{eqn:GB}.
	\end{pf}


\end{document}